\documentclass[journal,twocolumn]{IEEEtran}

\usepackage{soul}
\usepackage{epic,eepic,units}
\usepackage{psfrag}
\usepackage{latexsym}
\usepackage{longtable}
\usepackage{mathrsfs}
\usepackage{bigstrut}
\usepackage{graphicx}
\usepackage{pifont}
\usepackage{amsmath}
\usepackage{setspace}
\usepackage{subfigure}
\usepackage{mathtools}
\usepackage{mathcomp}
\usepackage{mathdots}
\usepackage{rotating}
\usepackage{array}
\usepackage{multirow}
\usepackage{multicol}
\usepackage{color}
\usepackage{url}
\usepackage{algorithm}
\usepackage{algorithmic}

\usepackage{cite,amssymb}
\usepackage[dvipsnames]{xcolor}

\usepackage{euscript,yfonts,dsfont,bbm,amstext,wasysym,pdfsync,amsfonts,arydshln,framed}

\usepackage{tikz}
\usetikzlibrary{arrows,shapes,trees,calc,positioning}
\usetikzlibrary{matrix}

\newtheorem{remark}{Remark}
\newtheorem{problem}{Problem}

\newcommand{\enma}[1]   {\ensuremath{#1}}

\newcommand{\beq}{\begin{equation}}
\newcommand{\eeq}{\end{equation}}
\newcommand{\bseq}{\begin{subequations}}
\newcommand{\eseq}{\end{subequations}}
\newcommand{\beqn}{\begin{eqnarray}}
\newcommand{\eeqn}{\end{eqnarray}}
\newcommand{\ba}{\begin{array}}
\newcommand{\ea}{\end{array}}
\newcommand{\bct}{\begin{center}}
\newcommand{\ect}{\end{center}}
\newcommand{\btmz}{\begin{itemize}}
\newcommand{\etmz}{\end{itemize}}
\newcommand{\benum}{\begin{enumerate}}
\newcommand{\eenum}{\end{enumerate}}

\newcommand{\cH}{\enma{\mathcal H}}
\newcommand{\cL}{\enma{\mathcal L}}

\newcommand{\norm}[1]{\| #1 \|}                 

\newcommand{\Ht}{ \mathcal{H}_{2}}

\newcommand{\trace}     {\enma{\mathrm{trace}}}

\newcommand{\inner}[2]{\left\langle #1,#2 \right\rangle}

\newcommand{\bv}{{\bf v}}

\newcommand{\matbegin}{
        \left[
}
\newcommand{\matend}{
        \right]
}

\newcommand{\tbo}[2]{
  \matbegin \begin{array}{c}
       #1 \\ #2
       \end{array} \matend }

\newcommand{\tbt}[4]{
  \matbegin \begin{array}{cc}
       #1 & #2 \\ #3 & #4
       \end{array} \matend }

\newcommand{\be}{\begin{equation}}
\newcommand{\ee}{\end{equation}}

\newcommand{\cplxs}{ C\kern -.35em \rule{0.03 em}{.7 ex}~   }

\def\complex{\hbox{C\kern -.45em \rule{0.03 em}{1.5 ex}}~}

\newcommand{\bi}{\begin{itemize}}
\newcommand{\ei}{\end{itemize}}

\newcommand{\DefinedAs}[0]{\mathrel{\mathop:}=}

\DeclareMathOperator*{\argmin}{argmin}

\DeclareMathOperator*{\minimize}{minimize}

\DeclareMathOperator*{\subject}{subject~to}

\newcommand{\vsp}{\vspace*{0.15cm}}

\newtheorem{mythm}{Theorem}
\newtheorem{myprop}{Proposition}
\newtheorem{mylem}{Lemma}

\newtheorem{myass}{Assumption}

\newcommand{\qedsymbol}{$\blacksquare$}

\newcommand{\cA}{{\cal A}}
\newcommand{\cl}{{\cal L}}
\newcommand{\cC}{{\cal C}}
\newcommand{\cD}{{\cal D}}
\newcommand{\cE}{{\cal E}}

\newcommand{\cI}{{\cal I}}
\newcommand{\cK}{{\cal K}}

\newcommand{\cB}{{\cal B}}
\newcommand{\cS}{{\cal S}}
\newcommand{\cM}{{\cal M}}

\newcommand{\hx}{{\hat x}}

\newcommand{\bbR}{\mathbb{R}}
\newcommand{\bbC}{\mathbb{C}}

\newcommand{\eps}{{\epsilon}}

\newcommand{\non}{\nonumber}
\newcommand{\ds}{\displaystyle}

\newcommand{\mrd}{\mathrm{d}}
\newcommand{\mre}{\mathrm{e}}

\newcommand{\bu}{{\bf u}}

\newcommand{\bpsi}{\mbox{\boldmath$\psi$}}

\newcommand{\bk}{{\bf k}}

\newcommand{\prox}{\mathbf{prox}}
\newcommand{\tY}{\tilde{Y}}
\newcommand{\tX}{\tilde{X}}
\newcommand{\La}{L_a}

\newcommand{\mua}{\mu_a}

\IEEEoverridecommandlockouts                              
\begin{document}

\title{
Proximal algorithms for large-scale\\ 
statistical modeling and sensor/actuator selection
}

\author{Armin Zare, {\em Member, IEEE\/}, Hesameddin Mohammadi, {\em Student Member, IEEE\/}, Neil K.\ Dhingra, {\em Member, IEEE\/}, Tryphon T.\ Georgiou, {\em Fellow, IEEE\/}, and Mihailo R.\ Jovanovi\'c, {\em Fellow, IEEE\/}
\thanks{Financial support from the National Science Foundation under Awards CMMI 1739243, ECCS 1509387, 1708906, 1809833, and 1839441, and the Air Force Office of Scientific Research under FA9550-16-1-0009, FA9550-17-1-0435, and FA9550-18-1-0422 is gratefully acknowledged.
}
\thanks{{A.\ Zare is with the Department of Mechanical Engineering, University of Texas at Dallas, Richardson, TX 75219.} H.\ Mohammadi and M.\ R.\ Jovanovi\'c are with the Ming Hsieh Department of Electrical and Computer Engineering, University of Southern California, Los Angeles, CA 90089. N.\ K.\ Dhingra is with Numerica Corporation, Fort Collins, CO 80528. T.\ T.\ Georgiou is with the Department of Mechanical and Aerospace Engineering, University of California, Irvine, CA 92697. E-mails: {armin.zare@utdallas.edu}, hesamedm@usc.edu, neil.k.dh@gmail.com, tryphon@uci.edu, mihailo@usc.edu.}
}

\maketitle
	\begin{abstract}
Several problems in modeling and control of stochastically-driven dynamical systems can be cast as regularized semi-definite programs. We examine two such representative problems and show that they can be formulated in a similar manner. The first, in statistical modeling, seeks to reconcile observed statistics by suitably and minimally perturbing prior dynamics. The second seeks to optimally select {a subset of available} sensors and actuators for control purposes. To address modeling and control of large-scale systems we develop a unified algorithmic framework using proximal methods. Our customized algorithms exploit problem structure and allow handling statistical modeling, as well as sensor and actuator selection, for substantially larger scales than what is amenable to current general-purpose solvers. {We establish linear convergence of the proximal gradient algorithm, draw contrast between the proposed proximal algorithms and alternating direction method of multipliers, and provide examples that illustrate the merits and effectiveness of our framework.}
	\end{abstract}
	

\vspace*{-2.ex}
\section{Introduction}
\label{sec.intro}

Convex optimization has had tremendous impact on many disciplines, including system identification and control design~\cite{boyelgferbal94,dulpag00,fazhinboy01,boyvan04,fazhinboy04,liuvan09,jovdhiEJC16}. The forefront of research points to broadening the range of applications as well as sharpening the effectiveness of algorithms in terms of speed and scalability. The present paper focuses on two representative control problems, statistical control-oriented modeling and sensor/actuator selection, that are cast as convex programs. A range of modern applications require addressing these over increasingly large parameter spaces, placing them outside the reach of standard solvers. A contribution of the paper is to formulate such problems as regularized semi-definite programs (SDPs) and to develop customized optimization algorithms that scale favorably with~size. 

Modeling is often seen as an inverse problem where a search in parameter space aims to find a parsimonious representation of data. For example, in the control-oriented modeling of fluid flows, it is of interest to improve upon dynamical equations arising from first-principles (e.g., linearized Navier-Stokes equations), in order to accurately replicate observed statistical features that are estimated from data. To this end, a perturbation of the prior model can be seen as a feedback gain that results in dynamical coupling between a suitable subset of parameters~\cite{zarchejovgeoTAC17,zarjovgeoJFM17,zargeojovARC20}. On the flip side, active control of large-scale and distributed systems requires judicious placement of sensors and actuators which again can be viewed as the selection of a suitable feedback or Kalman gain. In either modeling or control, the selection of such gain matrices must be guided by optimality criteria as well as simplicity (low rank or sparse architecture). We cast both types of problems as optimization problems that utilize suitable convex surrogates to handle complexity. The use of such surrogates is necessitated by the fact that searching over all possible architectures is combinatorially prohibitive.

Applications that motivate our study require scalable algorithms that can handle large-scale problems. While the optimization problems that we formulate are SDP representable, e.g., for actuator selection, worst-case complexity of generic solvers scales as the sixth power of the sum of the state dimension and the number of actuators. Thus, solvers that do not exploit the problem structure cannot cope with the demands of such large-scale applications. This necessitates the development of customized algorithms that are pursued herein.

Our presentation is organized as follows.
In Section~\ref{sec.applications}, we describe the modeling and control problems that we consider, provide an overview of literature and the state-of-the-art, and highlight the technical contribution of the paper. In Section~\ref{sec.prelim}, we formulate the {\em minimum energy covariance completion\/} (control-oriented modeling) and {\em sensor/actuator selection\/} (control) problems as nonsmooth SDPs. In Section~\ref{sec.algorithms}, we present a customized Method of Multipliers (MM) algorithm for covariance completion. An essential ingredient of MM is the Proximal Gradient (PG) method. {We also use the PG method for sensor/actuator selection and establish its convergence rate.} In Section~\ref{sec.example}, we offer two motivating examples for actuator selection and covariance completion and discuss computational experiments. We conclude with a brief summary of the results and future directions in Section~\ref{sec.remarks}.

\vspace*{-2ex}
\section{Motivating applications and contribution}
\label{sec.applications}

We consider dynamical systems with additive stochastic disturbances. In the first instance, we are concerned with a modeling problem where the statistics are not consistent with a prior model that is available to us. In that case, we seek to modify our model in a parsimonious manner (a sparse and structured perturbation of the state matrix) so as to account for the partially observed statistics. In the second, we are concerned with the control of such stochastic dynamics via a collection of judiciously placed sensors and actuators. Once again, the architecture of the (now) control problem calls for the selection of sparse matrix gains that effect control and estimation. These problems are explained next.

\vspace*{-2ex}
\subsection{Statistical modeling and covariance completion}

It is well-established that the linearized Navier-Stokes (NS) equations driven by stochastic excitation can account for qualitative~\cite{farioa93,bamdah01,jovbamJFM05,moajovJFM12,ranzarhacjovPRF19b} and quantitative~\cite{zarjovgeoJFM17,zargeojovARC20} features of shear flows. The value of such models has been to provide insights into the underlying physics as well as to guide control design. A significant recent step in this direction was to recognize~\cite{zarjovgeoJFM17} that {\em colored-in-time\/} excitation can account for features of the flow field that {\em white\/} noise in earlier literature cannot~\cite{jovbamCDC01}. Furthermore, it has been pointed out that the effect of colored-in-time excitation is {\em equivalent\/} to white-in-time excitation together with a {\em structural perturbation of the system dynamics\/}~\cite{zarchejovgeoTAC17,zarjovgeoJFM17}. Such structural perturbations may reveal salient dynamical couplings between variables and, thereby, enhance understanding of basic physics~\cite[Section 6.1]{zarjovgeoJFM17}; see~\cite{zargeojovARC20} for a review of covariance completion problems and its relevance in stochastic dynamical modeling of turbulent~flows.

\begin{figure}[t!]
\centering
%
%
%
%
%
\input{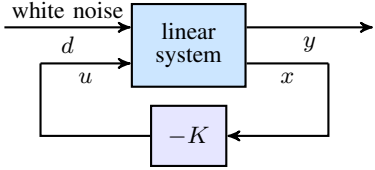}
\small

\noindent
\begin{tikzpicture}[scale=.8, auto, >=stealth']
    
     \node[block, minimum height = 1.1cm, top color=RoyalBlue!20, bottom color=RoyalBlue!20] (sys1) {\begin{tabular}{c}linear\\[-.05cm] system\end{tabular}};
     
     \node[block, minimum height = .8cm, top color=blue!10, bottom color=blue!10] (sys2) at ($(sys1.south) - (0cm,.8cm)$) {$-K$};
                          
     \node[] (R1) at ($(sys1.east) + (1.5cm,-0.3cm)$){};
     
%
     
    
    \draw [connector] ($(sys1.west) + (-2.1cm,.3cm)$) -- node [midway, above] {white noise} node [midway, below] {$d$} ($(sys1.west) + (0cm,.3cm)$);
                    	
    
    \draw [line] (sys2.west) -|  ($(sys1.west) + (-1.5cm,-.3cm)$);
    
    \draw [connector] ($(sys1.west) + (-1.5cm,-.3cm)$) -- node [midway, below] {$u$} ($(sys1.west) + (0cm,-.3cm)$);
    
    \draw [connector] ($(sys1.east) + (0cm,.3cm)$) -- node [midway, below] {$y$} ($(sys1.east) + (2.1cm,.3cm)$);
    
    \draw [line] ($(sys1.east) + (0cm,-.3cm)$) -- node[midway, below] {$x$} (R1.west);
    
    \draw [connector] (R1.west) |- (sys2.east);
    
                                       
\end{tikzpicture}
\vspace{-.1cm}
\caption{A feedback connection of an LTI system with a static gain matrix that is designed to account for the sampled steady-state covariance~$X$.}
\label{fig.sysfeedback}
\end{figure}

{These insights and reasoning motivate} an optimal state-feedback synthesis problem~\cite{zarjovgeoCDC16} to identify dynamical couplings that bring consistency between the model and the observed statistics. Model parsimony dictates a penalty on the complexity of structural perturbations and leads to an optimization problem that involves a composite cost function
\begin{align}
	f(X,K) 
	\;+\;
	\gamma\, g(K)
	\label{eq.compositeJ}
\end{align}
subject to stability of the system in Fig.~\ref{fig.sysfeedback}. Here, $X$ denotes a state covariance matrix and $K$ is a state-feedback matrix. The function $f(X,K)$ penalizes variance and control energy while $g(K)$ is a sparsity-promoting regularizer which penalizes the number of nonzero rows in $K$; sparsity in the rows of $K$ amounts to a {reduced number of feedback couplings that modify the system dynamics.} In addition, state statistics may be partially known, in which case a constraint $X_{ij} = G_{ij}$ for $(i,j) \in \cI$ is added, where the entries of $G$ represent known entries of $X$ for indices in $\cI$.

The resulting {\em minimum-control-energy covariance completion problem} can be cast as an SDP which, for small-size problems, is readily solvable using standard software. {A class of similar problems have been proposed in the context of stochastic control~\cite{HotSke87,yasskegri93,grikarske94,chegeopav16b} and of output covariance estimation~\cite{linjovTAC09,zorfer12} which, likewise and for small-size, are readily solvable by standard software.}

\vspace*{-2ex}
\subsection{Sensor and actuator selection}

The selection and proper placement of sensors/actuators impacts the performance of closed-loop control systems; making such a choice is a nontrivial task even for systems of modest size. Previous work on actuator/sensor placement either relies on heuristics or on greedy algorithms and convex relaxations.

The benefit of a particular sensors/actuator placement is typically quantified by properties of the resulting {observability/controllability} and the  selection process is guided by indicators of diminishing return in performance near optimality~\cite{sumcorlyg16,tzorahpapjad16}. However, metrics on the performance of Kalman filters and other control objectives have been shown to lack supermodularity~\cite{zhaayosun17,ols18}, which hampers the effectiveness of greedy approaches in sensor/actuator~selection.

The literature on different approaches includes convex formulations for sensor placement in problems with linear measurements~\cite{josboy09}, maximizing the trace of the Fisher information under constraints when dealing with correlated measurement noise~\cite{liuchefarmasleuvar16}, and a variation of optimal experiment design for placing measurement units in power networks~\cite{kekgiawol12}. Actuator selection via genetic algorithms has also been explored~\cite{rog00}. Finally, a non-convex formulation of the joint sensor and actuator placement was advanced in~\cite{konyatino90,hirdokobi00} and was recently applied to the linearized Ginzburg-Landau equation~\cite{cherowJFM11}.

{Herein, we cast our placement problem as one of optimally selecting a subset of potential sensors or actuators which, in a similar manner as our earlier modeling problem, involves the minimization of a nonsmooth composite function as in~\eqref{eq.compositeJ}.}
{More specifically, we utilize the sparsity-promoting framework developed in~\cite{farlinjovACC11,linfarjovACC12,linfarjovTAC13admm} to enforce block-sparse structured observer/feedback gains and select sensors/actuators.}

The algorithms developed in~\cite{linfarjovTAC13admm} have been used for sensor selection in target tracking~\cite{masfarvar12} and in periodic sensor scheduling in networks of dynamical systems~\cite{liufarmasvar14}. However, they were developed for general problems, without exploiting a certain hidden convexity in sensor/actuator selection. Indeed, for the design of row-sparse feedback gains, the authors of~\cite{polkhlshc13} introduced a convex SDP reformulation of the problem formulated in~\cite{linfarjovTAC13admm}. Inspired by~\cite{linfarjovTAC13admm}, the authors of~\cite{munpfiwol14} extended the SDP formulation to $\cH_2$ and $\cH_\infty$ sensor/actuator placement problems for discrete time LTI systems. Their approach utilizes standard SDP-solvers with re-weighted $\ell_1$-norm regularizers. In the present paper, we integrate several of these ideas. In particular, we borrow group-sparsity regularizers from statistics~\cite{yualin06} and develop efficient customized proximal algorithms for the resulting SDPs.

\vspace*{-1ex}
\subsection{Main contribution}

{In the present paper, we highlight the structural similarity between statistical modeling and sensor/actuator selection, and develop a unified algorithmic framework for handling large-scale problems. Proximal algorithms are utilized to address the non-differentiability of the sparsity-promoting term $g(K)$ in the objective function. We exploit the problem structure, implicitly handle the stability constraint on state covariances and controller gains by expressing one in terms of the other, and develop a customized proximal gradient algorithm that scales with the third power of the state-space dimension. We prove linear convergence for the proximal gradient algorithm with fixed step-size and propose an adaptive step-size selection method that can improve convergence. We also discuss initialization techniques and stopping criteria for our algorithms, and provide numerical experiments to demonstrate the effectiveness of our approach relative to existing methods.}

	\vspace*{-2ex}
\section{Problem formulation}
\label{sec.prelim}

Consider a linear time-invariant (LTI) system with state-space representation
\begin{align}
	\ba{rcl}
	\dot{x} 
	& \!\!\! = \!\!\! &
	A \, x  \,+\, B \, u \,+\, {d}
	\\[0.1cm]
	y 
	& \!\!\! = \!\!\! &
	C \, x
	\ea
	\label{eq.LTIsys}
\end{align}
where $x(t) \in \bbC^n$ is the state vector, $y(t) \in \bbC^p$ is the output, $u(t) \in \bbC^m$ is the control input, and $d(t)$ is a white stochastic process with zero-mean and the covariance matrix $V \succ 0$, {$\mathbf{E}(d(t)d^*(\tau)) = V \delta (t - \tau)$. Here, $\mathbf{E}$ is the expected value,} $B \in \bbC^{n\times m}$ is the input matrix with $m \leq n$, $C \in \bbC^{p\times n}$ is the output matrix, and the pair ($A, B$) is controllable. {The choice of the state-space is motivated by spatially distributed systems where the application of the spatial Fourier transform naturally leads to complex-valued quantities in~\eqref{eq.LTIsys}; e.g., see~\cite{bampagdah02}.}

We consider two specific applications, one that relates system identification and covariance completion, and another that focuses on {actuator} selection in a control problem. Both can be cast as the problem to select a stabilizing state-feedback control law, $u=-Kx$, that utilizes few input degrees of freedom in the sense that the matrix $K$ has a large number of zero rows. At the same time, the closed-loop system
\[
	\dot{x} 
	\; = \;
	\left( A \,-\,B\,K \right) x 
	\; + \; 
	d
\]
shown in Fig.~\ref{fig.sysfeedback} is consistent with partially available state-correlations and/or is optimal in a quadratic sense.

{
More specifically, if
	\[
	X \,\DefinedAs\, \lim_{t\,\rightarrow \,\infty} \mathbf{E} \left( x(t)\, x^*(t) \right)
	\]
denotes the stationary state-covariance of the controlled system,
the pertinent quadratic cost is
\begin{align}
\label{eq.H2norm}
	\ba{rrl}
		\!\!\!
		f(X,K)
		& \!\!\! \DefinedAs \!\!\! &
		\trace \left( Q\,X \,+\, K^*R\,K X \right)
		\\[.15cm]
		\!\!\!
		& \!\!\! = \!\!\! &
		\ds{\lim_{t\,\rightarrow \,\infty} \mathbf{E} \left( x^*(t)\,Q\, x(t) \, + \, u^*(t)\,R\, u(t) \right)}
	\ea
\end{align}
whereas $Q = Q^* \succ 0$ and $R = R^* \succ 0$ specify penalties on the state and control input, respectively. Both stability of the feedback dynamics and consistency with the state covariance $X$ reduce to an algebraic constraint on $K$ and $X$, namely,
\begin{align}
\label{eq:lyapKX}
	(A - B\, K)  X \;+\; X (A - B\, K)^* \,+\; V  \; = \; 0.
\end{align}
Finally, the number of non-zero rows of $K$ can be seen as the number of active degrees of freedom of the input $u=-Kx$. The choice of such a $K$, with few non-zero rows is sought via minimization of a non-smooth composite objective function in Problem~\ref{prob1}, where 
\begin{align}
\label{eq:g}
	g(K)
	\;\DefinedAs\;
	\ds{\sum^n_{i \, = \, 1} w_i\, \norm{\mre_i^* K}_2}
\end{align}
is a regularizing term that promotes row-sparsity of $K$~\cite{yualin06}, $w_i$ are positive weights, and $\mre_i$ is the $i$th unit vector in $\bbR^m$. 

	\begin{problem}
Minimize $f(X,K) + \gamma\, g(K)$, subject to \eqref{eq:lyapKX}, $X \succ 0$, and, possibly, constraints on the values of specified entries of $X$, $X_{ij}=G_{ij}$ for $(i,j) \in {\mathcal I}$, where a set of pairs ${\mathcal I}$ and the entries $G_{ij}$ are given.
\label{prob1}
\end{problem}

	\vsp

{In this problem, $\gamma>0$ specifies the importance of sparsity, and $\mathcal I$ specifies indices of available covariance data. A useful variant of the constraint on the entries of $X$, when, e.g., statistics of output variables are estimated, can be expressed~as
\begin{align}
	\left(C X C^* \right)_{ij} 
	\; = \;
	G_{ij}
	~\mbox{ for } 
	(i,j) \, \in \, {\mathcal I}
	\label{eq:entries}
\end{align}
We next explain how Problem~\ref{prob1} relates to the two aforementioned topics of covariance completion and actuator selection.}

	\vspace*{-2ex}
\subsection{Covariance completion and model consistency}
\label{sec.cc}

In many problems, it is often the case that a model is provided for a given process which, however, is inconsistent with new data. In such instances, it is desirable to revise the dynamics by a suitable perturbation to bring compatibility between model and data.  The data in our setting consists of statistics in the form of a state covariance $X$ for a linear model
\beq
\label{eq:unperturbed}
	\dot{x} 
	 \;=\;
	A\, x  \; + \; d
\eeq
with white noise input $d$. 

We postulate and deal with a further complication when the data is incomplete. More specifically, we allow $X$ to be only partially known. Such an assumption is motivated by fluid flow applications that rely on the linearized NS equations~\cite{zarjovgeoJFM17}. In this area both the numerical and experimental determination of all entries of $X$ is often prohibitively expensive. Thus, the problem to bring consistency between data and model can be cast in the form of Problem~\ref{prob1}, where we seek a completion of the missing entries of $X$ along with a perturbation $\Delta \DefinedAs -BK$ of the system dynamics~\eqref{eq:unperturbed}, into
\begin{align*}
	\dot{x} 
	 \; = \;
	(A \,+\, \Delta)\, x  \; + \; d.
\end{align*}
The assumed structure of $\Delta$ is without loss of generality, and the choice of $B$ may incorporate added insights into the strength and directionality of possible couplings between state variables. {It should be noted that a full-rank matrix $B$ that allows the perturbation signal $Kx$ to manipulate all degrees of freedom can lead to the complete cancellation of the original dynamics $A$; see~\cite[Section III]{zarchejovgeoTAC17} for details. Then, when seeking a suitable perturbation, it is also natural to impose a penalty on the average quadratic size of signals $Kx$. This brings us into the setting of Problem~\ref{prob1}, where the choice of most suitable perturbation is determined by the optimization criterion. Once again, the row-sparsity promoting penalty $g(K)$ impacts the choice of feedback couplings that need to be introduced to modify \mbox{the dynamical generator $A$~\cite{zarjovgeoCDC16}.}}

	\vspace*{-2ex}
\subsection{Actuator selection}
\label{sec.actsel}

As is well-known, the unique optimal control law that minimizes the steady-state variance~\eqref{eq.H2norm} of system~\eqref{eq.LTIsys} is a static state-feedback $u=-K x$. The optimal gain $K$ and the corresponding state covariance $X$ can be obtained by {minimizing $f(X,K)$, over $K \in \bbC^{m\times n}$, and positive definite $X =X^*\in \bbC^{n\times n}$.}  The solution can also be obtained by solving an algebraic Riccati equation arising from the KKT conditions of this optimization problem. In general, {$K$ is populated by non-zero entries, implying that all ``input channels'' (i.e., all entries of $u$) would be active.} Since the columns of $B$ encode the effect of individual ``input channels'', representing location of actuators, a subselection that is affected by the row-sparsity promoting regularizer in Problem~\ref{prob1}, amounts to actuator selection amongst available options. A dual formulation can be cast to address sensor selection and can be approached in a similar manner; see Appendix~\ref{sec.sensel}.}

	\vspace*{-2ex}
\subsection{Change of variables and SDP representation}
\label{sec.cov}

The constraint $X \succ 0$ in Problem~\ref{prob1} allows for a standard change of variables $Y\DefinedAs K X$ to replace $K$ in $f(X,K) = \trace \left( Q X + K^*R \, K X \right)$. This yields the function
\begin{align}
\label{eq.fhat}
	\ba{rcl}
	f(X,Y)
	& \!\!\! = \!\!\; &
	\trace\left( Q\, X \,+\, Y^* R\, Y X^{-1} \right)
	\ea
\end{align}
which is jointly convex in ($X,Y$).
Further, the row-sparsity of $K$ is equivalent to the row-sparsity of $Y$~\cite{polkhlshc13}. This observation leads to the convex reformulation of Problem \ref{prob1} (incorporating the more general version of constraints \eqref{eq:entries}) as follows.
	
	\vsp
	
\begin{problem}\label{prob2}
Minimize $f (X,Y) + \gamma\, {\sum_{i} w_i \norm{\mre_i^* Y}_2}$ over a Hermitian matrix $X \in \bbC^{n\times n}$ and $Y \in \bbC^{m\times n}$, subject to:
\begin{align*}
	\ba{rcl}
	A\, X \,+\, X\,A^* \,-\, B\, Y -\, Y^*B^* +\, V  
	&\!\!=\!\!& 
	0
	\\[.15cm]
	\left(1 \,-\, \delta \right)\left[\, \left(C X C^* \right)\circ E \,-\, G\, \right] 
	&\!\!=\!\!& 
	0
	\\[.15cm]
	 X 
	 &\!\! \succ \!\!& 
	 0
	 \ea
\end{align*}
where 
	\begin{align}
	\delta 
	\;=\;
	\left\{
	\ba{ll}
	0,
	&
	\text{for covariance completion} 
	\\[.1cm]
	1,
	&
	\text{for actuator selection.}
	\ea
	\right.
	\non
\end{align}
The symbol $\circ$ denotes elementwise matrix multiplication, and $E$ is the structural identity matrix,
\begin{align}
	E_{ij} \;=\;
	\left\{
	\ba{ll}
	1,
	&
	\text{if} ~ G_{ij} ~ \text{is available}
	\\[.1cm]
	0,
	&
	\text{if} ~ G_{ij} ~ \text{is unavailable.}
	\ea
	\right.
	\non
\end{align}
\end{problem}

As explained earlier, the matrices $A$, $B$, $C$, $G$, and $V$ are problem data. From the solution of Problem~\ref{prob2}, the optimal feedback gain matrix can be recovered as $K=Y X^{-1}$. We note that the optimization of $f$ can be expressed as an SDP. Specifically, the Schur complement can be used to characterize the epigraph of $\trace\left( R\, Y X^{-1}Y^*  \right)$ via the convex {constraint
\begin{align*}
	  \left[\begin{matrix} W & R^{1/2}\,Y
	  \\
	   Y^*R^{1/2} & X
	\end{matrix}\right] 
	\;\succeq\; 
	0
\end{align*}
and $\trace \, (W)$,} where $W$ is a matrix variable and the joint convexity of $\trace\left( R\, Y X^{-1}Y^*  \right)$ in $(X,Y)$ follows~\cite{boyvan04}.

We also note that although the row-sparsity patterns of $Y$ and $K$ are equivalent, the weights $w_i$ are not necessarily the same in the respective expressions in Problems~\ref{prob1} and~\ref{prob2}. In practice, the weights are iteratively adapted to promote row-sparsity; see Section~\ref{sec.reweighting_polishing}. Problem~\ref{prob2} can be solved efficiently using general-purpose solvers for small number of variables. To address larger problems, we next exploit the structure and develop optimization algorithms based on the proximal gradient algorithm and the method of multipliers.

\section{Customized algorithms}
\label{sec.algorithms}

In this section, we describe the steps through which we solve {Problem~\ref{prob2}}, identify the essential input channels in $B$, and subsequently refine the solutions based on the identified sparsity structure. For notational compactness, we write the linear constraints in {Problem~\ref{prob2}} as
\begin{align*}
	\ba{rcl}
	\cA_1 (X) \,-\, \cB(Y) +\, V  
	&\!\! = \!\!& 
	0
	\\[.1cm]
	(1 \, - \, \delta) 
	\left[
	\,
	\cA_2 (X) \,-\, G 
	\,
	\right]
	&\!\! = \!\!& 
	0
	 \ea
\end{align*}
where the linear operators $\cA_1$: $\bbC^{n\times n} \to \bbC^{n\times n}$, $\cA_2$: $\bbC^{n\times n} \to \bbC^{p\times p}$ and $\cB$: $\bbC^{m\times n} \to \bbC^{n\times n}$ are given by
\begin{align*}
	\ba{rcl}
		\cA_1 (X)
		&\!\! \DefinedAs \!\!&
		A\,X \,+\, X\, A^*
		\\[.15cm]
		\cA_2 (X)
		&\!\! \DefinedAs \!\!&
		\left(C X C^* \right)\circ E
		\\[.15cm]
		\cB (Y)
		&\!\! \DefinedAs \!\!&
		B\,Y +\, Y^* B^*.
	\ea
\end{align*}

	\vspace*{-2ex}
\subsection{Elimination of variable $X$}
\label{sec.elimX} 

For any $Y$, there is a unique $X$ that solves the equation
\begin{align}
\label{eq.lyap-X-Y}
	\cA_1 (X) \,-\, \cB(Y) \,+\, V  \, = \, 0
\end{align}
if and only if the matrices $A^*$ and $-A$ do not have any common eigenvalues~\cite{horjoh12}. When this condition holds, we can express the variable $X$ as {an affine} function of $Y$, 
\begin{align}
	X(Y)
	\;=\;
	\cA_1^{-1} (\cB(Y) \,-\, V )
	\label{eq.XofY}
\end{align}
and restate {Problem~\ref{prob2}} as
\begin{align}
	\ba{cl}
	\minimize\limits_{Y}
	&
	f(Y) \,+\, \gamma\, g(Y)
	\\[.25cm]
	\subject 
	&
	\left(1 \,-\, \delta \right) \left[\,\cA_2 (X(Y)) \,-\, G\, \right] \;=\; 0
	\\[.15cm]
	&
	X (Y) \,\succ\, 0.
	 \ea
	\label{eq.MCC-1}
\end{align}
The smooth part of the objective function in~\eqref{eq.MCC-1} is given by
\begin{align}
\label{eq.smooth_f}
	\ba{rcl}
	f(Y) 
	&\!\! \DefinedAs \!\!& 
	\trace\left( Q\, X(Y) \,+\, Y^*R\, Y X^{-1}(Y)  \right)
	\ea
\end{align}
and the regularizing term is
	\begin{align}
		\label{eq.gDef}
	\ba{rcl}
	g(Y)
	&\!\! \DefinedAs \!\!& 
	\ds{\sum^n_{i \, = \, 1} w_i\,\norm{\mre_i^* Y}_2}.
	\ea
\end{align}
Since optimization problem~\eqref{eq.MCC-1} is equivalent to {Problem~\ref{prob2} constrained} to the {affine equality}~\eqref{eq.XofY}, it remains convex.

{
When the matrix $A$ is Hurwitz, expression~\eqref{eq.XofY} can be cast in terms of the well-known integral representation, 		
\begin{align*}
		X (Y)
		\; = \;
		\int_{0}^{\infty}
		\mre^{A t}
		\,
		(V \, - \, B\,Y -\, Y^*B^*)
		\,
		\mre^{A^* t}
		\,
		\mrd t.
\end{align*}
Even for unstable open-loop systems, the operator $\cA_1$ is invertible if the matrices $A^*$ and $-A$ do not have any common eigenvalues. In our customized algorithms, we numerically evaluate the action of $\cA_1^{-1}$ on the current iterate by solving the corresponding Lyapunov equation which requires making the following assumption.

	\vsp

\begin{myass}
	The operator $\cA_1$ is invertible.
\end{myass}
	\vsp
	
Appendix~\ref{sec.noninvertA} provides a method to handle cases where this assumption does not hold.}

	\vspace*{-3ex}
\subsection{Proximal gradient method for actuator selection}
\label{sec.ProxGrad}

The proximal gradient (PG) method generalizes gradient descent to composite minimization problems in which the objective function is the sum of a differentiable and non-differentiable component~\cite{becteb09,parboy13}. It is most effective when the proximal operator associated with the nondifferentiable component is easy to evaluate; many common regularization functions, such as the $\ell_1$ penalty, nuclear norm, and hinge loss, satisfy this condition. {Herein, we present details of a customized variant of the PG method for solving~\eqref{eq.MCC-1} with $\delta=1$. In Algorithm~\ref{alg.PG}, we follow the recommendations of~\cite{becteb09,golstubar14} for choosing the step-size and stopping criterion.}

The PG method for solving~\eqref{eq.MCC-1} with $\delta=1$ is given by
\begin{align}
	\label{eq.PGiter1}
	Y^{k+1}
	~\DefinedAs~
	\prox_{\beta_k g} 
	\!
	\left( Y^k \,-\, \alpha_k\, \nabla f(Y^k) \right)
\end{align}
where $Y^k$ is the $k$th iterate, $\alpha_k >0$ is the step-size, {and $\beta_k \DefinedAs \gamma \alpha_k$. The proximal operator of a real-valued proper, closed, convex function $h$ is defined as~\cite{baucom11}
\begin{align}
	\label{eq.prox-op}
	 \prox_{h} (V)
	 ~\DefinedAs~
	 \argmin\limits_Y
	 ~ 
	 \left(
	 h(Y) \;+\; \dfrac{1}{2} \, \norm{Y \,-\, V}_F^2
	 \right).
\end{align}
where} $\norm{\cdot}_F$ is the Frobenius norm. For the row-sparsity regularizer, the proximal operator {of the function $\beta g$} is determined by the soft-thresholding operator which acts on the rows of the matrix~$V$,
{\begin{align*}
	\cS_\beta (\mre_i^* V) 
	\;=\;
	\left\{
	\ba{rrcl}
		\hspace{-.1cm}
		\left(1 - \beta w_i/\norm{\mre_i^* V}_2 \right) \mre_i^* V,
		&~
		\norm{\mre_i^* V}_2 
		&\!\!\! > \!\!\!\!& 
		\beta w_i
		\\[.15cm]
		\hspace{-.1cm}
		0,
		&
		\norm{\mre_i^* V}_2 
		&\!\!\! \leq \!\!\!\!& 
		\beta w_i.
	\ea
	\right. 
\end{align*}}

Proximal update~\eqref{eq.PGiter1} results from a local quadratic approximation of $f$ at iteration $k$, i.e.,
\begin{align}
	\nonumber
	\!\!\!
	Y^{k+1} 
	\; \DefinedAs \;
	\argmin \limits_Y
	&
	~~
	f(Y^k) \, + \, \inner{\nabla f(Y^k)}{Y-Y^k} \, +
	\\[-.15cm]
	\label{eq.pg-expansion}
	&
	~~
	\dfrac{1}{2 \alpha_k}\,\norm{Y \,-\, Y^k}_F^2 \, + \, \gamma\, g(Y)
\end{align}
followed by a completion of squares that brings the problem into the form of~\eqref{eq.prox-op} {with $h \DefinedAs \gamma \alpha_k g$}. Here, $\inner{\cdot}{\cdot}$ denotes the standard matricial inner product $\inner{M_1}{M_2} \DefinedAs \trace \, (M_1^* M_2)$ and the expression for the gradient of $f(Y)$ is provided in Appendix~\ref{sec.gradf}. 

\subsubsection{{Initialization and} choice of step-size in~\eqref{eq.PGiter1}}
\label{sec.PG-step-size}
{The PG algorithm is initialized with ${Y^0} = K^0 X^0$, where $K^0$ is a stabilizing feedback gain and $X^0$ is the corresponding covariance matrix that satisfies~\eqref{eq:lyapKX}. The optimal centralized controller resulting from the solution of the algebraic Riccati equation provides a stabilizing initial condition and the closed-loop stability is maintained via step-size selection in subsequent iterations of Algorithm~\ref{alg.PG}.} At each iteration of the PG method, we determine the step-size $\alpha_k$ via an adaptive Barzilai-Borwein (BB) initial step-size selection~\cite{golstubar14}, i.e.,
\begin{align}
	\label{eq.init-alpha}
	{\alpha_{k,0}}
	\;=\; 
	\left\{
	\ba{ll}
		\alpha_m~~
		& 
		\text{if}~\, \alpha_m/\alpha_s 
		\, > \, 
		1/2
		\\[.1cm]
		\alpha_s \,-\, \alpha_m/2 ~~
		&
		\text{otherwise}
	\ea 
	\right.
\end{align}
followed by backtracking to ensure {closed-loop} stability 
\begin{subequations}
\label{eq.backtracking-conds}
\begin{align}
	\label{eq.stability}
	X(Y^{k+1})
	\;\succ\;
	0
\end{align}
and sufficient descent of the objective function {$f(Y) + \gamma g(Y)$ resulting from} 
\begin{eqnarray}
	f(Y^{k+1}) 
	&\!\!\leq\!\!& 
	f(Y^{k}) \;+\; \inner{\nabla f(Y^k)}{Y^{k+1} - Y^k} 
	\, + 
	\non
	\\[.05cm]
	&&	
	\dfrac{1}{2 \alpha_k} \, \norm{Y^{k+1} - Y^k}_F^2.
	\label{eq.suff-descent}
\end{eqnarray}
\end{subequations}
Similar strategies as~\eqref{eq.suff-descent} were used {in~\cite[Section 3]{becteb09}}. Here, the ``steepest descent'' step-size $\alpha_s$ and the ``minimum residual" step-size $\alpha_m$ are given by,
\[
	\ba{rcl}
	\alpha_s
	&\!\!=\!\!& 
	\dfrac{\inner{Y^k - Y^{k-1}}{Y^k - Y^{k-1}}}{\inner{Y^k - Y^{k-1}}{\nabla f(Y^k) - \nabla f(Y^{k-1})}}
	\\[.4cm]
	\alpha_m
	&\!\!=\!\!& 
	\dfrac{\inner{Y^k - Y^{k-1}}{\nabla f(Y^k) - \nabla f(Y^{k-1})}}{\inner{\nabla f(Y^k) - \nabla f(Y^{k-1})}{\nabla f(Y^k) - \nabla f(Y^{k-1})}}.
	\ea
\]
If $\alpha_s < 0$ or $\alpha_m < 0$, the step-size from the previous iteration is used; see~\cite[Section 4.1]{golstubar14} for additional details.

\vsp
\subsubsection{Stopping criterion}
\label{sec.stopping-criterion}

We employ a combined condition that terminates the algorithm when either the relative residual
\begin{align*}
	r_r^{k+1}
	\;=\;
	\dfrac{\norm{r^{k+1}}}{\max\{\norm{\nabla f(Y^{k+1})}, \norm{(\hat{Y}^{k+1} - Y^{k+1})/\alpha_k} \} \,+\, \eps_r},
\end{align*}
or the normalized residual
\begin{align*}
	r_n^{k+1}
	\; = \;
	\dfrac{\norm{r^{k+1}}}{\norm{r^1} \;+\; \eps_n}
\end{align*}
are smaller than a desired tolerance. Here, $\eps_r$ and $\eps_n$ are small positive constants, the residual is defined as
\begin{align*}
	r^{k+1}
	\;\DefinedAs\;
	\nabla f(Y^{k+1}) \,+\, (\hat{Y}^{k+1} \,-\, Y^{k+1})/\alpha_k
\end{align*}
and $\hat{Y}^{k+1} \DefinedAs Y^k - \alpha_k \nabla f(Y^k)$. While achieving a small $r_r$ guarantees a certain degree of accuracy, its denominator nearly vanishes when $\nabla f(x^\star) = 0$, which happens when $0 \in \partial g(Y^\star)$. In such cases, $\| r_n \|$ provides an appropriate stopping criterion; see~\cite[Section 4.6]{golstubar14} for additional details.

\begin{algorithm}
\caption{Customized PG Algorithm}
\label{alg.PG}
\begin{algorithmic}
\STATE \textbf{input:} $A$, $B$, $V$, $Q$, $R$, $\gamma > 0$, positive constants {$\eps_r$, $\eps_n$}, tolerance $\eps$, and backtracking constant $c\in (0, 1)$.
\STATE \textbf{initialize:} $k=0$, $\alpha_{0,0} = 1$, $r_r^0 = 1$, $r_n^0 = 1$, choose {${Y^0} = K^0 X^0$ where $K^0$ is a stabilizing feedback gain with corresponding covariance matrix $X^0$.}
\vspace*{0.1cm}
\STATE \textbf{while:} $r_r^k > {\eps}$ or $r_n^k > \eps$
\\[.1cm]
~\,\quad compute $\alpha_k$: largest feasible step in $\{c^j \alpha_{k,0}\} _{j=0,1,\ldots}$
\\[.05cm]
~\,\quad such that {$Y^{k+1}$ satisfies~\eqref{eq.backtracking-conds}}
\\[.05cm]
~\,\quad compute $r_r^{k+1}$ and $r_n^{k+1}$
\\[.05cm]
~\,\quad $k = k + 1$
\\[.05cm]
~\,\quad choose $\alpha_{k,0}$ based on~\eqref{eq.init-alpha} 
\vspace*{0.1cm}
\STATE \textbf{endwhile}
\STATE \textbf{output:} $\eps$-optimal solutions, $Y^{k+1}$ and $X(Y^{k+1})$.
\end{algorithmic}
\end{algorithm}

\subsection{{Convergence} of the proximal gradient algorithm}
\label{sec.convergence-PG}

We next analyze the convergence of the PG algorithm for the strongly convex nonsmooth composite optimization problem,
\begin{align}
\ba{cl}
\minimize\limits_{Y}
&
f(Y) \,+\, \gamma\, g(Y).
\ea
\label{eq.general-composite}
\end{align}
The PG algorithm~\eqref{eq.PGiter1} with suitable step-size converges with the linear rate $O(\rho^k)$ for some $\rho \in (0,1)$~if: (i) the function $f$ is strongly convex and smooth (i.e., it has a Lipschitz continuous gradient) {\em uniformly over the entire domain\/}; and (ii) the function $g$ is proper, closed, and convex~\cite[Theorem~10.29]{bec17}. In problem~\eqref{eq.MCC-1}, however, condition (i) does not hold over the function domain
\begin{align}
	\label{eq.Ds}
	\cD_s
	\DefinedAs 
	\{ Y \in \bbC^{m\times n}|\,\cA_1(X(Y))-\cB(Y) = -V,\; X(Y)\succ 0\}
\end{align}
corresponding to stabilizing feedback gains $K = YX^{-1}$. 
To address this issue, we exploit the coercivity~\cite[Definition 11.10]{baucom11} of common regularization functions and establish linear convergence of the PG method for a class of problems~\eqref{eq.general-composite} in which the function $f$ satisfies the following assumption.

	\vsp
	
\begin{myass}\label{ass.sublevel}
	For all scalars $a$, the proper closed convex function $f$ defined over an open convex domain $\cD$ has
	\begin{enumerate}
		\item[(i)] compact sublevel sets
		$
		\cD(a)\DefinedAs\{ Y \in \cD \,|\, f(Y) \le \, a \}
		$; 
		\vspace{.1cm}
		\item[(ii)] an $\La$-Lipschitz continuous gradient over $\cD(a)$;
		\vspace{.1cm}
		\item[(iii)] a strong convexity modulus $\mua>0$ over $\cD(a)$.
	\end{enumerate}
\end{myass}

	\vsp
Proposition~\ref{prop.megaprop} establishes linear convergence of the PG algorithm with sufficiently small fixed step-size. Proofs of all technical results presented here are provided in Appendix~\ref{sec.proofs-sec-conv}.

\vsp

\begin{myprop}
\label{prop.megaprop}
	Let the function $g$ be coercive, proper, closed, and convex and let the function $f$ in~\eqref{eq.general-composite} satisfy conditions~(i) and (ii) in Assumption~\ref{ass.sublevel}. Then, for any initial condition $Y^0 \in \cD$ the iterates $\{Y^k\}$ of the PG algorithm~\eqref{eq.PGiter1} with step-size $\alpha \in [0,1/\La]$ remain in the sublevel set $D(a)$, with $a > f(Y^0) + \gamma\, (g(Y^0)-g(Y))$, for all $Y$. Furthermore, if condition (iii) in Assumption~\ref{ass.sublevel} also holds, then
\begin{align}
\label{eq.linConv}
	\norm{Y^{k+1} \, - \, Y^\star}_F^2
	\; \leq \;
	\left( 1 \,-\, \mua \alpha \right) \norm{Y^k \, - \, Y^\star}_F^2
\end{align}
where $Y^\star$ is the globally optimal solution of~\eqref{eq.general-composite}.
\end{myprop}

We next establish strong-convexity and smoothness for the function $f$ in~\eqref{eq.smooth_f} over its sublevel sets. 
These properties allow us to invoke Proposition~\ref{prop.megaprop} and prove linear convergence for the PG algorithm applied to problem~\eqref{eq.MCC-1} with $\delta=1$.
\vsp
\begin{myprop}
\label{prop.strongConvexity}
The function $f$ in~\eqref{eq.smooth_f} with the convex domain $\cD_s$ given by~\eqref{eq.Ds} satisfies Assumption~\ref{ass.sublevel}.
\end{myprop}
	
	\vsp
Our main result is presented in Theorem~\ref{thm.LQRconvergence}.
	\vsp
	
\begin{mythm}
\label{thm.LQRconvergence}
	For any stabilizing initial condition $Y^0 \in \cD_s$, the iterates of the PG algorithm~\eqref{eq.PGiter1} with step-size $\alpha \in [0,1/\La]$ applied to problem~\eqref{eq.MCC-1} with $\delta=1$ satisfy~\eqref{eq.linConv},
where $\mua$ and $\La$ are the strong convexity modulus and smoothness parameter of the function $f$ over $\cD(a)$ with $a > f(Y^0) + \gamma\, g(Y^0)$.
\end{mythm}
	\vsp
\begin{proof}
	In addition to being proper, closed, and convex, it is straightforward to verify that the function $g$ given by~\eqref{eq.gDef} is coercive, i.e.,
		\begin{align*}
			\lim_{\norm{Y}_F \, \rightarrow \, +\infty} g(Y) \; = \; +\infty.
		\end{align*}
		 Moreover, from the nonnegativity of $g(Y)$, it follows that $a > f(Y^0) + \gamma\, (g(Y^0)-g(Y))$ for all $Y$. Thus, the result follows from combining Propositions~\ref{prop.megaprop}~and~\ref{prop.strongConvexity}.
\end{proof}
	\vsp
	
\begin{remark}

Proposition~\ref{prop.megaprop} proves that the PG algorithm with fixed step-size $\alpha \in (0,1/\La]$ converges at the linear rate $O((1-\mua \alpha)^k)$. A linear rate $O (\rho^k)$ with $\rho = 1 - \min\{1/(\sqrt{2}L_a),c/L_a\}$ can also be guaranteed using the adaptive step-size selection method of Section~\ref{sec.PG-step-size}; see Appendix~\ref{sec.backtracking}.
\end{remark}
	\vsp
The next lemma provides an expression for the smoothness parameter of the function $f$ over its sublevel sets. We note that this parameter depends on problem data.
	\vsp	
\begin{mylem}
	\label{lem.L}
	Over any non-empty sublevel set $\cD(a)$, the gradient $\nabla f(Y)$ is Lipschitz continuous with parameter
	\begin{subequations}
		\begin{align}
		\La
		\;=\;
		\dfrac{2\, \lambda_{\max}(R)}{\nu} \,  \left( 1 \,+\, \dfrac{\sqrt{a} \, \norm{  \cA^{-1}_1 \cB }_2}{\sqrt{\nu\,\lambda_{\min}(R)}} \right)^2
		\label{eq.L}
		\end{align}
		where the positive scalar
		\begin{align}
		\nu
		\;\DefinedAs\;
		\dfrac{\lambda_{\min}^2(V)}{4\,a} \left(\dfrac{\norm{A}_2}{\sqrt{\lambda_{\min}(Q)}} \,+\, \dfrac{\norm{B}_2}{\sqrt{\lambda_{\min}(R)}} \right)^{-2}
		\label{eq.nu}
		\end{align}
		\label{eq.Lnu}
	\end{subequations}
	gives the lower bound $\nu I \preceq X(Y)$ on the covariance matrix.
\end{mylem}
\vsp
\begin{remark}
	While Lemma~\ref{lem.L}  provides an expression for the smoothness parameter, we have recently established an explicit expression for the strong convexity modulus~\cite{mohzarsoljovCDC19}
\begin{align*}
		\mua
		= 
		\dfrac{2\, \lambda_{\min}(R) \lambda_{\min}(Q)}{\left(a^{1/2} + a^{2} \norm{\cB}_2\big(\lambda_{\min}(Q) \lambda_{\min}(V) \sqrt{\nu \lambda_{\min}(R)} \big)^{-1} \right)^2}.
\end{align*}
Based on Theorem~\ref{thm.LQRconvergence}, the explicit expressions for parameters $\La$ and $\mua$ determine a theoretical bound of $1-\mua/\La$ on the linear convergence rate of the PG algorithm with step-size $\alpha_k = 1/\La$. It should be noted that this bound depends on the initial condition $Y^0$ and problem data.
\end{remark}

	\vspace*{-2ex}
\subsection{Method of multipliers for covariance completion}
\label{sec.MM}

We handle the additional constraint in the covariance completion problem by employing the Method of Multipliers (MM). MM is the dual ascent algorithm applied to {a smooth variant of the dual problem} and it is widely used for solving constrained nonlinear programming problems~\cite{ber82,ber99,nocwri06}.

The MM algorithm for constrained optimization problem~\eqref{eq.MCC-1} with $\delta = 0$ is given by, 
 \begin{subequations}
	\label{eq.MM_steps}
	\begin{eqnarray}
	\hspace{-.5cm} Y^{k+1} &\!\!\DefinedAs\!\!& \argmin\limits_{Y} \,  \cl_{\rho_k} (Y;\, \Lambda^k)
	\label{eq.Ymin_MM}
	\\
	\hspace{-.5cm} \Lambda^{k+1} &\!\!\DefinedAs\!\!& \Lambda^k \,+\, \ds{\rho_k \left( \cA_2 (X(Y^{k+1})) \,-\, G \right)}
	\label{eq.dualupdate_MM}
	\end{eqnarray}
\end{subequations}
where $ \cl_\rho$ is the associated augmented Lagrangian,
\[
	\ba{l}
	\cL_\rho (Y; \Lambda) 
	\;=\; 
	f(Y) \;+\; \gamma\, g(Y) ~+
	\\[.15cm]
	~\quad \quad \quad \inner{\Lambda}{\cA_2 (X(Y)) \,-\, G} \;+\; \dfrac{\rho}{2}\, \norm{\cA_2 (X(Y)) \,-\, G}_F^2
	\ea
\]
$\Lambda \in \bbC^{p \times p}$ is the Lagrange multiplier and $\rho$ is a positive scalar. 
The algorithm terminates when the primal and dual residuals are small enough. The primal residual is given as
\begin{subequations}
\label{eq.MM-termination}
\begin{align}
	\label{eq.MM-termination1}
	\Delta_p
	\;=\;
	\norm{\cA_2 (X(Y^{k+1})) \,-\, G}_F
\end{align}
and the dual residual corresponds to the stopping criterion on subproblem~\eqref{eq.Ymin_MM}
\begin{align}
	\label{eq.MM-termination2}
	\Delta_d
	\;=\;
	\min\{r_r, r_n\}
\end{align}
\end{subequations}
where the relative and normal residuals, $r_r$ and $r_n$, are described in Section~\ref{sec.ProxGrad}. 

\vsp
\subsubsection{Solution to the $Y$-minimization problem~\eqref{eq.Ymin_MM}}
For fixed $\{\rho_k, \Lambda^k\}$, minimizing the augmented Lagrangian with respect to $Y$ amounts to finding the minimizer of $\cL_{\rho_k} (Y; \Lambda^k)$ subject to $X(Y) \succ 0$. Since $g(Y)$ is nonsmooth, we cannot use standard gradient descent methods to find the update $Y^{k+1}$. However, similar to Section~\ref{sec.ProxGrad}, a PG method can be used to solve this subproblem iteratively
\begin{align}
	\label{eq.PGiter2}
	Y^{j+1}
	\;=\;
	\prox_{\beta_j g} \left( Y^j \,-\, \alpha_j \nabla F(Y^j) \right)
\end{align}
where $j$ is the inner PG iteration counter, $\alpha_j>0$ is the step-size, $\beta_j \DefinedAs \alpha_j \gamma$, and $F(Y)$ denotes the smooth part of the augmented Lagrangian~$\cL_{\rho_k} (Y; \Lambda^k)$,
\begin{align*}
	\ba{rrl}
	F(Y)
	& \!\!\! \DefinedAs \!\!\! &
	f(Y) \;+\; \inner{\Lambda^k}{\cA_2 (X(Y)) \,-\, G} ~+
	\\[.1cm]
	& \!\!\! \!\!\! &
	\dfrac{\rho_k}{2}\, \norm{\cA_2 (X(Y)) \,-\, G}_F^2.
	\ea
\end{align*}
The expression for the gradient of $F(Y)$ is provided in Appendix~\ref{sec.gradF}. Similar to Section~\ref{sec.ProxGrad}, we combine BB step-size initialization with backtracking {to satisfy conditions~\eqref{eq.backtracking-conds}.}


\subsubsection{Lagrange multiplier update and choice of step-size in~\eqref{eq.dualupdate_MM}}
\label{sec.dualupdate-MM}

Customized MM for covariance completion is summarized as Algorithm~\ref{alg.MM}. We follow the procedure outlined in~\cite[Algorithm 17.4]{nocwri06} for the adaptive update of $\rho_k$. This procedure allows for inexact solutions of subproblem~\eqref{eq.Ymin_MM} and a more refined update of the Lagrange multiplier $\Lambda$ through the adjustment of convergence tolerances on $\Delta_p$ and $\Delta_d$. {Note that standard convergence results for MM depend on the level of accuracy in solving subproblem~\eqref{eq.Ymin_MM}~\cite[Sections 5.3 and 5.4]{ber82}. While we establish linear convergence of the PG algorithm for solving this subproblem, we relegate a detailed convergence analysis for the MM algorithm to future work.}

	\vspace*{-2ex}
\subsection{Computational complexity}

Computation of the gradient in both algorithms involves {evaluation} of $X$ from $Y$ based on~\eqref{eq.XofY}, a matrix inversion, and solution to the Lyapunov equation. {Each of these} take $O(n^3)$ operations as well as an $O(m n^2)$ matrix-matrix multiplication. The proximal operator for the function $g$ amounts to computing the $2$-norm of all $m$ rows of a matrix with $n$ columns, which takes $O(m n)$ operations. These steps are embedded within an iterative backtracking procedure for selecting the step-size $\alpha$. If the step-size selection takes $q_1$ inner iterations the total computation cost for a single iteration of the PG algorithm is $O(q_1 n^3)$. On the other hand, if it takes $q_2$ iterations for the {PG} method to converge, the total computation cost for a single iteration of our customized MM algorithm is $O(q_1 q_2 n^3)$. {In practice, the backtracking constant $c$ is chosen such that $q_1<50$. The computational efficiency of the PG algorithm relative to standard SDP solvers whose worst-case complexity is $O(n^6)$ is thus evident. However, in MM, $q_2$ depends on the required level of accuracy in solving~\eqref{eq.Ymin_MM}. While there is a clear trade-off between this level of accuracy and the number of MM steps, careful analysis of such effects is beyond the scope of the current paper. Nonetheless, in Section~\ref{sec.example-cc}, we demonstrate that relative to ADMM and SDPT3, customized MM can provide significant speedup.}

\begin{algorithm}
\caption{Customized MM Algorithm}
\label{alg.MM}
\begin{algorithmic}
\STATE \textbf{input:} $A$, $B$, $C$, $E$, $G$, $V$, $\gamma > 0$, and tolerances $\eps_p$ and $\eps_d$.
\vspace*{-0.4cm}
\STATE \textbf{initialize:} $k = 0$, $\rho_0 = 1$, $\rho_{\max} = 10^{9}$, $\eps_0 = 1/\rho_0$, $\eta_0 = \rho_0^{-0.1}$, choose {${Y^0} = K^0 X^0$ where $K^0$ is a stabilizing feedback gain with corresponding covariance matrix $X^0$.}
\vspace*{0.1cm}
\STATE \textbf{for} $k = 0,1,2,\dots$
\\[.05cm]
~\,\quad solve~\eqref{eq.Ymin_MM} using a similar PG algorithm to Algorithm~\ref{alg.PG}
\\
~\,\quad such that $\Delta_d \leq \eps_k$.
\\[.1cm]
~\,\quad \textbf{if}  $\Delta_p \leq \eta_k$
\\[.1cm]
~~~~\,\quad \textbf{if}  $\Delta_p \leq \eps_p$ and $\Delta_d \leq \eps_d$
\\[.1cm]
~~\quad\quad\quad \textbf{stop} with approximate solution $Y^{k+1}$
\\
~~~~\,\quad \textbf{else}
\\[-.1cm]
{\small
\[
	\ba{rcl}
		\Lambda^{k+1} 
		&\!\!\!=\!\!\!& 
		\Lambda^k \,+\, \rho_k \left( \cA_2 (X(Y^{k+1})) \,-\, G \right)
		 \\[.05cm]
		 \rho_{k+1}
		 &\!\!\!=\!\!\!& 
		 \rho_k, \quad~
		 \eta_{k+1} \;=\; \max \{ \eta_k\, \rho_{k+1}^{-0.9}, \eps_p \}
		 \\[.05cm]
		 \eps_{k+1} 
		 &\!\!\!=\!\!\!& 
		 \max \{ \eps_k/\rho_{k+1}, \eps_d \}
	\ea
\]
}
\\[-.1cm]
~~~~\,\quad \textbf{endif}
\\
~\,\quad \textbf{else}
\\[-.3cm]
{\small
\[
	\ba{rcl}
		\Lambda^{k+1} 
		&\!\!\!=\!\!\!&
		\Lambda^k
		\\[.05cm]
		\rho_{k+1} 
		&\!\!\!=\!\!\!& 
		\{ 5 \rho_k, \rho_{\max} \}, 
		\quad~
		\eta_{k+1} \;=\; \max \{ \rho_{k+1}^{-0.1}, \eps_p \}
		 \\[.05cm]
		 \eps_{k+1}
		 &\!\!\!=\!\!\!& 
		 \max \{ 1/\rho_{k+1}, \eps_d \}
	\ea
\]
}
\\[-.1cm]
~\,\quad \textbf{endif} 
\vspace*{0.05cm}
\STATE \textbf{endfor}
\STATE \textbf{output:} optimal solutions, $Y^{k+1}$ and $X(Y^{k+1})$.
\end{algorithmic}
\end{algorithm}

	\vspace*{-2ex}
\subsection{Comparison with other methods}
\label{sec.compare_splitting}

One way of dealing with the lack of differentiability of the objective function in~\eqref{eq.MCC-1} is to split the smooth and nonsmooth parts over separate variables and to add an additional equality constraint to couple these variables. This allows for the minimization of the augmented Lagrangian via {the} Alternating Direction Method of Multipliers (ADMM)~\cite{boyparchupeleck11}.

In contrast to splitting methods, the algorithms considered in this paper use the PG method to solve the nonsmooth problem in terms of the primal variable $Y$, thereby avoiding the necessity to update additional auxiliary variables and their corresponding Lagrange multipliers. Moreover, it is important to note that the performance of augmented Lagrangian-based methods is strongly influenced by the choice of $\rho$. In contrast to ADMM, there are principled adaptive rules for updating the step-size $\rho_k$ in MM. Typically, in ADMM, either a constant step-size is used or the step-size is adjusted to keep the norms of primal and dual residuals within a constant factor of one another~\cite{boyparchupeleck11}. Our computational experiments demonstrate that the customized proximal algorithms considered in this paper significantly outperform ADMM.

	\vsp
	
	\begin{remark}
In~\cite{dhijovluoCDC14}, a customized ADMM algorithm was proposed for solving the optimal sensor and actuator selection problems. In this, the structural Lyapunov constraint on $X$ and $Y$ is dualized via the augmented Lagrangian. While this approach {does not rely on the invertibility of operator $\cA_1$ (cf.~\eqref{eq.XofY}),
it involves subproblems that are difficult to solve. Furthermore, as we show in Section~\ref{sec.example}, it performs poorly in practice, especially for large-scale systems.} This is because of higher computational complexity ($O(n^5)$ per iteration) of the ADMM algorithm developed in~\cite{dhijovluoCDC14}. 
	\end{remark}

	\vspace*{-2ex}
\subsection{Iterative reweighting and polishing}
\label{sec.reweighting_polishing}

To obtain sparser structures at lower values of $\gamma$, we follow~\cite{canwakboy08} and implement a reweighting scheme in which we run the algorithms multiple times for each value of $\gamma$ and update the weights as
	$
	w_i^{j+1}
	=
	1/(\norm{\mre_i^* Y^j}_2 + \eps).
	$
Here, $Y^j$ is the solution in the $j$th reweighting step and the small parameter $\eps$ ensures that the weights are well-defined.

After we obtain the solution to problem~\eqref{eq.MCC-1}, we conduct a {\em polishing\/} step to refine the solution based on the identified sparsity structure. For this, we consider the system
\[
 	\dot{x}
	\;=\;
	(A \, - \, B_{\mathrm{sp}} \, K)\, x \;+\; d
 \]
where the matrix $B_{\mathrm{sp}}\in \bbC^{n \times q}$ is obtained by eliminating the columns of $B$ corresponding to the identified row sparsity structure of $Y$, and $q$ denotes the number of retained input channels. For this system, we solve optimization problem~\eqref{eq.MCC-1} with $\gamma=0$. This step allows us to identify {the optimal matrices $Y \in \bbC^{q \times n}$ and $K \in \bbC^{q \times n}$} for a system with a lower number of input channels.

	\vspace*{-2ex}
\section{Computational experiments}
\label{sec.example}

We provide two examples to demonstrate the utility of the optimization framework for optimal {actuator} selection and covariance completion problems and highlight the computational efficiency of our customized algorithms.

	\vspace*{-2ex}
\subsection{Actuator selection}
\label{sec.example-actsel}

The Swift-Hohenberg equation is a partial differential equation that has been widely used as a model for studying pattern formations in hydrodynamics and nonlinear optics~\cite{crohoh93}. Herein, we consider the linearized Swift-Hohenberg equation around its time independent spatially periodic solution~\cite{burkno06}
\[
	\ba{rcl}
	\partial_t\, \psi (t,\xi)
	&\!\!=\!\!&
	- \left( \partial_x^2 \,+\, 1 \right)^2 \psi (t,\xi) \,-\,c\,\psi (t,\xi) \,+\, f\,\psi (t,\xi)
	\\[.15cm] 
	&&
	+\; u(t,\xi) \,+\, d(t,\xi)
	\ea
\]
with periodic boundary conditions on a spatial domain $\xi \in [ 0,\, 2\pi ]$. Here, {the state $\psi(t,\xi)$ denotes the fluctuation field, $u(t,\xi)$ is a spatio-temporal control input, $d(t,\xi)$ is a zero-mean additive white noise,} $c$ is a {constant bifurcation parameter}, and we assume that $f(\xi) \DefinedAs \alpha \cos(\omega \xi)$ with $\alpha \in \bbR$. Finite dimensional approximation {using the spectral collocation method} yields the following state-space representation
\begin{align}
	\label{eq.SH.dynamics}
	\ba{rcl}
		\dot{\psi}
		&\!\!=\!\!&
		A\, \psi \;+\; u \;+\; d.
	\ea
\end{align}

For $c = -0.2$, $\alpha = 2$, and $\omega = 1.25$, the linearized dynamical generator has two unstable modes. We set $Q = I$ and $R = 10 I$ and solve the actuator selection problem (problem~\eqref{eq.MCC-1} with $\delta=1$) for 32, 64, 128 and 256 discretization points and for various values of the regularization parameter $\gamma$. For $\gamma = 10$, Table~\ref{table.SH_CVX_PG_ADMM} compares the proposed PG algorithm against {SDPT3~\cite{SDPT3}} and the ADMM algorithm of~\cite{dhijovluoCDC14}. Both PG and ADMM were initialized with $Y^0 = K_c X_c$, where $K_c$ and $X_c$ solve the algebraic Riccati equation which specifies the optimal centralized controller. This choice guarantees that $X (Y^0 ) \succ 0$.  All algorithms were implemented in Matlab and executed on a $2.9$ GHz Intel Core i5 processor with $16$ GB RAM. {The parser CVX~\cite{cvx} was used to call the solver SDPT3.} The algorithms terminate when an iterate achieves a certain distance from optimality, i.e., $\norm{X^k - X^\star}_F/\norm{X^\star}_F < \eps$ and $\norm{Y^k - Y^\star}_F/\norm{Y^\star}_F < \eps$. The choice of $\eps = 10^{-3}$ guarantees that the value of the objective function is within $0.01\%$ of optimality. For $n=256$, CVX failed to converge. In this case, iterations are run until the relative or normalized residuals defined in Section~\ref{sec.stopping-criterion} become smaller than $10^{-2}$.

\begin{figure}
\begin{tabular}{rcc}
	\hspace{-.6cm}
	\begin{tabular}{c}
	\vspace{.4cm}
	\rotatebox{90}{\footnotesize $\norm{Y^k-Y^\star}_F/\norm{Y^\star}_F$}
	\end{tabular}
	&
	\hspace{-.95cm}
	\begin{tabular}{c}
         	\includegraphics[width=.25\textwidth]{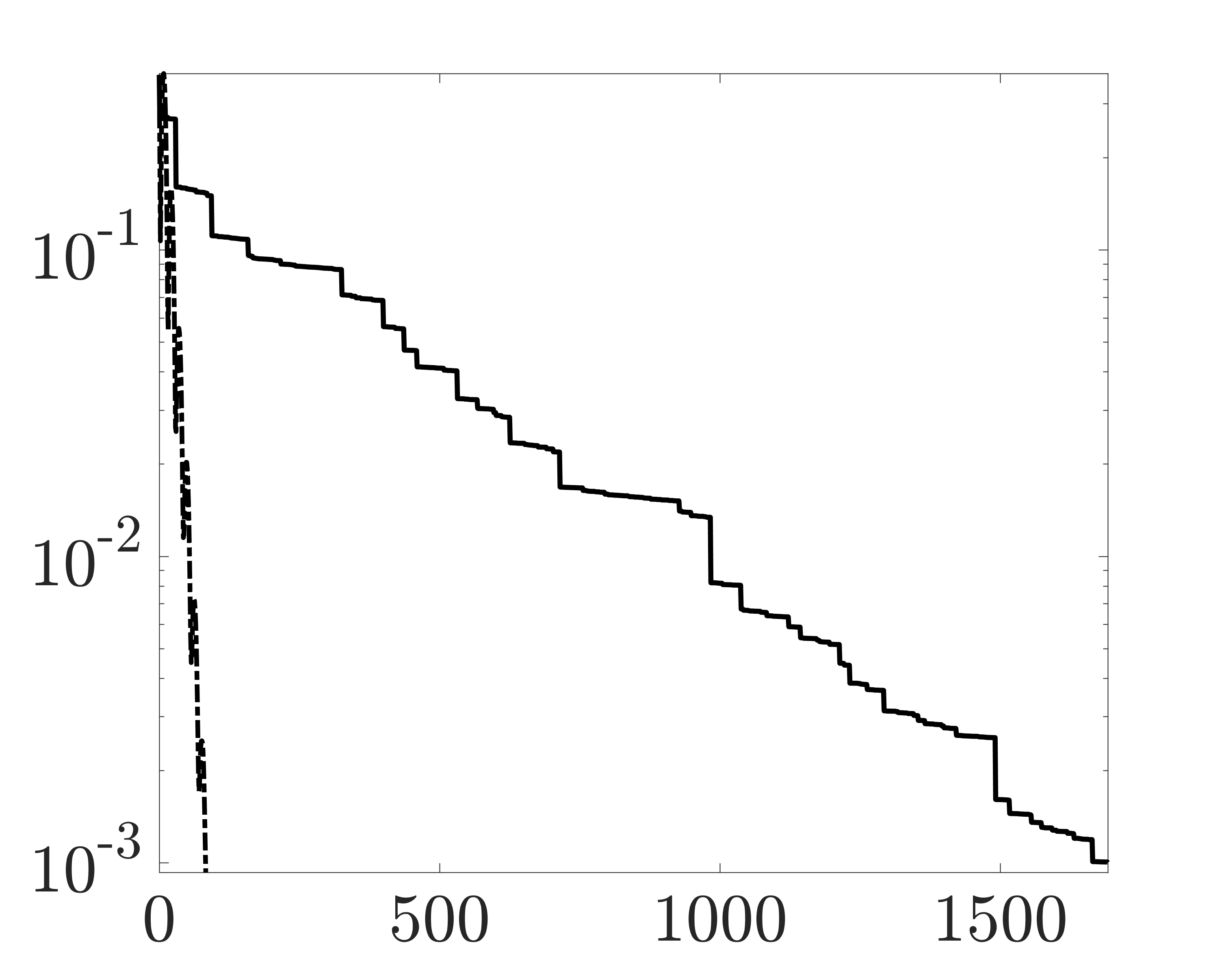}
		\\[-.1cm]
		{\small iteration}
	\end{tabular}
	&
	\hspace{-1.35cm}
	\begin{tabular}{c}
         	\includegraphics[width=.25\textwidth]{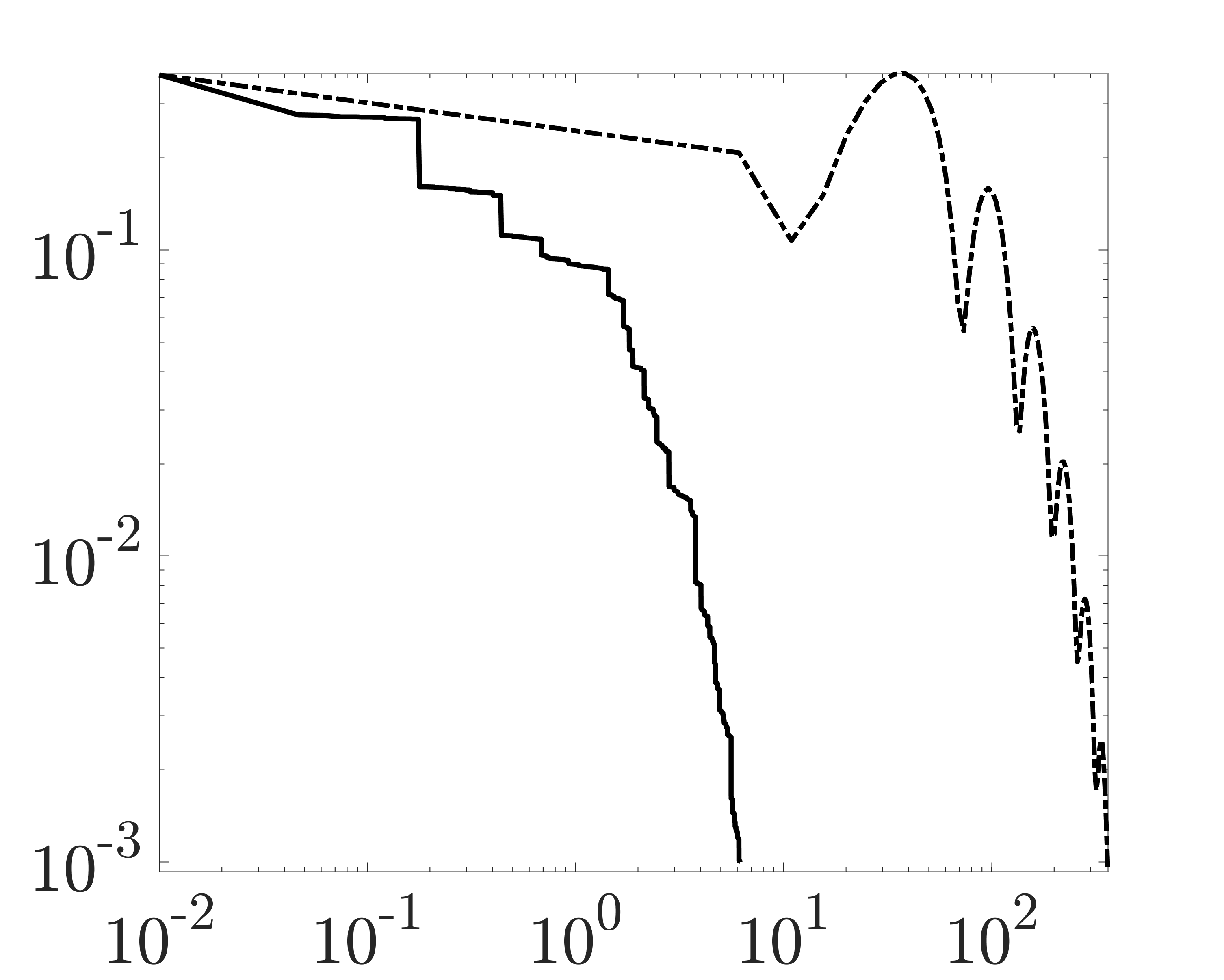}
		\\[-.1cm]
		{\small solve time (sec)}
	\end{tabular}
	\\[-.2cm]
	&
	\hspace{-1cm}
	\subfigure[]{\label{fig.iter-SH}}
	&
	\hspace{-1cm}
	\subfigure[]{\label{fig.time-SH}}
\end{tabular}
\vspace{-.4cm}
\caption{Convergence curves showing performance of PG ($-$) and ADMM ($-\cdot-$) vs. (a) the number of outer iterations; and (b) solve times for the Swift-Hohenberg problem with $n=32$ discretization points and $\gamma=10$. Here, $Y^\star$ is the optimal value for $Y$.
}
\label{fig.convergence-SH}
\end{figure}

For $n=128$ and $256$, ADMM did not converge to desired accuracy in reasonable time. Typically, the ADMM algorithm of~\cite{dhijovluoCDC14} computes low-accuracy solutions quickly but obtaining higher accuracy requires precise solutions to subproblems. The iterative reweighting scheme of Section~\ref{sec.reweighting_polishing} can be used to improve the sparsity patterns that are identified by such low-accuracy solutions. Nonetheless, Fig.~\ref{fig.convergence-SH} shows that even for larger tolerances, PG is faster than ADMM.

\begin{table}
\label{table.SH_CVX_PG_ADMM}
\caption{Comparison of different algorithms (in seconds) for different number of discretization points $n$ and $\gamma=10$.}
\vspace{-.1cm}
\centering
{\footnotesize
\begin{tabular}{c c c c}
\hline
\hline
\\[-.15cm]
$n$ & CVX & PG & ADMM
\\[.1cm]
\hline
\\[-.15cm]
$32$ & $12.39$ & $6.2$ & $362.4$ 
\\[.1cm]
$64$ & $268.11$ & $51.9$ & $4182.6$
\\[.1cm]
$128$ & $8873.3$ & $875.8$ & $-$
\\[.1cm]
$256$ & $-$ & $3872.1$ & $-$
\\[.1cm]
\hline
\hline
\end{tabular}
}
\end{table}

\begin{figure}
\begin{center}
\begin{tabular}{rcrc}
\hspace{-.6cm}
	\begin{tabular}{c}
		\vspace{1.1cm}
		\rotatebox{90}{\small number of actuators}
	\end{tabular}
	&
	\hspace{-.92cm}
	\begin{tabular}{c}
         	\includegraphics[width=.225\textwidth]{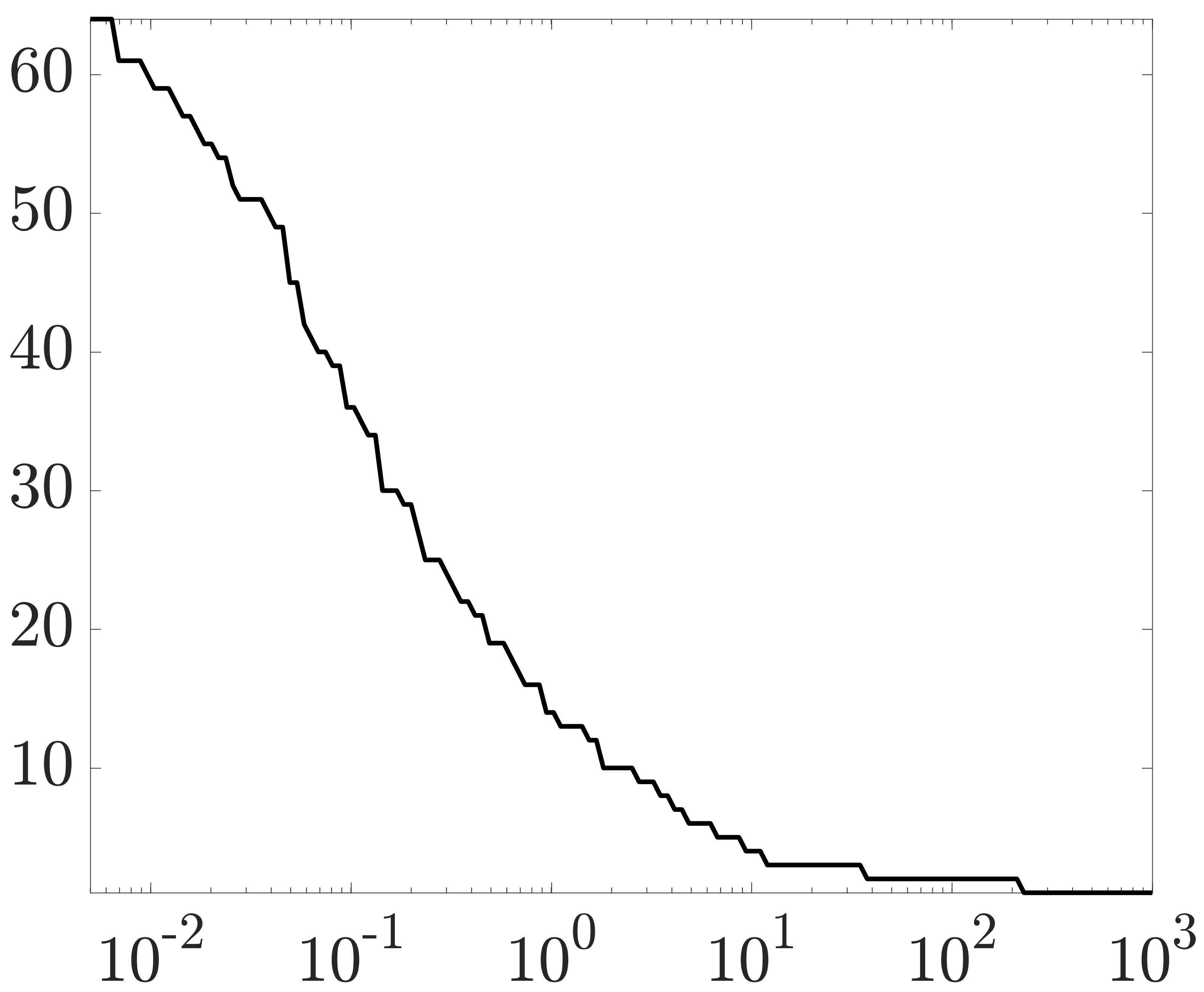}
		\\[-.1cm]
		{\small $\gamma$}
		\\[-.1cm]
		\subfigure[]{\label{fig.nnZ_gamma-SH}}
	\end{tabular}
	&
	\hspace{-1cm}
	\begin{tabular}{c}
		\vspace{1.1cm}
		\rotatebox{90}{\small $(J-J_c)/J_c (\%)$}
	\end{tabular}
	&
	\hspace{-.96cm}
	\begin{tabular}{c}
         	\includegraphics[width=.22\textwidth]{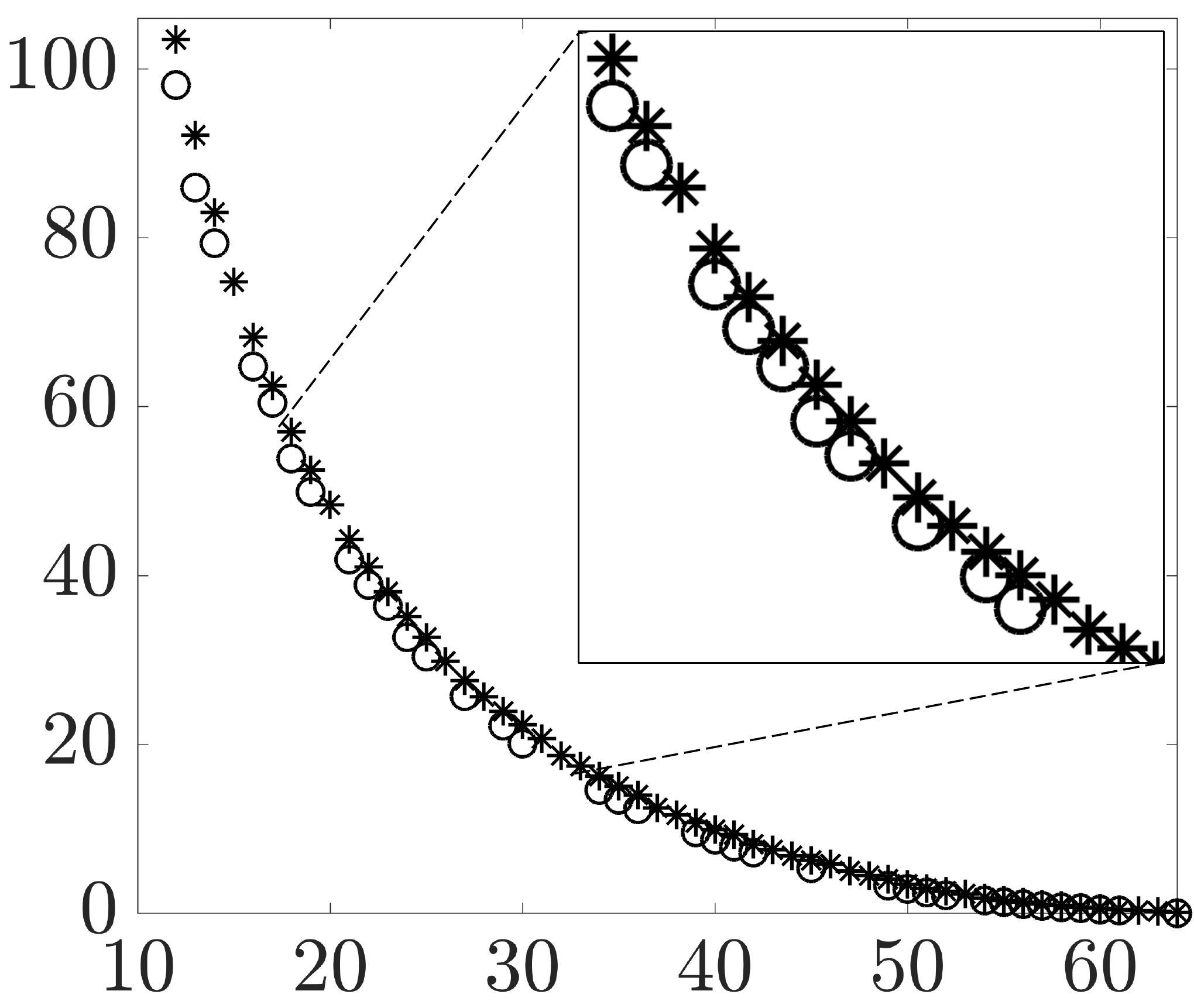}
		\\[-.1cm]
		{\small number of actuators}
		\\[-.1cm]
		\subfigure[]{\label{fig.performance_nnZ-SH}}
	\end{tabular}
	\\[-.4cm]
\end{tabular}
\end{center}
\caption{(a) Number of actuators as a function of the sparsity-promoting parameter $\gamma$; and (b) performance comparison of the optimal feedback controller resulting from the regularized actuator selection problem ($\Circle$) and from the greedy algorithm ($*$) for the Swift-Hohenberg problem with $n=64$.}
\label{fig.performance-SH}
\end{figure}

As $\gamma$ increases in {Problem~\ref{prob2}}, more and more actuators are dropped and the performance degrades monotonically. For $n=64$, Fig.~\ref{fig.nnZ_gamma-SH} shows the number of retained actuators as a function of $\gamma$ and Fig.~\ref{fig.performance_nnZ-SH} shows the percentage of performance degradation as a function of the number of retained actuators. Figure~\ref{fig.performance_nnZ-SH} also illustrates that for various numbers of retained actuators, the solution to convex optimization problem~\eqref{eq.MCC-1} with $\delta=1$ consistently yields performance degradation that is no larger than the performance degradation of a greedy algorithm (that drops actuators based on their contribution to the $\cH_2$ performance index). For example, the greedy algorithm leads to $24.6\%$ performance degradation when $30$ actuators are retained whereas our approach yields $20\%$ performance degradation for the same number of actuators. This greedy heuristic is summarized in Algorithm~\ref{alg.greedy}, where $S$ is the set of actuators and $f(S)$ denotes the performance index resulting from the actuators within the set $S$. When the individual subproblems for choosing fixed numbers of actuators can be executed rapidly, greedy algorithms provide a viable alternative. There has also been recent effort to prove the optimality of such algorithms for certain classes of problems~\cite{sum16}. However, in our example, the greedy algorithm does not always provide the optimal set of actuators with respect to the $\cH_2$ performance index. Relative to the convex formulation, similar greedy techniques yield suboptimal sensor selection for a flexible aircraft wing~\cite[Section 5.2]{jovdhiEJC16}. 

\begin{figure}
\begin{tabular}{cc}
	\hspace{-.6cm}
         \includegraphics[height =0.21\textwidth, width=0.26\textwidth]{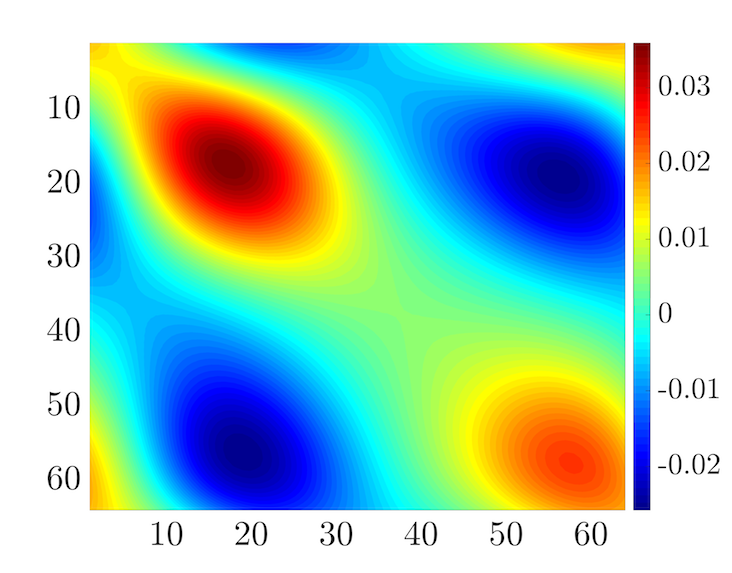}
	&
	\hspace{-.56cm}
         \includegraphics[height =0.21\textwidth, width=0.26\textwidth]{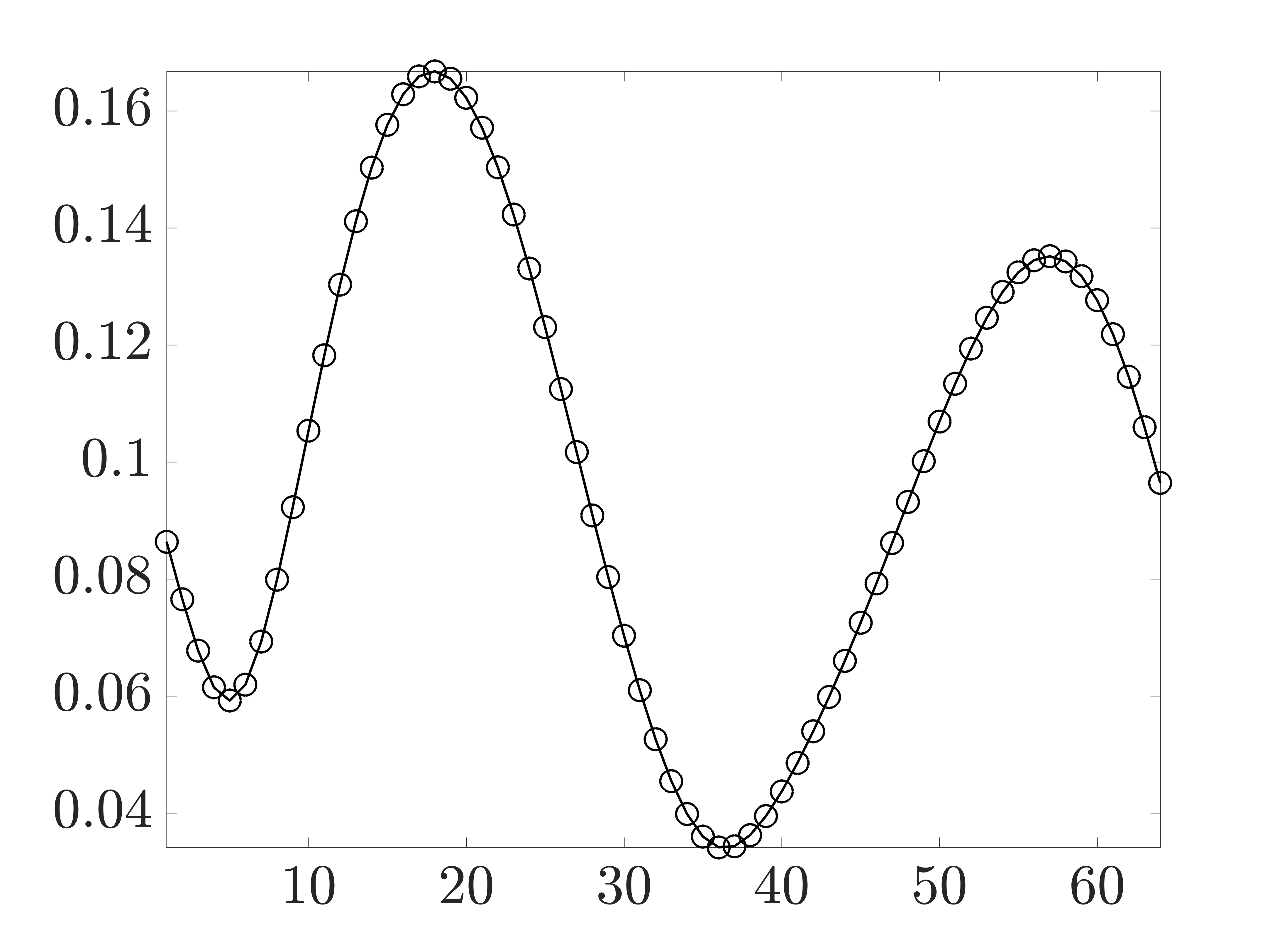}
         \\[-.2cm]
	\hspace{-.6cm}
	&
	\hspace{-.3cm}
	{\small row number}
	\\[-.2cm]
	\hspace{-.6cm}
	\subfigure[]{\label{fig.Klqr}}
	&
	\hspace{-.56cm}
	\subfigure[]{\label{fig.Klqr_norm}}
	\\[-.1cm]
	\hspace{-.6cm}
         \includegraphics[height =0.21\textwidth, width=0.26\textwidth]{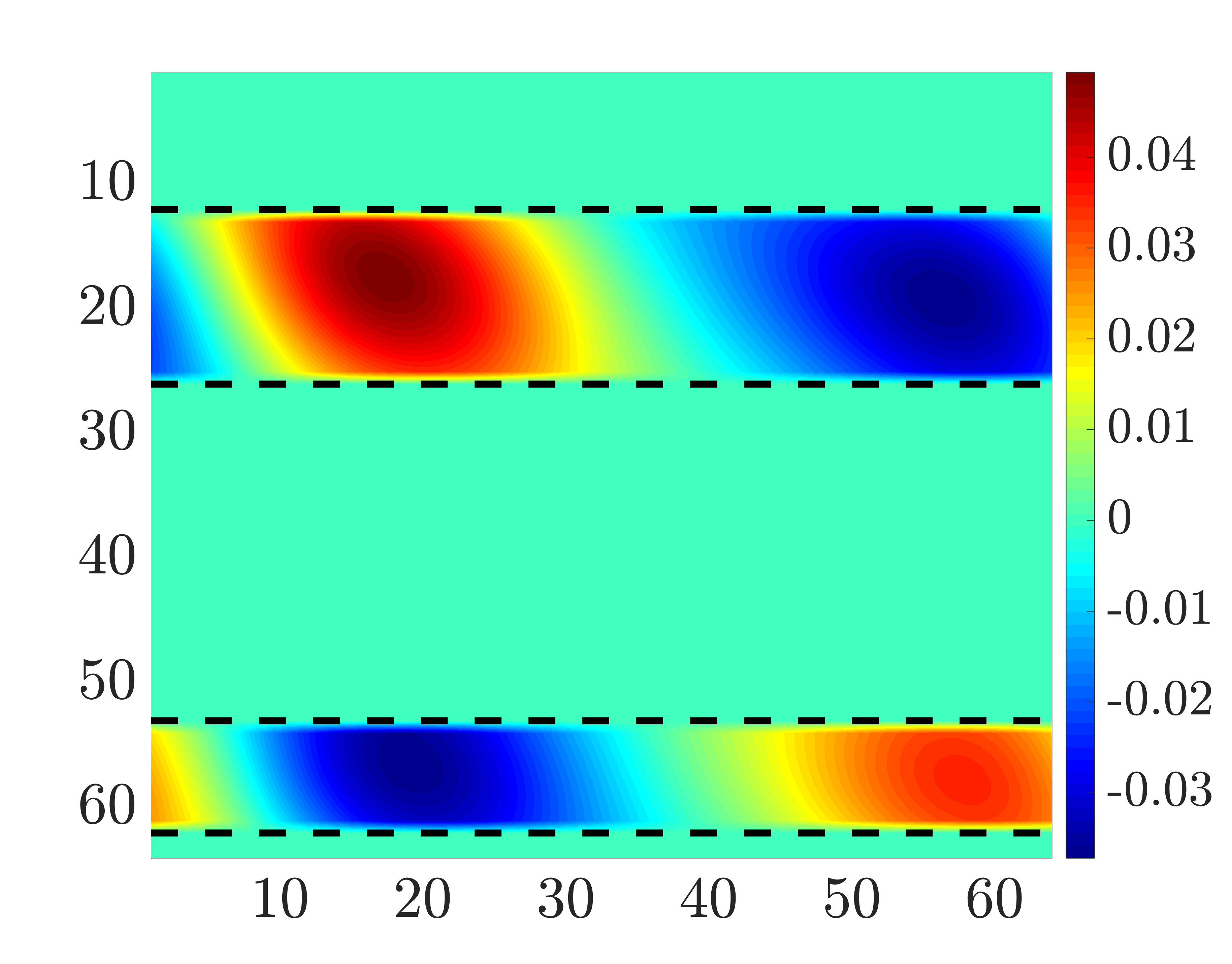}
	&
	\hspace{-.56cm}
         \includegraphics[height =0.21\textwidth, width=0.26\textwidth]{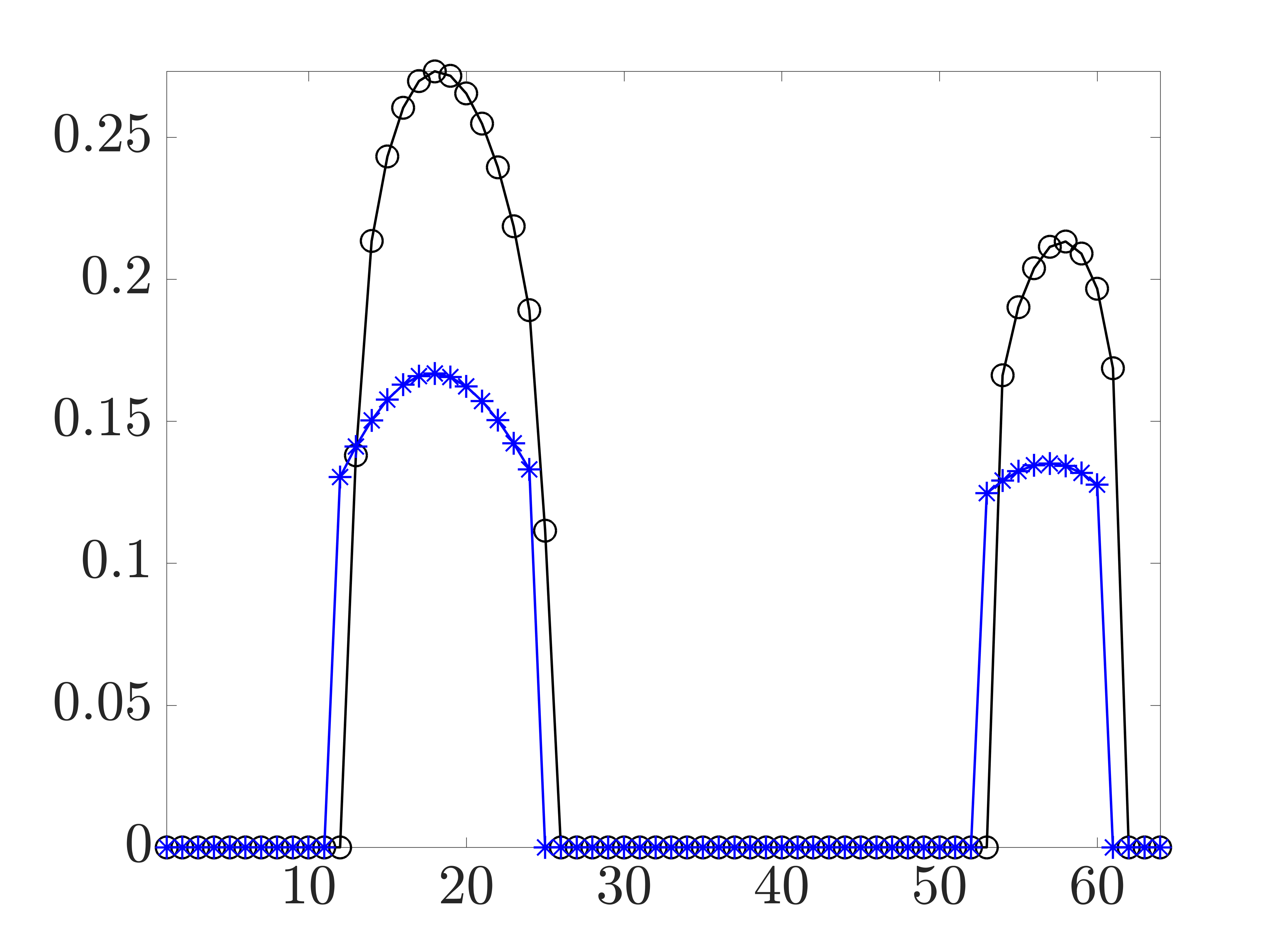}
          \\[-.2cm]
	\hspace{-.6cm}
	&
	\hspace{-.3cm}
	{\small row number}
	\\[-.2cm]
	\hspace{-.6cm}
	\subfigure[]{\label{fig.Ksp_gamma0p4}}
	&
	\hspace{-.56cm}
	\subfigure[]{\label{fig.Ksp_norm_gamma0p4}}
\end{tabular}
\caption{(a) Optimal centralized feedback gain matrix and (b) its row-norms corresponding to the Swift-Hohenberg dynamics~\eqref{eq.SH.dynamics} with $n=64$. (c) The optimal feedback gain matrix and (d) its row-norms {($\circ$)} resulting from solving {Problem~\ref{prob2}} with $\delta=1$ and $\gamma=0.4$ in which the rows between the dashed lines have been retained and polished via optimization. {The result of truncating the centralized feedback gain matrix based on its row-norms is shown using blue $*$ symbols.}}
\label{fig.Klqr_Ksp}
\end{figure}

The absence of the sparsity promoting regularizer in {Problem~\ref{prob2}} leads to the optimal centralized controller which can be obtained from the solution to the algebraic Riccati equation. For $n = 64$, Figs.~\ref{fig.Klqr} and~\ref{fig.Klqr_norm} show this centralized feedback gain and the two norms of its rows, respectively. For {$\gamma=0.4$, $21$} of $64$ possible actuators are retained and the corresponding optimal feedback gain matrix and row norms are shown in Figs.~\ref{fig.Ksp_gamma0p4} and~\ref{fig.Ksp_norm_gamma0p4}. {Figure~\ref{fig.Ksp_norm_gamma0p4} also shows that a truncation of the centralized feedback gain matrix based on its row-norms (marked by blue $*$ symbols) yields a different subset of actuators than the solution to Problem~\ref{prob2}.}

\begin{algorithm}
\caption{A greedy heuristic for actuator selection}
\label{alg.greedy}
\begin{algorithmic}
\STATE \textbf{input:} $A$, $B$, $V$, $Q$ $R$.
\vspace*{0.05cm}
\STATE \textbf{initialize:} ${\Pi} \leftarrow \{1,\dots,m\}$.
\vspace*{0.1cm}
\STATE \textbf{while:} $| {\Pi} | > 0$ and $f(S) < \infty$
\\[.1cm]
\quad\quad\quad $e^* = \argmin\limits_{e \, \in \, {\Pi}} f({\Pi}) - f({\Pi} \backslash\{e\})$
\\[.1cm]
\quad\quad\quad ${\Pi} \leftarrow {\Pi} \backslash \{e\}$
\vspace*{0.05cm}
\STATE \textbf{endwhile}
\STATE \textbf{output:} the set of actuators represented by the set ${\Pi}$.
\end{algorithmic}
\end{algorithm}

	\vspace*{-3ex}
\subsection{Covariance completion}
\label{sec.example-cc}

We provide an example to demonstrate the utility of our {approach} for the purpose of completing partially available second-order statistics of a three-dimensional channel flow. In an incompressible channel-flow, the dynamics of infinitesimal fluctuations around the parabolic mean velocity profile, $\bar{\bu} = [\,U(x_2)\,~0\,~0\,]^T$ with $U(x_2) = 1-x_2^2$, are governed by the Navier-Stokes equations linearized around $\bar{\bu}$. The streamwise, wall-normal, and spanwise coordinates are represented by $x_1$, $x_2$, and $x_3$, respectively; see Fig.~\ref{fig.channel} for geometry. Finite dimensional approximation via application of the Fourier transform in horizontal dimensions ($x_1$ and $x_3$) and spatial discretization of the wall-normal dimension ($x_2$) using $N$ collocation points, yields the state-space representation
\begin{subequations}
	\label{eq.filter-system}
\begin{align}
\begin{array}{rcl}
    \dot{\bpsi} (\bk, t)
    & \!\! = \!\! &
    A (\bk)
    \,
    \bpsi (\bk, t)
    \; + \;
    \xi (\bk, t)
    \\[0.15cm]
    \bv (\bk,\, t)
    & \!\! = \!\! &
    C (\bk)
    \,
    \bpsi (\bk, t).
\end{array}
    \label{eq.channel-lnse1}
\end{align}
Here, $\bpsi = [\,v_2^T\,~\eta^T\,]^T \in \bbC^{2N}$ is the state of the linearized model, $v_2$ and $\eta = \partial_{x_3} v_1 - \partial_{x_1} v_3$ are the normal velocity and vorticity, the output $\bv = [\,v_1^T\,~v_2^T\,~v_3^T\,]^T \in \bbC^{3N}$ denotes the fluctuating velocity vector, $\xi$ is a stochastic forcing disturbance, $\bk = [\,k_{1}\,\;k_{3}\,]^T$ denotes the vector of horizontal wavenumbers, and the input matrix is the identity $I_{2N \times 2N}$. The dynamical matrix $A \in \bbC^{2N \times 2N}$ and output matrix $C \in \bbC^{3N \times 2N}$ are described in~\cite{jovbamJFM05}.

\begin{figure}
\centering
\begin{tabular}{cc}
	\begin{tabular}{c}
         \includegraphics[width=0.275\textwidth]{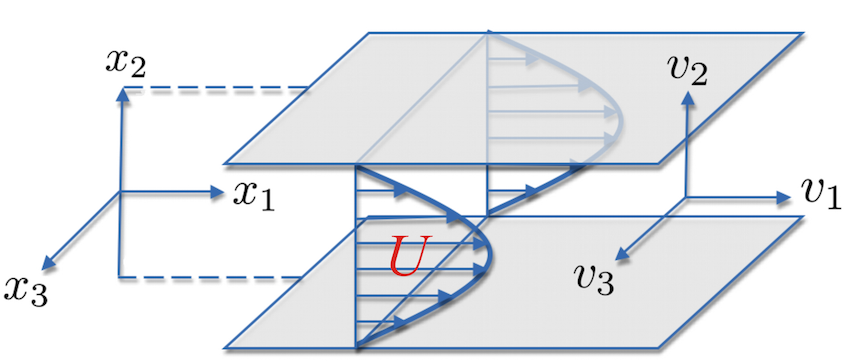}
         \end{tabular}
	&
	\hspace*{-0.5cm}
	\begin{tabular}{c}
	\includegraphics[width=0.125\textwidth]{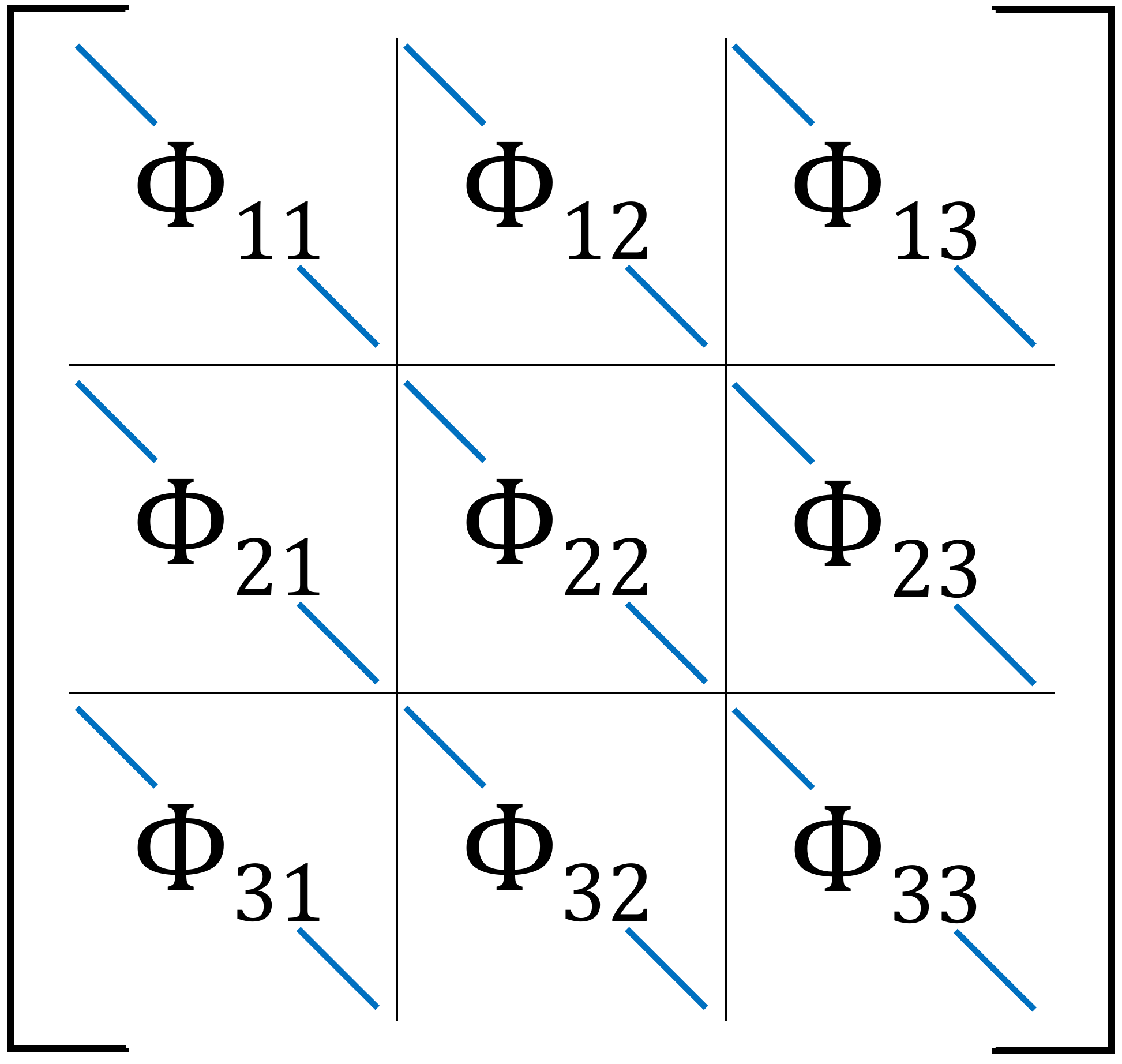}
	\end{tabular}
	\\[-.2cm]
	\subfigure[]{\label{fig.channel}}
	&
	\hspace{-0.5cm}
	\subfigure[]{\label{fig.output-covariance}}
\end{tabular}
\caption{(a) Geometry of a three-dimensional pressure-driven channel flow. (b) Structure of the matrix $\Phi= \lim_{t \rightarrow \infty} \mathbf{E} \left( \bv(t) \bv^*(t) \right)$, where $\Phi_{ij}$ denotes the cross-correlation matrix of components $v_i$ and $v_j$ of the velocity vector $\bv$ across the discretization points in the wall-normal direction. Available diagonal entries of the blocks in the velocity covariance matrix $\Phi$ determine correlations at the same discretization~point.}
\label{fig.channel-covariance}
\end{figure}

We assume that the stochastic disturbance $\xi$ is generated by a low-pass filter with state-space representation
\begin{align}
	\dot{\xi}(\bk,t)
	\;=\;
	- {\xi}(\bk,t) \,+\, w(t)
	\label{eq.filter-1}
\end{align}
\end{subequations}
where $w$ denotes a zero-mean {white process with identity covariance matrix.} The steady-state covariance of system~\eqref{eq.filter-system} can be obtained as {the} solution to the Lyapunov equation
\begin{align*}
	\ba{c}
	\tilde{A}\, \Sigma
	\;+\;
	\Sigma\, \tilde{A}^*
	\;+\;
	\tilde{B}\, \tilde{B}^*
	\;=\;
	0
	\\[.15cm]
	\tilde{A}
	\;=\;
	\tbt{A}{\,~~I}{O}{-I},
	~
	\tilde{B}
	\;=\;
	\tbo{0}{I},
	~
	\Sigma
	\;=\;
	\tbt{\Sigma_{11}}{\Sigma_{12}}{\Sigma_{12}^*}{\Sigma_{22}}.
	\ea
\end{align*}
For any $\bk$, the matrix $\Sigma_{11} = \lim_{t \rightarrow \infty} \mathbf{E} \left( \bpsi(t) \bpsi^*(t) \right)$ denotes the steady-state covariance of system~\eqref{eq.channel-lnse1} and is related to the steady-state covariance matrix of the output $\bv$ via
	$
	\Phi (\bk)
	= 
	C (\bk) \Sigma_{11}(\bk) C^* (\bk).
	$
Figure~\ref{fig.output-covariance} shows the structure of the output covariance matrix~$\Phi$.

In this example, we assume that all one-point velocity correlations, i.e., the diagonal entries of all submatrices $\Phi_{ij}$ in Fig.~\ref{fig.output-covariance}, are known. Owing to experimental and computational limitations, one-point correlations are easier to measure and compute than two-point spatial correlations~\cite{zargeojovARC20}. While the colored-in-time input process $\xi$ enters across all channels, not all input channels equally impact the state statistics $\Sigma_{11}$ as the input to state gain differs across different inputs. Herein, we seek a minimal set of input channels with dominant contribution that can lead to a parsimonious perturbation $A-BK$ of the system dynamics. The identified structure represents important feedback mechanisms that are responsible for generating the available statistics when the system is driven by white noise $d$. Finally, we note that due to the parameterization of system dynamics~\eqref{eq.filter-system} over wavenumbers $\bk$, modification $BK$ also depends on $\bk$.

\begin{figure}
\begin{tabular}{rcc}
	\hspace{-.6cm}
	\begin{tabular}{c}
	\vspace{.4cm}
	\rotatebox{90}{\footnotesize $\norm{Y^k-Y^\star}_F/\norm{Y^\star}_F$}
	\end{tabular}
	&
	\hspace{-1cm}
	\begin{tabular}{c}
         	\includegraphics[width=.25\textwidth]{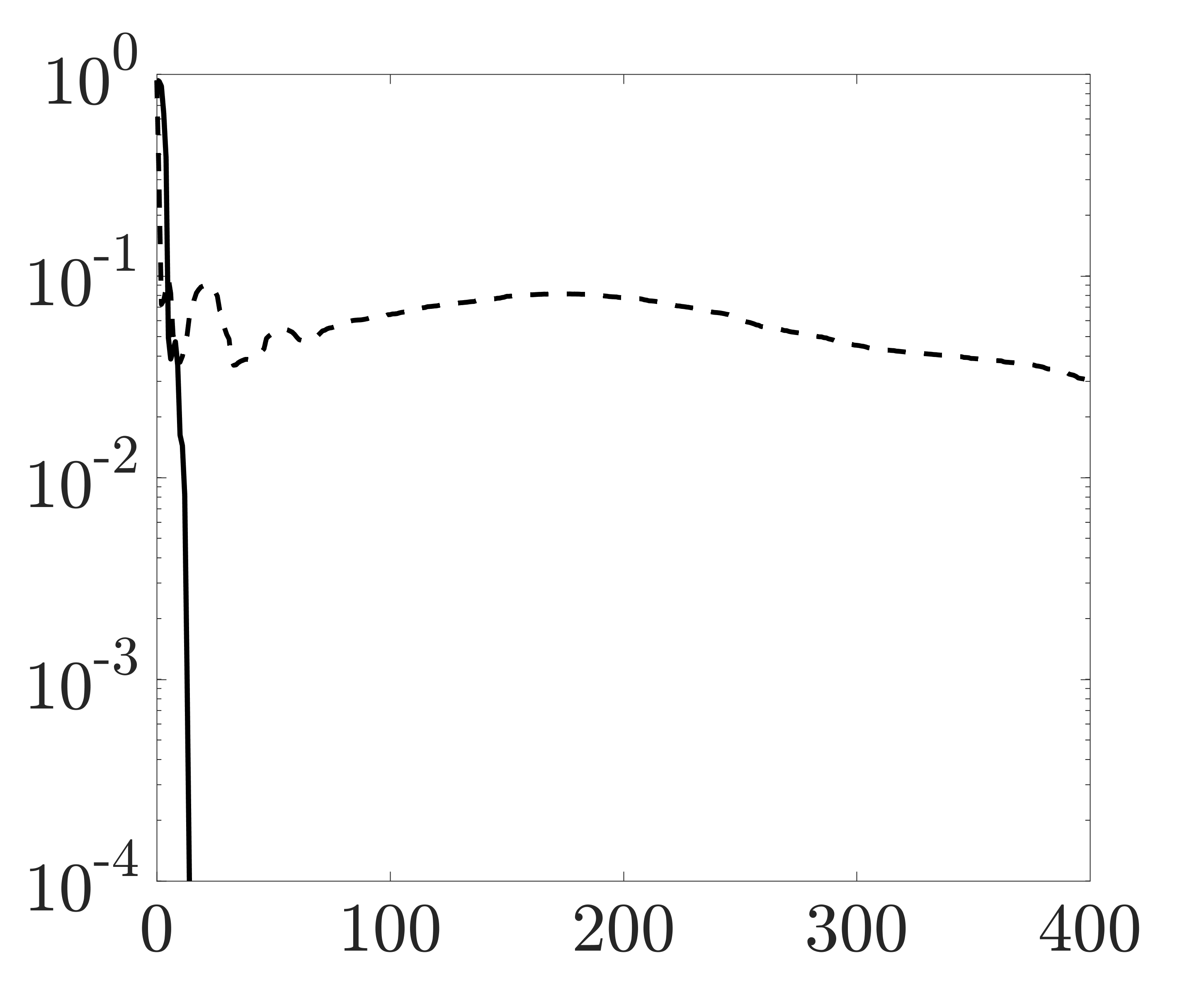}
		\\[-.1cm]
		{\small iteration}
	\end{tabular}
	&
	\hspace{-1.15cm}
	\begin{tabular}{c}
         	\includegraphics[width=.25\textwidth]{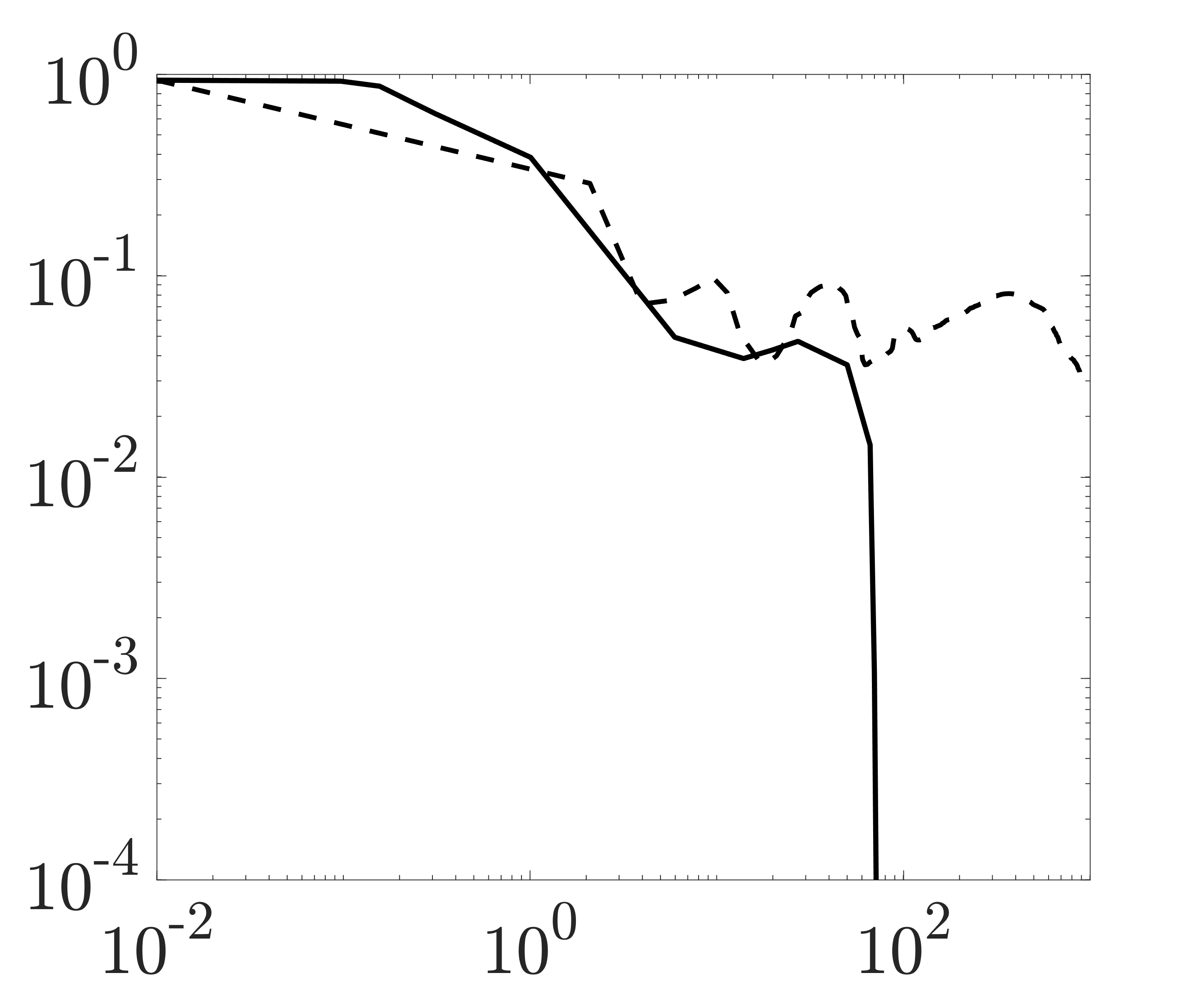}
		\\[-.1cm]
		{\small solve time (sec)}
	\end{tabular}
	\\[-.2cm]
	&
	\hspace{-1cm}
	\subfigure[]{\label{fig.iter}}
	&
	\hspace{-1cm}
	\subfigure[]{\label{fig.time}}
\end{tabular}
\vspace{-.4cm}
\caption{Convergence curves showing performance of MM ($-$) and ADMM ($- -$) versus (a) the number of outer iterations; and (b) solve times for $N=31$ collocation points in the normal direction $x_2$ and $\gamma=10$. Here, $Y^\star$ is the optimal value for $Y$.
}
\label{fig.convergence}
\vspace{-1.5cm}
\end{figure}

\begin{figure}
\begin{tabular}{rcrc}
	\hspace{-.6cm}
	\begin{tabular}{c}
	\vspace{.4cm}
	\rotatebox{90}{\small $\Delta_p/\norm{G}_F$}
	\end{tabular}
	&
	\hspace{-.9cm}
	\begin{tabular}{c}
         	\includegraphics[width=.235\textwidth]{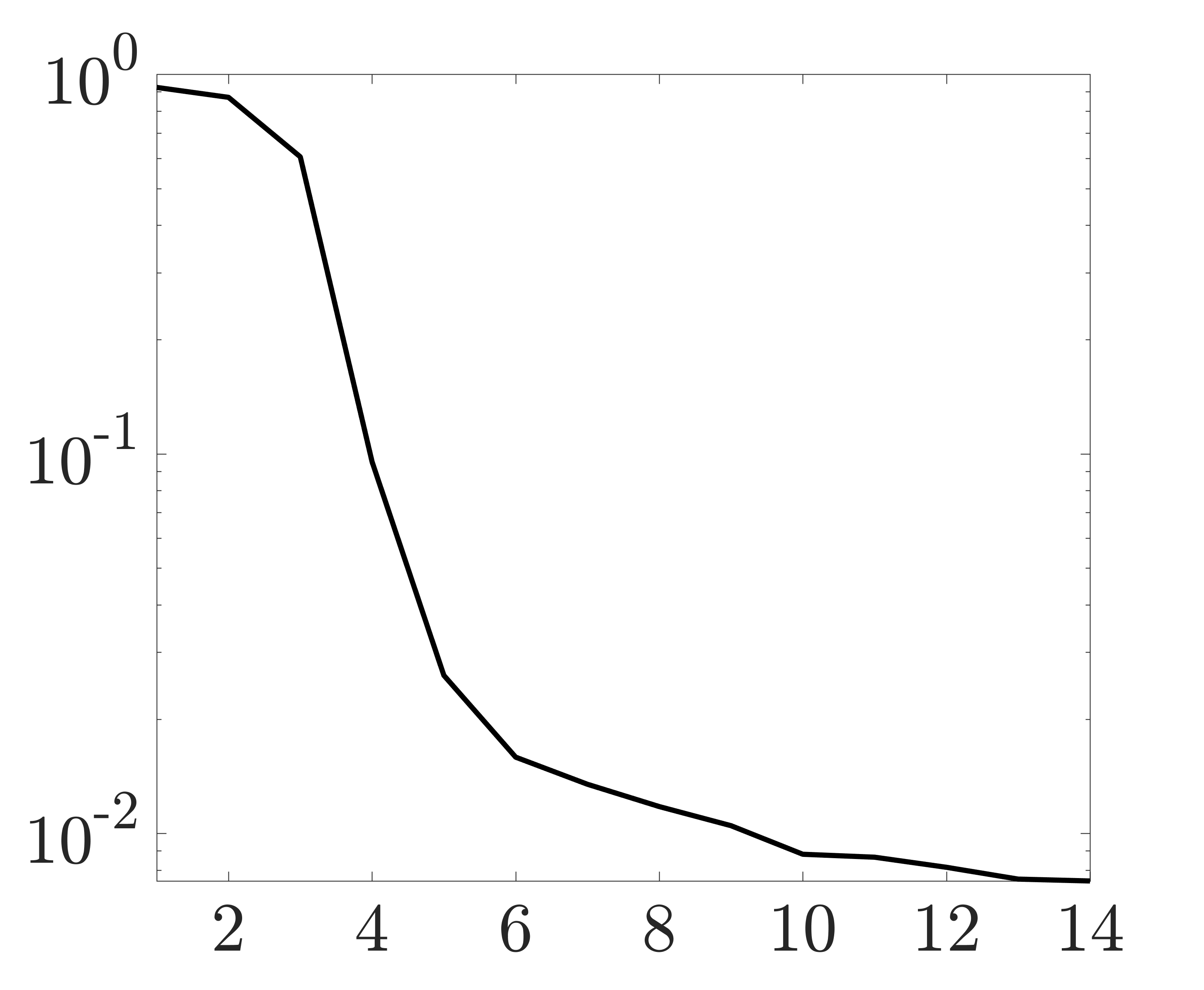}
		\\[-.1cm]
		{\small iteration}
	\end{tabular}
	\hspace{-.45cm}
	&
	\hspace{-2cm}
	\begin{tabular}{c}
	\vspace{.4cm}
	\rotatebox{90}{\small $\Delta_d$}
	\end{tabular}
	&
	\hspace{-.9cm}
	\begin{tabular}{c}
         	\includegraphics[width=.235\textwidth]{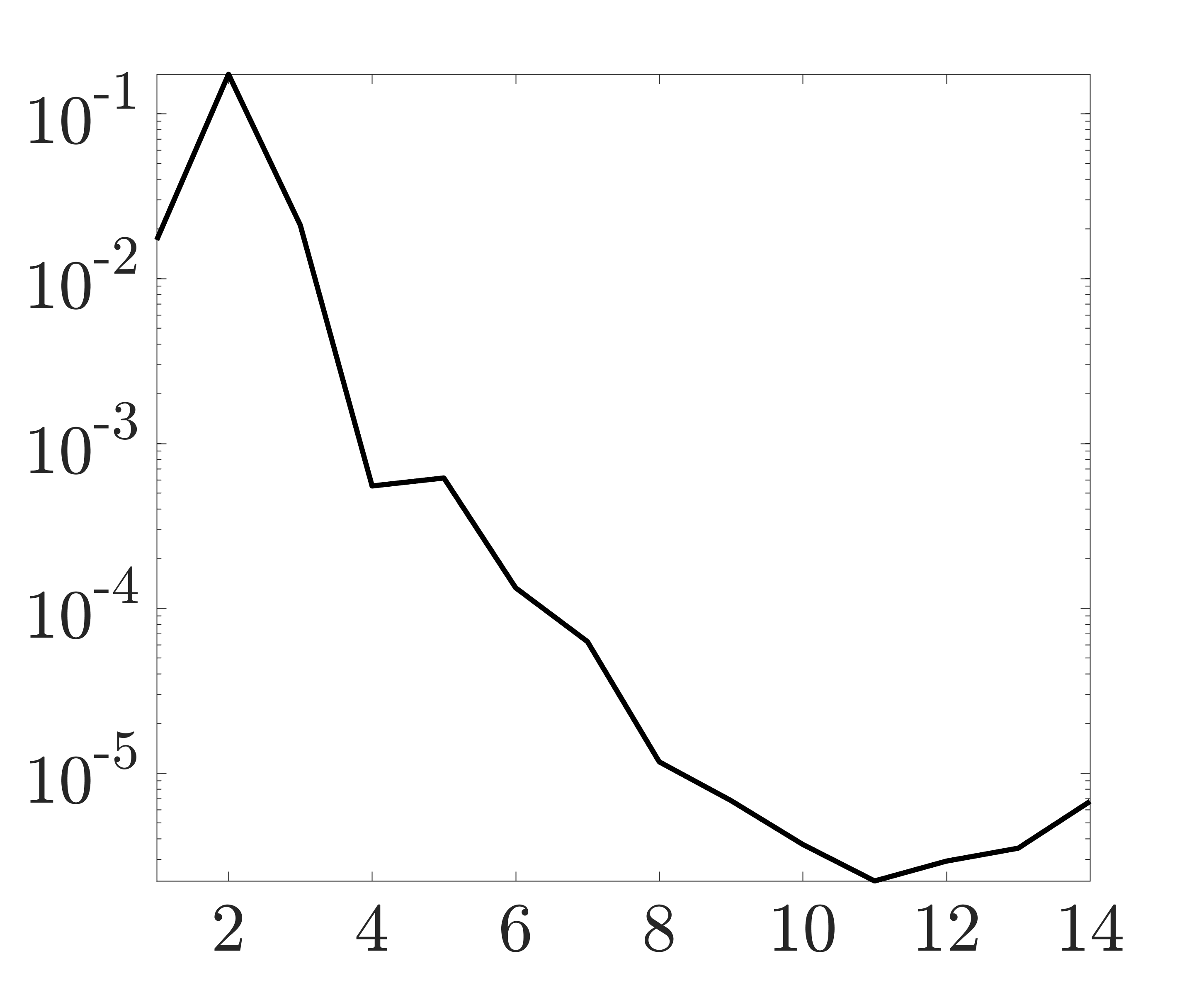}
		\\[-.1cm]
		{\small iteration}
	\end{tabular}
	\\[-.2cm]
	&
	\hspace{-1cm}
	\subfigure[]{\label{fig.iter-primalresidual}}
	&&
	\hspace{-1cm}
	\subfigure[]{\label{fig.iter-dualresidual}}
\end{tabular}
\vspace{-.4cm}
\caption{Performance of MM for the fluids example with $N=31$ collocation points in the normal direction $x_2$ and $\gamma=10$. (a) normalized primal residual; and (b) dual residual based on~\eqref{eq.MM-termination}.}
\label{fig.performance-fluids-primal-dual-residual}
\end{figure}

Computational experiments are conducted for a flow with Reynolds number $10^3$, the wavenumber pair $(k_1,k_3) = (0,1)$, for various number of collocation points $N$ in the wall-normal direction (state dimension $n=2 N$), {$R = I$, $Q=0$}, and for various values of the regularization parameter $\gamma$. {Moreover, we assume that system~\eqref{eq.LTIsys} is driven by white process $d$ with covariance $V=I$. We initialize Algorithm~\ref{alg.MM}} with the optimal centralized controller, $Y^0 \DefinedAs K_c X_c$. Our MM algorithm is compared against SDPT3 and ADMM where CVX is used to call SDPT3. When CVX can compute the optimal solution of Problem~\ref{prob2}, for each method, iterations are run until the solutions are within $5\%$ of the CVX solution. For larger problems, iterations are run until the primal and dual residuals satisfy certain tolerances; $\eps_p$, $\eps_d = 10^{-2}$. For $\gamma = 10$, Table~\ref{table.1} compares various methods based on run times (sec). For $N = 51$ and $101$, CVX failed to converge and ADMM did not converge in a reasonable time. Clearly, MM outperforms ADMM. This can also be deduced from Fig.~\ref{fig.convergence}, which shows convergence curves for $14$ steps of MM and $500$ steps of ADMM for $N=31$ and $\gamma=10$. For this example, Fig.~\ref{fig.performance-fluids-primal-dual-residual} shows the convergence of MM based on the normalized primal residual $\Delta_p/\norm{G}_F$ and the dual residual $\Delta_d$ in~\eqref{eq.MM-termination}.

\begin{table}
\label{table.1}
\caption{Comparison of different algorithms (in seconds) for different number of discretization points $N$ and $\gamma=10$.}
\vspace{-.1cm}
\centering
{\footnotesize
\begin{tabular}{c c c c}
\hline
\hline
\\[-0.15cm]
N \!\! & \!\! CVX \!\! & \!\! MM \!\! & \!\! ADMM
\\[.1cm]
\hline
\\[-.15cm]
 $11$ \!\! & \!\! $9.3$ \!\! & \!\! $0.19$ \!\! & \!\! $3.10$
\\[.1cm]
$21$ \!\! & \!\! $97.67$ \!\! & \!\! $5.6$ \!\! & \!\! $113.4$
\\[.1cm]
$31$ \!\! & \!\! $900$ \!\! & \!\! $7.19$ \!\! & \!\! $574.44$
\\[.1cm]
$51$ \!\! & \!\! $-$ \!\! & \!\! $34.76$ \!\! & $-$
\\[.1cm]
$101$ \!\! & \!\! $-$ \!\! & \!\! $146.51$ \!\! & $-$
\\[.1cm]
\hline
\hline
\end{tabular}
}
\end{table}

\begin{figure}
\centering
\begin{tabular}{rc}
	\begin{tabular}{c}
		\vspace{.4cm}
		\rotatebox{90}{number of inputs}
	\end{tabular}
	&
	\hspace{-1cm}
	\begin{tabular}{c}
         	\includegraphics[width=.25\textwidth]{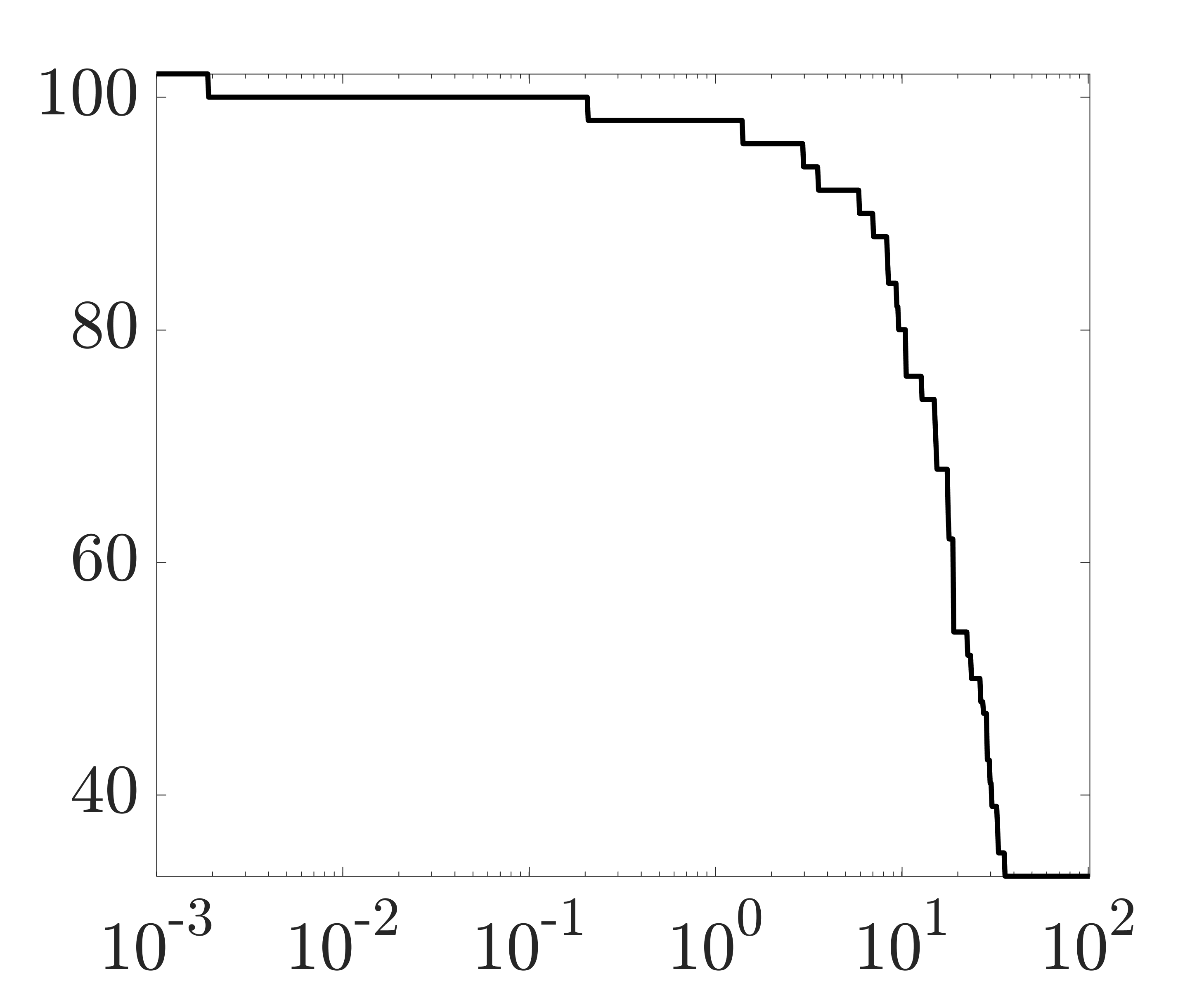}
		\\[-.2cm]
		$\gamma$
	\end{tabular}
\end{tabular}
\caption{The $\gamma$-dependence of the number of input channels that are retained after solving problem~\eqref{eq.MCC-1} for the channel flow problem with $m=102$ inputs.}
\label{fig.gamma-inputs-CC}
\end{figure}

\begin{figure}
\centering
\begin{tabular}{ccccc}
	\begin{tabular}{c}
	\begin{tabular}{c}
		\vspace{1.5cm}
		$\;u_1$
	\end{tabular}
	\\[1.4cm]
	\begin{tabular}{c}
		\vspace{.8cm}
		$u_2$
	\end{tabular}
	\end{tabular}
	&
	\hspace{-.8cm}
	\begin{tabular}{c}
         	\includegraphics[width=1cm,height=5cm]{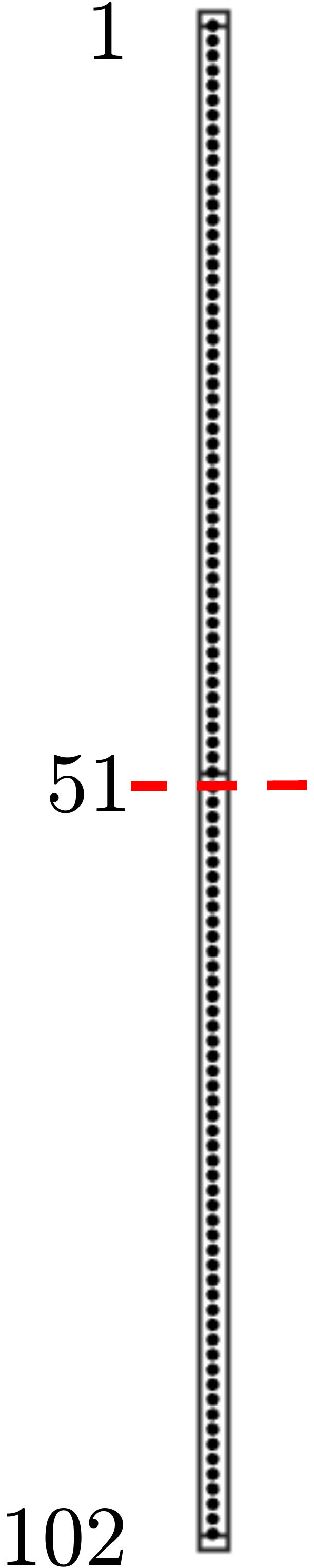}
		\\
		\subfigure[]{\label{fig.input-gamma0}}
	\end{tabular}
	&
	\hspace{-.5cm}
	\begin{tabular}{c}
         	\includegraphics[width=0.9cm,height=5.0cm]{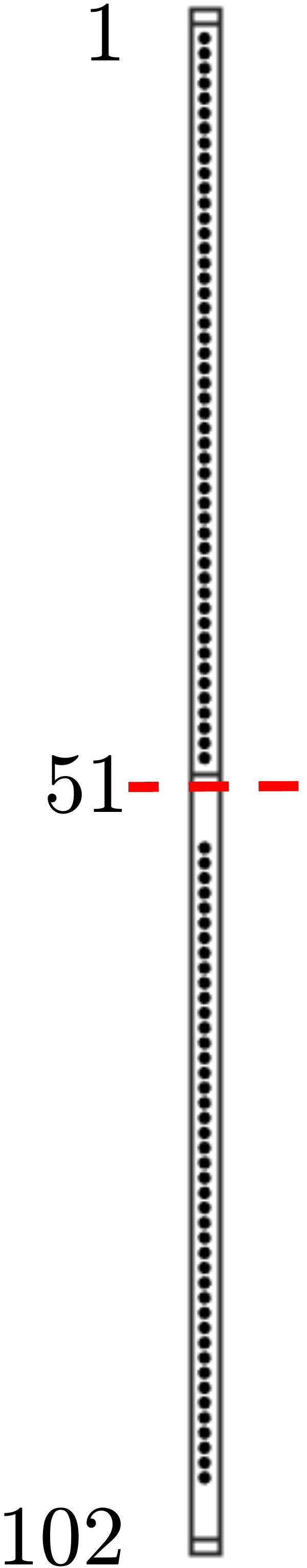}
		\\
		\subfigure[]{\label{fig.input-gamma1e-1}}
	\end{tabular}
	&
	\hspace{-.5cm}
	\begin{tabular}{c}
         	\includegraphics[width=1cm,height=5cm]{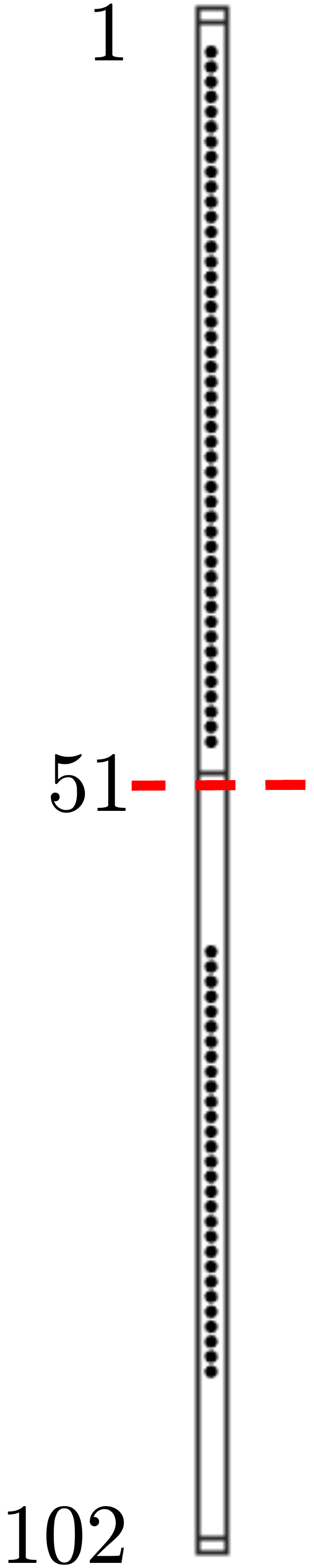}
		\\
		\subfigure[]{\label{fig.input-gamma10}}
	\end{tabular}
	&
	\hspace{-.5cm}
	\begin{tabular}{c}
         	\includegraphics[width=1cm,height=5cm]{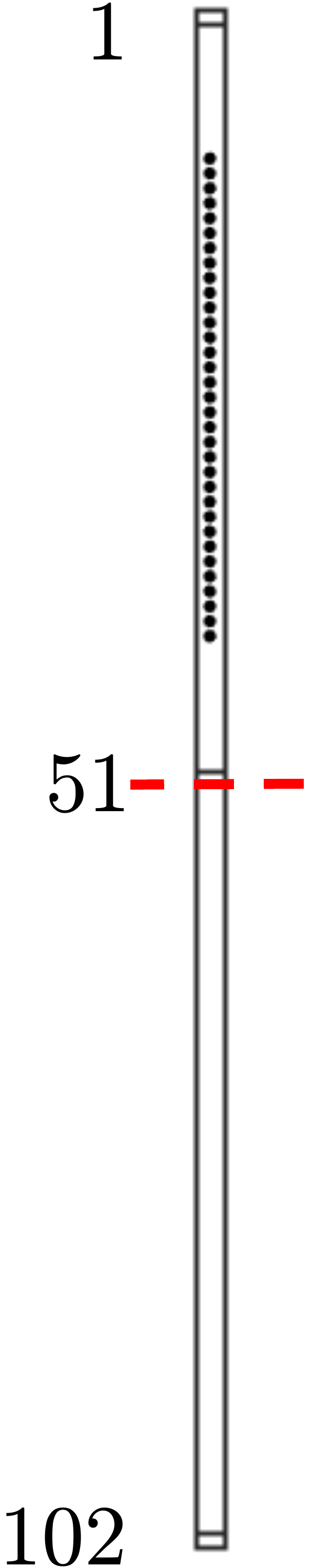}
		\\
		\subfigure[]{\label{fig.input-gamma100}}
	\end{tabular}
\end{tabular}
\vspace{-.2cm}
\caption{{Active input channels in $u \in \bbC^{102}$ (black dots) corresponding to row-sparsity of $Y^*$ in problem~\eqref{eq.MCC-1} with (a) $\gamma=0$, (b) $\gamma=0.1$, (c) $\gamma=10$, and (d) $\gamma=100$, for the channel flow problem.}}
\label{fig.retained-inputs-CC}
\end{figure}

We now focus on $N=51$ collocation points and solve Problem~\ref{prob2} for various values of $\gamma$. Since {$B = I$}, the number of inputs {$u$} in this case is $m=102$. Figure~\ref{fig.gamma-inputs-CC} shows the $\gamma$-dependence of the number of retained input channels that result from solving {Problem~\ref{prob2}}. As $\gamma$ increases, more and more input channels are dropped. A feature of our framework is that the solution $Y^\star$ determines which inputs in $u$ play a role in matching the available statistics in a way that is consistent with the underlying dynamics. Figure~\ref{fig.retained-inputs-CC} shows the input channels that are retained via optimization for different values of $\gamma$. This figure illustrates the dominant role of input channels that enter the dynamics of normal velocity {$v_2$} and away from the boundaries of the channel. In favor of brevity, we do not expand on the physical interpretations of such findings.

\begin{figure}
		\begin{tabular}{cccc}
	\hspace{-.8cm} \subfigure[]{\label{fig.Siguu}}
	& 
	&
	\hspace{-1.4cm} \subfigure[]{\label{fig.Xuu_gammaopt}}
	&
	\\[-.4cm]
	\begin{tabular}{c}
		\vspace{.4cm}
		\hspace{-.6cm}
		\small{$x_2$}
	\end{tabular}
	&
	\hspace{-0.9cm}
	\begin{tabular}{c}
		\includegraphics[width=4cm]{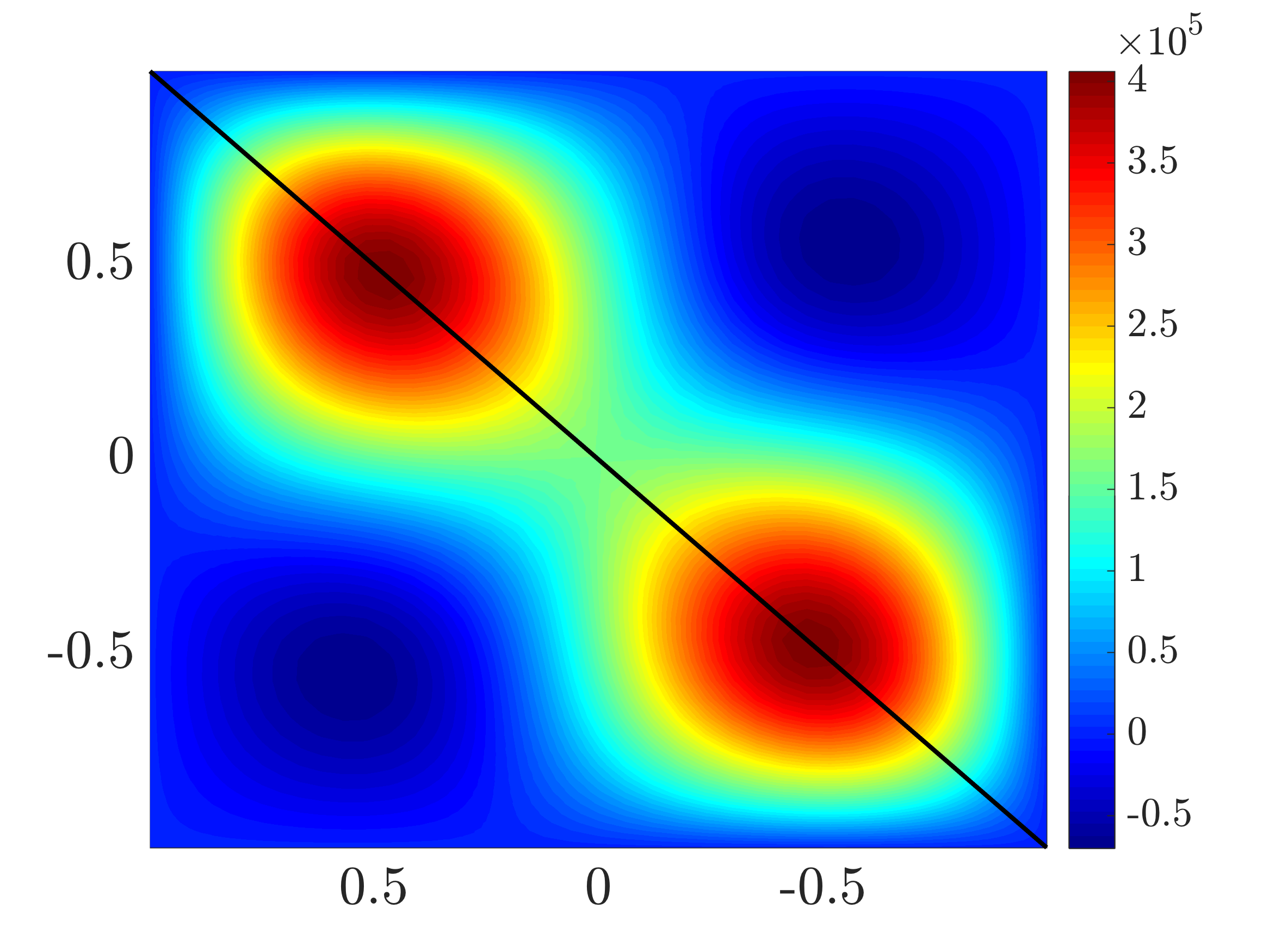}
	\end{tabular}
	&
	&
	\hspace{-1.1cm}
	\begin{tabular}{c}
		\includegraphics[width=4cm]{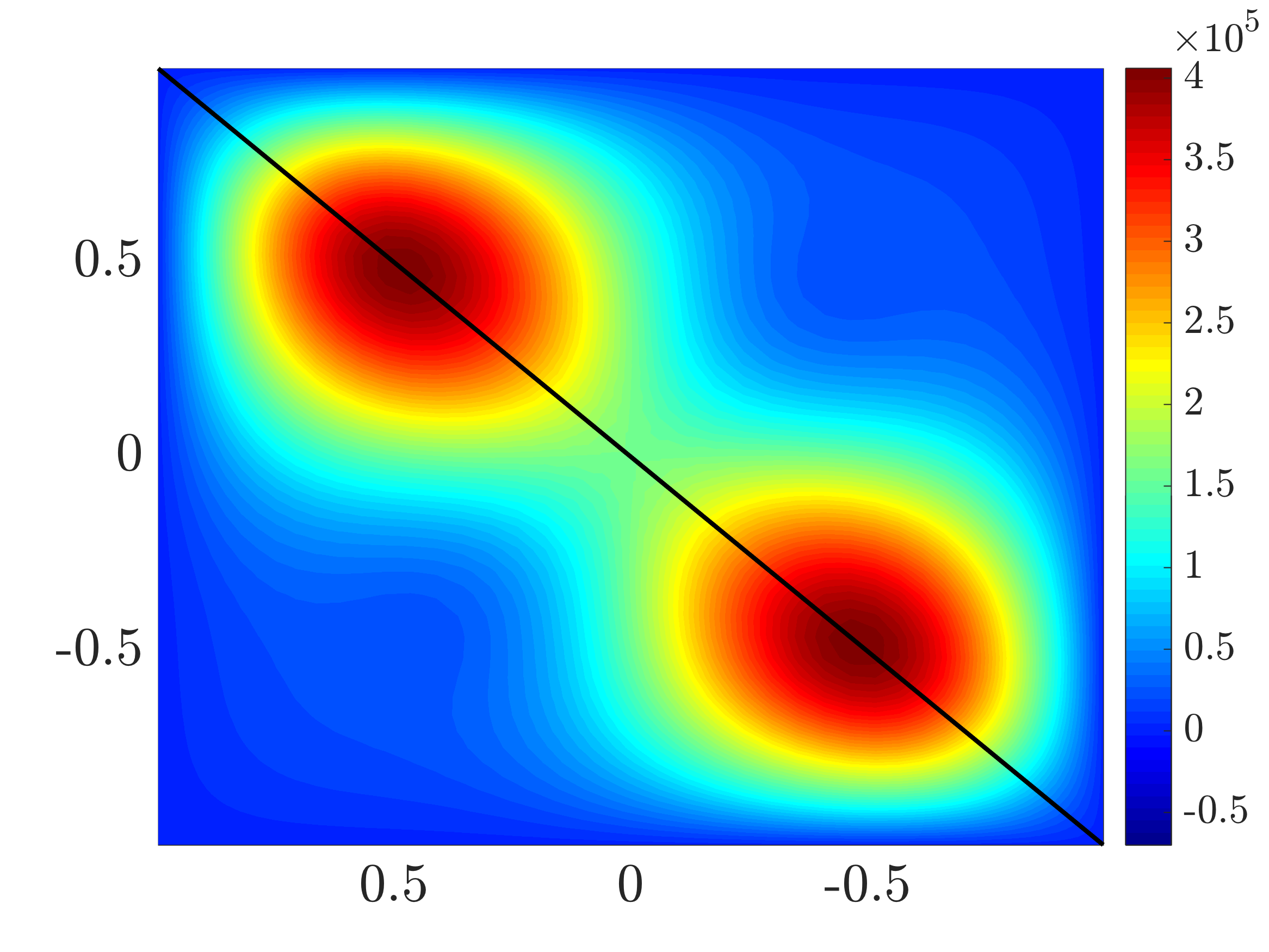}
	\end{tabular}
	\\[-.2cm]
	\hspace{-.8cm} \subfigure[]{\label{fig.Sigvv}}
	& 
	&
	\hspace{-1.4cm} \subfigure[]{\label{fig.Xvv_gammaopt}}
	&
	\\[-.4cm]
	\begin{tabular}{c}
		\vspace{.4cm}
		\hspace{-.6cm}
		\small{$x_2$}
	\end{tabular}
	&
	\hspace{-0.9cm}
	\begin{tabular}{c}
		\includegraphics[width=4cm]{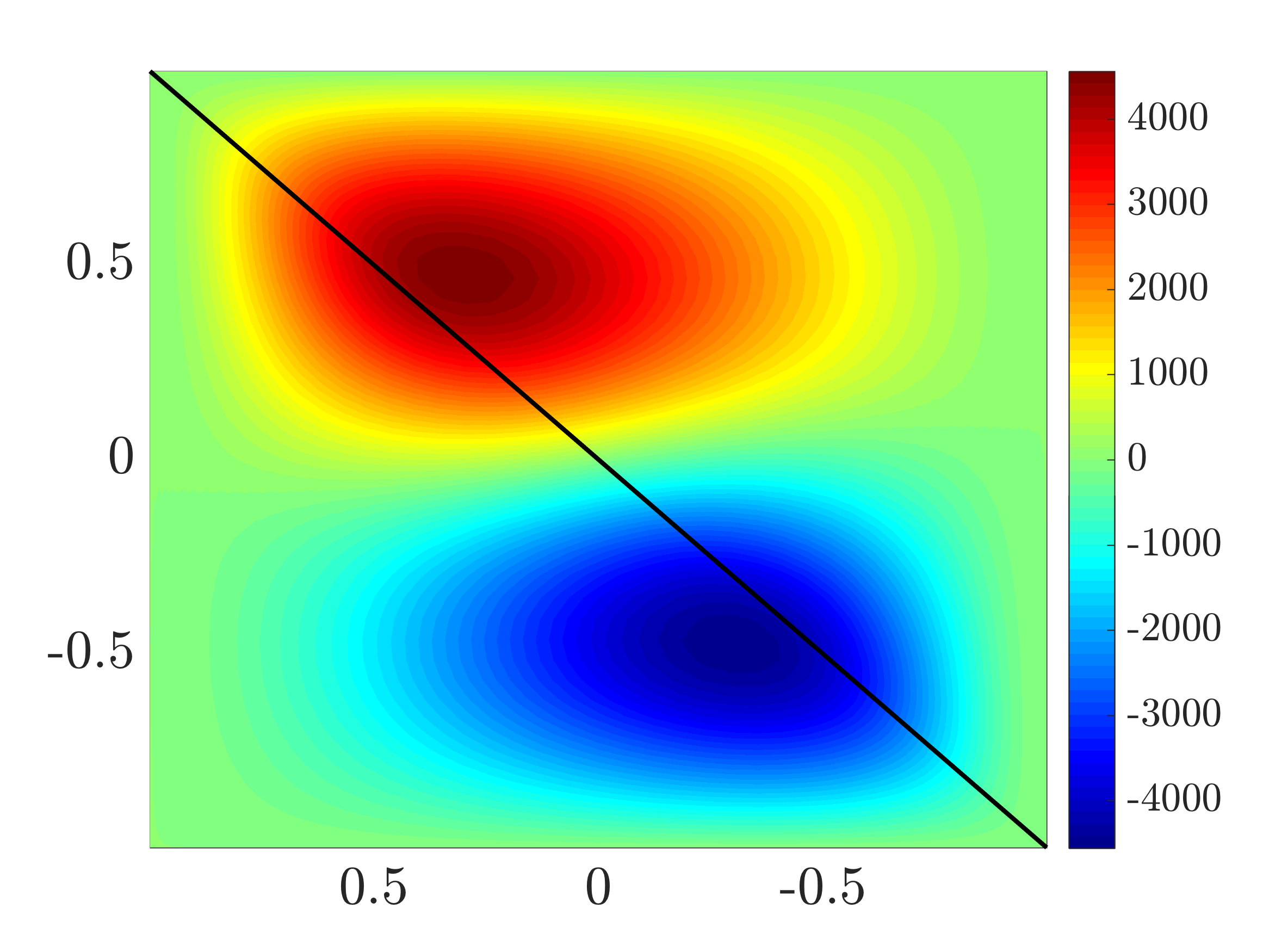}
	\end{tabular}
	&&
	\hspace{-1cm}
	\begin{tabular}{c}
		\includegraphics[width=4cm]{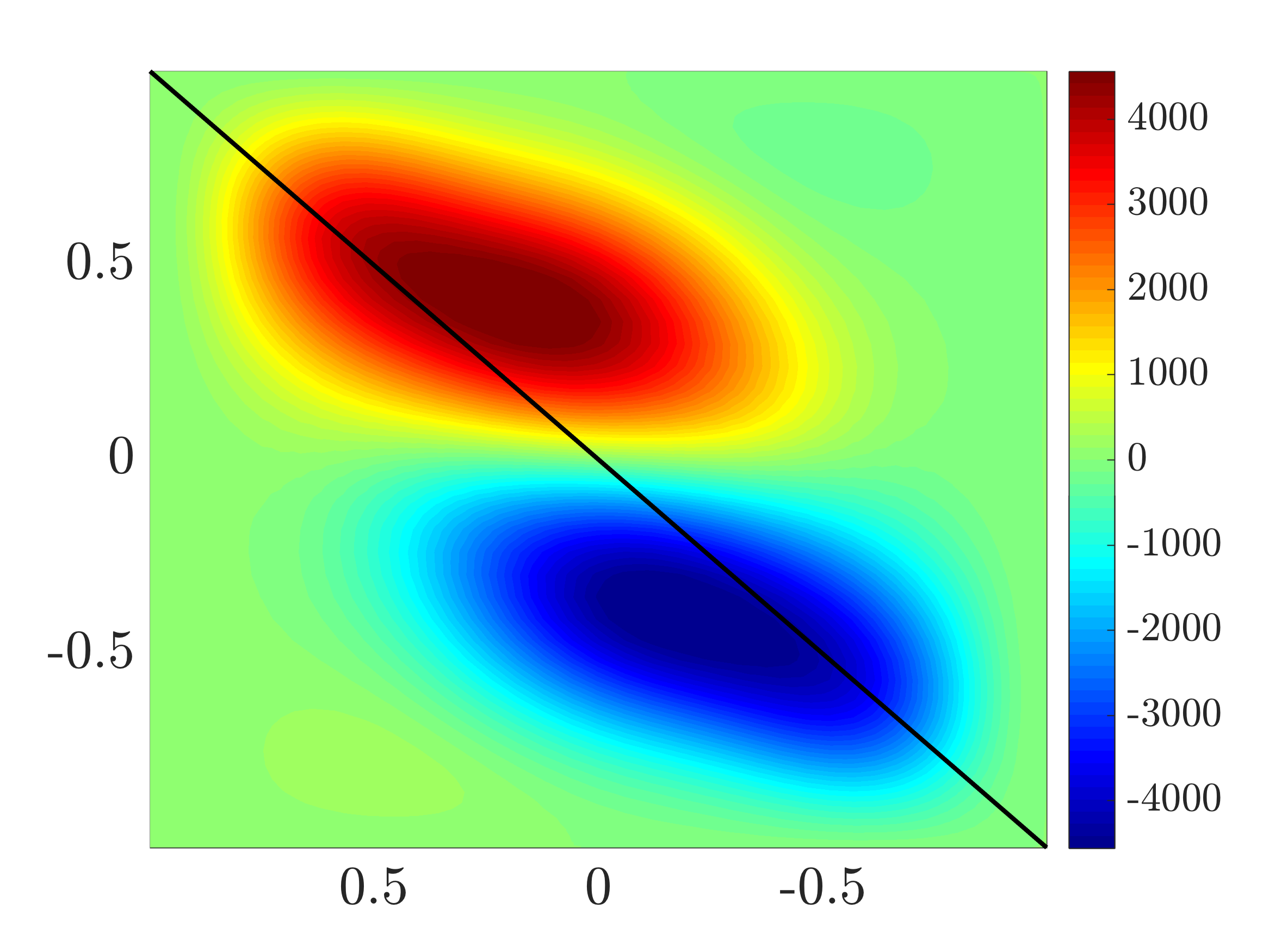}
	\end{tabular}
	\\[-.2cm]
	&
	\hspace{-1cm}
	\small{$x_2$}
	&&
	\hspace{-1cm}
	\small{$x_2$}
	\end{tabular}
		\caption{True covariance matrices of the output velocity field (a, c), and covariance matrices resulting from solving problem~\eqref{eq.MCC-1} (b, d) with $\gamma=100$ and $N=51$. (a, b) Streamwise $\Phi_{11}$, and (c, d) streamwise/normal $\Phi_{12}$ two-point correlation matrices at $\bk = (0,1)$. One-point correlation profiles that are used as problem data are marked along the main diagonals.}
	\label{fig.fluidscompletion}
\end{figure}

Figures~\ref{fig.fluidscompletion}(b,d) show the streamwise, and the streamwise/normal two-point correlation matrices ($\Phi_{11}$ and $\Phi_{12}$ in Fig.~\ref{fig.output-covariance}) resulting from solving~\eqref{eq.MCC-1} with $\gamma=100$. 
Even though only one-point velocity correlations along the main diagonal of these matrices were used in Problem~\ref{prob2}, we observe reasonable recovery of off-diagonal terms of the full two-point velocity correlation matrices and $82\%$ of the original output covariance matrix $\Phi$ is recovered. This quality of completion is consistently observed for various values of $\gamma$ that do not result in the elimination of the critical input channels in the direction of normal velocity, and is an artifact of including the Lyapunov constraint in our formulation. This allows us to simultaneously retain the relevance of the system dynamics and match the partially available statistics of the underlying dynamical system. Additional details regarding the stochastic modeling of turbulent flow statistics and the importance of predicting two-point velocity correlations can be found in~\cite{zarjovgeoJFM17}.

	\vspace*{-1ex}
\section{Concluding remarks}
\label{sec.remarks}

We have examined two problems that arise in modeling and control of stochastically driven dynamical systems. The first addresses the modeling of second-order statistics by a parsimonious perturbation of system dynamics, while the second deals with the optimal selection of sensors/actuators for estimation/control purposes. We have shown that both problems can  be viewed as the selection of suitable feedback gains, guided by similar optimality metrics and subject to closed-loop stability constraints. We cast both problems as optimization problems and use convex surrogates from group-sparsity paradigm to address combinatorial complexity of searching over all possible architectures. While these are SDP representable, the applications that drive our research give rise to the need for scalable algorithms that can handle large problem sizes. We develop a unified algorithmic framework to address both problems using proximal methods. Our algorithms allow handling statistical modeling, as well as sensor and actuator selection, for substantially larger scales than what is amenable to current general-purpose solvers. 

In this work, we promote row sparsity by penalizing a weighted sum of row norms of the feedback gain matrix. While we note that iterative reweighting~\cite{canwakboy08} can improve the row-sparsity patterns determined by this approach, the efficacy of more refined approximations, namely low-rank inducing norms~\cite{gruzarjovranCDC16,gru17}, for which proximal operators can be efficiently computed, is a subject of future research. Moreover, we will investigate solving these problems via primal-dual algorithms based on the proximal augmented Lagrangian~\cite{dhikhojovTAC19,dhikhojovTAC17}, and proximal Newton-type methods~\cite{leesunsau14,stethepat17}.

\section*{Acknowledgments}

We thank Meisam Razaviyayn for useful discussions.

	\vspace*{-3ex}
\appendix

	\vspace*{-2ex}
\subsection{Sensor selection}
\label{sec.sensel}

Consider the LTI system
\begin{align*}
	\ba{rcl}
	\dot{x} 
	& \!\! = \!\! &
	A_s\, x  \; + \; d
	\\[0.1cm]
	y 
	& \!\! = \!\! & 
	C \, x \; + \; \eta
	\ea
\end{align*}
where $y$ denotes measurement data which is corrupted by additive white noise $\eta$. If $(A, C)$ is observable, the observer
\[
	\ba{rcl}
	\dot{\hx} 
	& \!\! = \!\! & 
	A_s \, \hx \;+\; L \, C \left( x \,-\, \hat{x} \right) \;+\; L \, \eta
	\ea
\]
provides an estimate $\hat x$ of the state $x$, where $L$ is the observer gain. When $A_s - LC$ is Hurwitz, the zero-mean estimate of $x$ is given by $\hx$. The Kalman gain minimizes the steady-state variance of $x - \hat{x}$, it is obtained by solving a Riccati equation, and, in general, has no particular structure and uses all available measurements.

Designing a Kalman filter which uses a subset of the available sensors is equivalent to designing a column-sparse Kalman gain matrix $L$. Based on this, the optimal sensor selection problem can be addressed by solving the following regularized optimization problem
 \begin{align}
	\label{eq.senssel}
	\minimize\limits_{L,\, X}
	&
	~~
	\trace \left( X V_d \,+\, L\, V_\eta L^* X \right) \,+\, \gamma\, \ds{\sum^n_{i \, = \, 1} w_i\, \norm{ L\, \mre_i}_2}
	\non
	\\[-.05cm]
	\subject 
	&
	~~
	(A_s - L\,C)^* X \;+\; X (A_s - L\,C) \,+\, C^* C  \;=\; 0
	\non
	\\
	&
	~~
	X \,\succ\, 0
\end{align}
where $\gamma$, $w_i$, $\mre_i$ are as described in Problem~\ref{prob1}, $V_d \succ 0$ is the covariance of $d$, and $V_\eta \succ 0$ is the covariance of $\eta$. By setting the problem data in Problem~\ref{prob1} to
\begin{align*}
	\ba{rclrclrcl}
	A
	&\!\!=\!\!& 
	A_s^*,
	&\quad
	B 
	&\!\!=\!\!&
	C^*,
	&\quad
	Q 
	&\!\!=\!\!& 
	V_d
	\\[.1cm]
	V 
	&\!\!=\!\!& 
	C^* C,
	&\quad
	R 
	&\!\!=\!\!& 
	V_\eta
	&\quad
	&& 
	\ea
\end{align*}
the solution to problem~\eqref{eq.senssel} can be obtained from the solution to the actuator selection problem as $X$ and $L=K^*$. 

\vspace*{-2ex}
\subsection{Non-invertibility of $\cA_1$}
\label{sec.noninvertA}

In cases where the matrix $X$ cannot be expressed via~\eqref{eq.XofY}, since $(A, B)$ is a controllable pair we can center the design variable around a stabilizing controller $K_0$, i.e., by letting $K \DefinedAs K_0 + K_1$, where $K_0$ is held fixed and $K_1$ is the design variable. Based on this, the change of variables introduced in Section~\ref{sec.cov} yields $Y= K_0 X + K_1 X \DefinedAs K_0 X + Y_1$ and $X(Y_1) = \hat{\cA}_1^{-1} (\cB(Y_1) - V)$ with
\begin{align}
	\label{eq.hatA1}
	\hat \cA_1 (X)
	\;\DefinedAs\;
	\left(A - B K_0\right)X \;+\; X\left(A - B K_0\right)^*.
\end{align}
{The resulting optimization problem, 
\begin{align*}
	\ba{cl}
	\minimize\limits_{Y_1}
	&
	f(Y_1) \,+\, \gamma\, g(Y_1 + K_0X(Y_1))
	\\[.25cm]
	\subject 
	&
	\left(1 \,-\, \delta \right) \left[\, \cA_2 (X(Y_1)) \,-\, G\,\right] \;=\; 0
	\\[.15cm]
	&
	X (Y_1) \,\succ\, 0
	 \ea
\end{align*}
involves a nonsmooth term $g$ which is not separable in $Y_1$, and the smooth term is given by}
\begin{multline*}
	f (Y_1)
	\, \DefinedAs \;
	 \trace \, ( Q\,X(Y_1) )  ~ +
	 \\[.1cm]
	 \trace 
	 \, ( 
	 \left(Y_1 + K_0X(Y_1)\right)^* \! R \left(Y_1 + K_0X(Y_1)\right) X^{-1} 
	 ).
	\end{multline*}
Although convex, $g(Y_1 + K_0 X(Y_1))$ does not have an easily computable proximal operator, making it difficult to apply algorithms that are based on proximal methods.

In this case, one may begin with an input matrix $B^0$ such that the pair $(A, B^0)$ is stabilizable and the nonzero columns of $B^0$ correspond to a subset of input channels $\cI$ that {\em always\/} remain active. It would thus be desired to search over input channels from the complement of $\cI$ via the following optimization problem
\begin{align*}
	\ba{cl}
	\minimize\limits_{Y_1}
	&
	f(Y_1) \,+\, \gamma\, \hat g(Y_1)
	\\[.1cm]
	\subject 
	&
	\left(1 \,-\, \delta \right) \left[\,\cA_2 (X(Y_1)) \,-\, G \,\right] \;=\; 0
	\\[.15cm]
	&
	X (Y_1) \,\succ\, 0.
	 \ea
\end{align*}
The operator $\hat{\cA_1}$ in~\eqref{eq.hatA1} would now be defined using $B^0$ and the fixed feedback gain matrix $K_0$ that abides the row-sparsity structure corresponding to $\cI$. The regularization term $\hat{g} (Y_1) \DefinedAs \sum_{i \, \not{\in} \, {\cI}} w_i \| \mre_i^* Y_1 \|_2$ is used to impose row-sparsity on the {\em remaining\/} input channels $i \, {\not \in} \, \cI$ and has an easily computable proximal operator, thus facilitating the use of proximal methods. It is noteworthy that this approach may also be employed to obtain an operator $\hat \cA_1$ which is better conditioned than $\cA_1$.

The alternative approach would be to avoid this problem altogether by not expressing $X$ as a function of $Y$ and directly dualizing the Lyapunov constraint on $X$ and $Y$ via augmented Lagrangian based methods, e.g., ADMM~\cite{dhijovluoCDC14}. However, as we show in Section~\ref{sec.example}, such approaches do not lead to algorithms that are computationally efficient for large problems.

	\vspace*{-2ex}
\subsection{Gradient of $f(Y)$ in~\eqref{eq.PGiter1}}
\label{sec.gradf}

To find $\nabla f(Y)$ in~\eqref{eq.PGiter1}, we expand $f(Y + \eps\, \tilde{Y})$ around $Y$ for the variation $\eps \tilde{Y}$, and collect {first-order terms in $\eps$}. We also account for the variation of $X$ as a result of the variation of $Y$ from
\begin{align*}
	(X + \eps\, \tilde{X})^{-1}
	\;=~
	X^{-1} \,-\; \eps\, X^{-1} \tilde{X}\,X^{-1} \,+\; {o(\eps)}
\end{align*}
and the linear dependence of $\tilde{X}$ on $\tilde{Y}$, i.e., $\tilde{X} = \cA^{-1}_1 (\cB(\tilde{Y}))$. Here, {$o(\eps)$ contains higher-order terms in $\eps$. Thus,} at the $k$th iteration, the gradient of $f$ with respect to $Y$ is given by,
\[
	\nabla f(Y^k)
	~=~
	2\, R\,Y^k X^{-1} 
	\;-\;
	2\, B^* (W_2 \,-\, W_1 )
\]
where $W_1$ and $W_2$ solve the Lyapunov equations
\begin{align*}
	\ba{rcl}
		A^*W_1 \;+\; W_1 A \;+\; X^{-1} Y^{k*} R\, Y^k X^{-1}
		&\!\!=\!\!&
		0
		\\[.15cm]
		A^*W_2 \;+\; W_2 A \;+\; Q
		&\!\!=\!\!&
		0
	\ea
\end{align*}
and $X^{-1}$ denotes the inverse of $X(Y^k)$.

	\vspace*{-1ex}
\subsection{Proofs of Section~\ref{sec.convergence-PG}}
\label{sec.proofs-sec-conv}

\subsubsection{Proof of {Proposition~\ref{prop.megaprop}}}
\label{sec.lemInvarProof} 
Without loss of generality, let $\gamma=1$ {and  $a=b-c$, where $b=f(Y^0)+g(Y^0)$ and $c<g(Y)$ is a lower bound on the function $g$.}
Consider the sublevel set
\begin{align*}
\cE(b) 
\;\DefinedAs\;
\{\, Y \in {\cD} \;|\; f(Y) \,+\, g(Y) \,\le\, b \,\}.
\end{align*}
{ It is easy to verify that $Y^0\in\cE(b)\subset\cD(a)$.}
For a given $Y\in\cE(b)$, let $P$$:\bbR^+\rightarrow \bbC^{m\times n}$ be defined as
\begin{align*}
P(\alpha) 
\;=\; 
\prox_{\alpha g} \! \left(Y \,-\, \alpha\nabla f(Y)\right).
\end{align*}
 In what follows, we show that $P(\alpha)\in\cE(b)$ for all $\alpha\in[0,1/\La]$, with $\La$ being the Lipschitz continuity parameter of $\nabla f(Y)$ over the sublevel set $\mathcal{D}(a)$. {Since $P(0) = Y$, this holds trivially for $\alpha = 0$. For $\alpha>0$, consider the}  quadratic function $l_\alpha:\bbC^{m\times n}\rightarrow \bbR$,
\begin{align*}
        l_\alpha(\hat{Y}) 
        \;\DefinedAs\;
        f(Y) \,+\, \inner{\nabla f(Y)}{\hat{Y} \,-\,Y} \,+\, \dfrac{1}{2\alpha} \, \norm{\hat{Y}\,-\,Y}^2
\end{align*}
{which satisfies}
\begin{align}
        \label{eq.temp6}
        f(\hat{Y})
        \;\le\;
        l_\alpha(\hat{Y})
\end{align}
{for all $\hat{Y}\in \cD(a)$ and $\alpha \in (0,1/\La]$. Inequality~\eqref{eq.temp6}} follows from the $\La$-Lipschitz continuity of $\nabla f(Y)$ over $\cD(a)$ (Descent Lemma). Moreover, by definition,   
\begin{align}
        \label{eq.Palpha}
        P(\alpha) 
        \;=\; 
        \ds{\argmin_{\hat{Y} \, \in \, \bbC^{m\times n}}} \; l_\alpha(\hat{Y}) \,+\, g(\hat{Y})
\end{align}
and $l_\alpha(Y)=f(Y)$, which yields
\begin{align}
        \label{eq.temp5}
        l_\alpha(P(\alpha)) 
        \,+\, 
        g(P(\alpha)) 
        \;\le\; 
        f (Y) 
        \,+\,
        g(Y)
        \;\le\; 
        b
\end{align}
for all {positive $\alpha$. We next show that $P(\alpha) \in \cD(a)$ for all $\alpha \in (0,1/\La]$, which allows us to substitute $P(\alpha)$ for $\hat{Y}$ in~\eqref{eq.temp6} and complete the proof by combining~\eqref{eq.temp6} and~\eqref{eq.temp5}.}

{Since the functions $g$ and $\norm{\cdot}^2$ are coercive, it follows from~\cite[Theorem 26.20]{baucom11} that the map $P(\alpha)$ is continuous. Let $\alpha_1\in(0,+\infty]$ be the smallest scalar such that $ f(P(\alpha_1)) \ge a$. Such $\alpha_1$ exists and $f(P(\alpha_1)) = a$ because the set $\cD$ is open, the function $f(P(\alpha))$ is continuous, and $f(P(0))=f(Y)<a$. We next show that $\alpha_1 > 1/\La$. For the sake of contradiction, suppose $\alpha_1 \le 1/\La$. By substituting $P(\alpha_1)$ for $\hat{Y}$ in Eq.~\eqref{eq.temp6}, using~\eqref{eq.temp5}, and $c<g(P(\alpha_1))$, we arrive at
\begin{align*}
        a \;=\; f(P(\alpha_1)) \;\le\; l_{\alpha_1}(P(\alpha_1)) \;<\; b \, - \, c
\end{align*}
which contradicts with $a=b-c$.} Thus, $\alpha_1 > 1/\La$ {and} $P(\alpha)\in \cD(a)$ for all $\alpha\in[0,1/\La]$. Furthermore, based on this, substituting $P(\alpha)$ in~\eqref{eq.temp6} and utilizing~\eqref{eq.temp5} gives
\begin{align*}
f(P(\alpha)) \,+\, g(P(\alpha))
\;\le\;
b
\end{align*}
which in turn implies $P(\alpha)\in \cE(b)$. 
	
{Based on the fact that we can restrict the domain of the optimization problem~\eqref{eq.general-composite} to the sublevel set $\cD(a)$, the rest of the proof about the convergence rate follows from the proof of~\cite[Theorem 10.29]{bec17}. }

\subsubsection{Proof of {Proposition}~\ref{prop.strongConvexity}}
\label{sec.thmStrongConvProof}
{It is straightforward to verify that the set $\cD_s$ is open.}
We first utilize previously established properties of the set of stabilizing feedback gains to prove that the sublevel sets $\cD(a)$ of the function $f(Y)$ are compact. We then prove that for any convex compact set $\cC \subset \cD_s$ there exist a strong convexity modulus $\mu>0$ and a smoothness parameter $L > 0$ for $f(Y)$ over $\cC$.

Consider the function 
	$
	Y(K) 
	\DefinedAs
	K X(K)
	$
where $K$ belongs to the set of stabilizing feedback gains $\cK_s$ and $X(K)\succ 0$ is the unique solution to the algebraic Lyapunov equation~\eqref{eq:lyapKX}. The function $X(K)$ is continuous and the sublevel sets of the function $f(Y(K))$ 
\begin{align*}
	\mathcal{K}(a) 
	\, \DefinedAs \,
	\{K \in \cK_s \;|\; f(Y(K)) \,\le\, a\}
\end{align*} 
are compact~\cite{toi85}. Since the sublevel set $\cD(a)$ is the image of the compact set $\cK(a)$ under the continuous map $Y(K)$, it follows that $\cD(a)$ is also compact. 
 
The next lemma provides an expression for the second-order approximation of the function $f(Y)$.

	\vsp

\begin{mylem} \label{lem.quadHess}
	The Hessian of the function $f(Y)$ satisfies
	\begin{align*}
            	\inner {\tY}{\nabla^2 f(Y; \tY)}    
            	\, = \,
            	2\,\norm{R^{\tfrac{1}{2}} (\tilde{Y} -Y X^{-1} \cM(\tilde{Y})) X^{-\tfrac{1}{2}}}_F^2
	\end{align*}
	where $X = \cA_1^{-1}(\cB(Y)-V)$ and $\cM(\tilde{Y}) \DefinedAs \cA^{-1}_1(\cB(\tilde{Y}))$.
\end{mylem}

	\vsp
	
\textit{Proof:} For any $Y\in \cD_s$ and $X = \cA_1^{-1}(\cB(Y)-V)$, the function $f(X,Y)$ in Problem~\ref{prob2} reduces to $f(Y)$. The second-order approximation of $f (Y)$ is determined by
	\[
	f (Y + \tY)
	\, \approx \, 
	f (Y)
	\, + \,
	\inner{\nabla f(Y)}{\tY}
	\, + \,  
	\dfrac{1}{2}
	\inner {\tY}{\nabla^2 f(Y; \tY)} 
	\] 
where the matrix $\nabla^2 f(Y; \tY)$ depends linearly on $\tY$. 

The gradient $\nabla f(X,Y)$ can be found by expanding $f(X+\eps\,\tX , Y + \eps\, \tY)$ around the ordered pair $(X,Y)$ for the variation $ (\eps\tX,\eps\tY)$ and collecting {first-order terms in $\eps$.} This yields,
\begin{align*}
	\ba{rcl}
	\nabla_X f(X,Y) 
	& \!\!\! = \!\!\! & 
	Q \, - \, X^{-1}Y^*R\,YX^{-1}
	\\[.1cm]
	\nabla_Y f(X,Y) 
	& \!\!\! = \!\!\! &
	2 \, R \,Y X^{-1}.
	\ea
\end{align*}
To find the Hessian, we expand $\nabla f(X + \eps \, \tX , Y + \eps \, \tY)$,
\begin{align*}
	\ba{lclcl}
	\nabla_X f(X \,+\, \eps \tilde{X},Y) 
	&\!\!\!-\!\!\!& 
	\nabla_X f(X,Y) 
	&\!\!=\!\!&
	\eps N_1 \;+\; {o(\eps)}
	\\[.15cm]
	\nabla_X f(X,Y \,+\, \eps \tilde{Y}) 
	&\!\!\!-\!\!\!&
	\nabla_X f(X,Y) 
	&\!\!=\!\!&
	\eps N_2 \;+\; {o(\eps)}
	\\[.15cm]
	\nabla_Y f(X \,+\, \eps \tilde{X},Y) 
	&\!\!\!-\!\!\!&
	\nabla_Y f(X,Y) 
	&\!\!=\!\!& 
	\eps N_3 \;+\; {o(\eps)}
	\\[.15cm]
	\nabla_Y f(X,Y \,+\, \eps \tilde{Y}) 
	&\!\!\!-\!\!\!&
	\nabla_Y f(X,Y) 
	&\!\!=\!\!& 
	\eps N_4 \;+\; {o(\eps)}
	\ea
\end{align*}
where the matrices
\begin{align*}
	\ba{rcl}
            	N_1
            	& \!\!\!\! \DefinedAs \!\!\!\! &
            	X^{-1}Y^*R \, Y X^{-1}\tilde{X}X^{-1} 
		\,+\, 
		X^{-1}\tilde{X}X^{-1}Y^*R \, YX^{-1}
            	\\[.1cm]
            	 N_2 
            	 & \!\!\!\! \DefinedAs \!\!\!\! &
            	-X^{-1}\tilde{Y}^* R\,Y\,X^{-1}
		\,-\,
		X^{-1}Y^* R \, \tilde{Y} X^{-1}
            	\\[.1cm]
            	 N_3
            	 & \!\!\!\! \DefinedAs \!\!\!\! &
            	-2R\,Y X^{-1}\tilde{X} \,X^{-1}
            	\\[.1cm]
            	 N_4
            	 & \!\!\!\! \DefinedAs \!\!\!\! &
            	2R\,\tilde{Y} X^{-1}
	\ea
\end{align*}
depend linearly on $\tX$ and $\tY$. Thus, we arrive at  
\begin{multline*}
	\inner{(\tilde{X},\tilde{Y})}{\nabla^2 f(X,Y; \tilde{X},\tilde{Y})}
	\\[-0.1cm]
	\,=\,
	{
	\inner{\tilde{X}}{N_1 + N_2}
	\; +\;
	\inner{\tilde{Y}}{N_3 + N_4}
	}
	\\[-0.1cm]
	\,=\;
	2\,\norm{R^{\tfrac{1}{2}} (\tilde{Y} \,-\, YX^{-1}\tilde{X}) X^{-\tfrac{1}{2}}}_F^2.
\end{multline*}
The result follows from $\cA_1(\tilde{X}) = \cB(\tilde{Y})$.
\hfill \qedsymbol

\vsp

Let us define $\zeta$: $\cD_s \times \cS_1 \rightarrow \bbR$ as
\begin{align*}
	\zeta(Y,\tilde{Y})
	\;=\; 
	 \inner {\tY}{\nabla^2f(Y, \tY)}
\end{align*}
where $\cS_1 \DefinedAs \{\tY \in \bbC^{m\times n}  \,|\, \norm{\tY}_F = 1\}$. To establish strong convexity of $f(Y)$ and Lipschitz continuity of its gradient over a compact set $\cC$, we find a positive  lower bound $\mu$ and an upper bound $L$ on $\zeta$,  
	$
	\mu 
	\le  
	\zeta(Y,\tilde{Y}) 
	\le
	L,
	$    
for all $(Y,\tilde{Y}) \in \cC\times \cS_1$.

Using the expression in Lemma~\ref{lem.quadHess}, it is straightforward to show that the function $\zeta$ is continuous. From the continuity of $\zeta(Y,\tilde{Y})$ and the compactness of $\cC\times \cS_1$, it follows that $\zeta$ is bounded on $\cC\times \cS_1$. This implies the existence of an upper bound $L$. To find a positive lower bound, let   $(Y_o,\tilde{Y}_o)$ be a minimizer of the function $\zeta(Y,\tilde{Y})$ over the set $\cC\times \cS_1$. The existence of $(Y_o,\tilde{Y}_o)$ follows from the compactness of \mbox{$\cC\times \cS_1$} and the continuity of the function $\zeta$. We next show  that \mbox{$\mu \,\DefinedAs\, \zeta(Y_o,\tilde{Y}_o) >0$.} 

Suppose, for the sake of contradiction, that $\zeta(Y_o,\tilde{Y}_o)=0$. From Lemma~\ref{lem.quadHess}, we have
\begin{align}
\label{eq.temp1}
	\tilde{Y}_o 
	\;=\; 
	K_o \, \tilde{X}_o   
\end{align}
where $K_o = Y_o X_o^{-1}$, $X_o = X(Y_o)$, and
\begin{align}\label{eq.temp2}
	\tilde{X}_o 
	\;=\;
	\cM(\tilde{Y}_o).
\end{align}
Combining~\eqref{eq.temp2} and the Lyapunov equation{~\eqref{eq.lyap-X-Y}} yields 
 \begin{align}
 \label{eq.temp3}
 	\cA_1(X_o + \tilde{X}_o) \,-\, \cB(Y_o + \tilde{Y}_o) 
	\, = \,
	- V.
 \end{align}
From~\eqref{eq.temp1}, we also have
\begin{align}
\label{eq.temp4}
	Y_o \,+\, \tilde{Y}_o 
	\;=\; 
	K_o \, (X_o  \,+\, \tilde{X}_o).
\end{align} 
Substituting for $Y_o + \tilde{Y}_o$ in~\eqref{eq.temp3} from~\eqref{eq.temp4}, we arrive at
\begin{align*}
	\cA_1(X_o + \tilde{X}_o) \,-\, \cB\left(K_o \, (X_o  +  \tilde{X}_o)\right) 
	\, = \,
	- V.
\end{align*}
Consequently, both $X_o$ and $X_o+\tilde{X}_o$ solve the Lyapunov equation with stabilizing feedback gain $K_o$, which is a contradiction. Thus, $\zeta(Y_o,\tilde{Y}_o)$ is positive. This completes the proof.

	\vspace*{-2ex}
\subsubsection{Proof of Lemma~\ref{lem.L}}
\label{sec.lemL}
We first show that the positive definite matrix $X = \cA_1(\cB(Y)-V)$ satisfies 
\begin{align}
\label{eq.lower-bound-X}
\nu I \, \preceq \, X
\end{align}
with $\nu$ given by~\eqref{eq.nu}. Let $v$ be the normalized eigenvector corresponding to the smallest eigenvalue of $X$. Multiplying {Lyapunov equation~\eqref{eq.lyap-X-Y}} from left and right by $v^*$ and $v$ gives
\begin{align*}
\ba{rcl}
v^* (DX^{\tfrac{1}{2}} \,+\, X^{\tfrac{1}{2}}D^* ) \,v 
& \!\!\! = \!\!\! &
\sqrt{\lambda_{\min}(X)} \; v^*(D \,+\, D^*)\,v 
\\[0.1cm]
& \!\!\! = \!\!\! &
- v^*V\, v
\ea
\end{align*}	
where $D\DefinedAs A X^{1/2} - B Y X^{-1/2}$. We thus have
\begin{align}
\label{eq.lambda-min-X}
\lambda_{\min}(X)
\;=\;
\dfrac{(v^*V\, v)^2}{(v^*(D \,+\, D^*)\,v)^2}
\;\ge\;
\dfrac{\lambda_{\min}^2 (V)}{4\,\norm{D}_2^2}
\end{align}
where we have applied the Cauchy-Schwarz inequality on the denominator. For $Y\in\cD(a)$, we have 
\begin{align*}
\trace \left(Q\, X \,+\, Y^* R\, Y X^{-1} \right)
\;\le\; 
a.
\end{align*}
This inequality along with $\trace \left(Q\, X \right) \geq \lambda_{\min}(Q) \, \norm{X^{1/2}}_F^2$ and $\trace \left(R\, Y X^{-1} Y^* \right) \geq \lambda_{\min} (R) \, \norm{YX^{-1/2}}_F^2$ yields
\begin{subequations}
	\label{eq.bounds-temp}
	\begin{align}
	\norm{X^{1/2}}_F^2 
	\; \le \; 
	a/\lambda_{\min}(Q)
	\\
	\norm{YX^{-1/2}}_F^2 
	\; \le \; 
	a/\lambda_{\min}(R).
	\label{eq.YXinvhalf}
	\end{align}
\end{subequations}
Combination of the triangle inequality, submultiplicative property of the $2$-norm, and~\eqref{eq.bounds-temp} leads to
\begin{align}
\label{eq.bound-D}
\norm{D}_2 
\;\le\; 
\sqrt{a} \left(\dfrac{\sigma_{\max}(A) }{\sqrt{\lambda_{\min}(Q)}} \,+\, \dfrac{\sigma_{\max}(B)}{\sqrt{\lambda_{\min}(R)}} \right).
\end{align}
Inequality~\eqref{eq.lower-bound-X}, with $\nu$ given by~\eqref{eq.nu}, follows from combining~\eqref{eq.lambda-min-X} and~\eqref{eq.bound-D}.

We now show that $L_a$ given by~\eqref{eq.L} is a Lipschitz continuity parameter of $\nabla f$. Form~\eqref{eq.YXinvhalf} and~\eqref{eq.lower-bound-X}, we have
	\begin{align}\label{eq.boundK}
	\norm{YX^{-1}}_F^2
	\;\le\;
	\dfrac{a}{\lambda_{\min}(R)\lambda_{\min}(X)}
	\;\le\;
	\dfrac{a}{\nu\,\lambda_{\min}(R)}.
	\end{align}
This allows us to upper bound the quadratic form provided in Lemma~\ref{lem.quadHess},
\begin{align*}
\inner {\tY}{\nabla^2f(Y, \tY)}    
\;=\;
2\,\norm{R^{\tfrac{1}{2}}(\tilde{Y} \,-\,Y X^{-1} \cM(\tilde{Y}))  X^{-\tfrac{1}{2}}}_F^2.
\end{align*}
In particular, for  $Y\in\cD(a)$ and $\tilde{Y}$ with $\norm{\tilde{Y}}_F=1$, we have
	\begin{multline*}
	2\,\norm{R^{\tfrac{1}{2}}(\tilde{Y} \,-\,Y X^{-1} \cM(\tilde{Y}))  X^{-\tfrac{1}{2}}}_F^2
	\\[.1cm]
	~\le\;
	2\,\lambda_{\max}(R)\lambda_{\max}(X^{-1})\norm{\tY \,-\, YX^{-1}\cM(\tY)}_F^2
	\\[.1cm]
	~\le\;
	2\,\lambda_{\max}(R)\lambda_{\max}(X^{-1})\big(\norm{\tY}_F \,+\,\norm{YX^{-1}\cM(\tY)}_F\big)^2
	\\[.1cm]
	~\le\;
	\dfrac{2\, \lambda_{\max}(R)}{\nu} \,  \Big( 1 \,+\, \dfrac{\sqrt{a} \norm{\cM}_2}{\sqrt{\nu\,\lambda_{\min}(R)}} \Big)^2 
	=\; 
	L_a
	\end{multline*}	
	where the last inequality follows from~\eqref{eq.lower-bound-X},~\eqref{eq.boundK}, and the sub-multiplicative property. This completes the proof.

\vspace*{-4ex}
\subsection{Linear convergence with adaptive step-size selection}
\label{sec.backtracking}

{We show that iterates $\{Y^k\}$ of the PG algorithm with the backtracking scheme of Section~\ref{sec.PG-step-size} remain in $\cD(a)$ and achieve linear convergence. The main challenge in proving the first part of Proposition~\ref{prop.megaprop} is to show that~\eqref{eq.temp6} holds for $\hat{Y}=P(\alpha_k)$, $\alpha_k \in(0,1/\La]$, where $P(\alpha_k)$ is given by~\eqref{eq.Palpha}. However, condition~\eqref{eq.suff-descent} is itself equivalent to~\eqref{eq.temp6} with $\hat{Y}=P(\alpha_k)$. Thus, from the proof of Proposition~\ref{prop.megaprop}, it is easy to verify that the iterates $\{Y^k\}\subset\cD(a)$ and
\begin{align}
\label{eq.linConv-1}
	\norm{Y^{k+1} \, - \, Y^\star}_F^2
	\; \leq \;
	\left( 1 \,-\, \mua \alpha_k \right) \norm{Y^k \, - \, Y^\star}_F^2.
\end{align}
Here, we show that the adaptive backtracking method generates a sequence $\{\alpha_k\}$ that is lower bounded by a fixed positive scalar. Together with~\eqref{eq.linConv-1}, this lower bound yields linear convergence for the PG method with backtracking.}
	
As we discussed in the proof of Proposition~\ref{prop.megaprop}, the step-size $\alpha_k=1/\La$ satisfies conditions~\eqref{eq.backtracking-conds}. Thus, backtracking from a constant initial step-size $\alpha_{k,0}$ would result in a step-size $\alpha_k \geq \min\{\alpha_{k,0},c/L_a\}$, where $c$ is the backtracking parameter in Algorithm~\ref{alg.PG}. {While the initialization $\alpha_{k,0}$ proposed by~\eqref{eq.init-alpha} is not constant, we show that
$\alpha_{k,0} \ge 1/(\sqrt{2}L')$, for any 
\begin{align}
	\label{eq.temp9}
	L' 
	\;\ge\; 
	\norm{{\Delta_2}}_F/\norm{\Delta_1}_F
\end{align}
where $\Delta_1 \,\DefinedAs\, Y^k-Y^{k-1}$ and $\Delta_2 \,\DefinedAs\, \nabla f(Y^k) - \nabla f(Y^{k-1})$.} 
Assuming $ \inner{\Delta_1}{\Delta_2}>0$, the steepest descent and minimum residual step-sizes are given by $\alpha_s = \norm{\Delta_1}^2_F / \inner{\Delta_1}{\Delta_2}$ and $\alpha_m = \inner{\Delta_1}{\Delta_2} / \norm{\Delta_2}_F^2$, respectively. If $ \alpha_m/\alpha_s > 1/2$, then $\sqrt{2}\inner{\Delta_1}{\Delta_2} > \norm{\Delta_1}\norm{\Delta_2}$, which yields
\begin{align*}
	\alpha_{k,0}
	\;=\;
	\dfrac{\inner{\Delta_1}{\Delta_2}}{\norm{\Delta_2}_F^2}
	\;>\;
	\dfrac{\norm{\Delta_1}_F}{\sqrt{2}\norm{\Delta_2}_F}
	\;\ge\; 
	\dfrac{1}{\sqrt{2} \, L'}.
\end{align*}
On the other hand, if $\alpha_m/\alpha_s \le 1/2$, then $\sqrt{2}\inner{\Delta_1}{\Delta_2} \le \norm{\Delta_1}_F\norm{\Delta_2}_F$, which yields
\begin{align*}
	\alpha_{k,0}
	&
	\;=\;
	\dfrac{\norm{\Delta_1}^2_F}{\inner{\Delta_1}{\Delta_2}}
	-
	\dfrac{\inner{\Delta_1}{\Delta_2}}{2\norm{\Delta_2}_F^2}
	\;\ge\;
	\dfrac{3}{2\sqrt{2}} \dfrac{\norm{\Delta_1}_F}{\norm{\Delta_2}_F}
	\;\ge\; 
	\dfrac{3}{2\sqrt{2} \, L'}.
	\end{align*}
		
	Since $Y^k,\, Y^{k-1} \in \cD(a)$, inequality~\eqref{eq.temp9} holds with $L'=L_a$ the Lipschitz continuity factor of $\nabla f(Y)$ over $\cD(a)$. Thus, {the resulting step-size satisfies $\alpha_k \geq \min\{1/(\sqrt{2}L_a),c/L_a\}$.}

	\vspace*{-2ex}
\subsection{Gradient of $F(Y)$ in~\eqref{eq.PGiter2}}
\label{sec.gradF}

Similar to Appendix~\ref{sec.gradf}, we expand $F(Y + \eps\, \tilde{Y})$ around $Y$ for the variation $\eps \tilde{Y}$, and collect {first-order terms in $\eps$}. At the $k$th iteration, the gradient of $F$ with respect to $Y$ is given by,
\begin{align*}
	\nabla F(Y^k)
	~=~
	2\,Y^k X^{-1} 
	\;-\;
	2\, B^* (W_2 \,+\, \rho_k\,W_3 \,-\, W_1 ),
\end{align*}
where $W_1$, $W_2$, and $W_3$ solve the Lyapunov equations
\begin{align*}
	\ba{rcl}
		A^*W_1 \;+\; W_1 A \;+\; X^{-1} Y^{k*} Y^k X^{-1}
		&\!\!=\!\!&
		0
		\\[.15cm]
		A^*W_2 \;+\; W_2 A \;+\; \cA_2^\dagger \left(\Lambda^k \right)
		&\!\!=\!\!&
		0
		\\[.15cm]
		A^*W_3 \;+\; W_3 A \;+\; \cA_2^\dagger \left(\cA_2 \left(X(Y^k)\right) \,-\, G \right)
		&\!\!=\!\!&
		0
	\ea
\end{align*}
Here, $X^{-1}$ denotes the inverse of $X(Y^k)$ and the adjoint of the operator $\cA_2$ is given by
	$
	\cA_2^\dagger \left(\Lambda \right)
	\DefinedAs 
	C^* \left(E \circ \Lambda \right) C.
	$

	\vspace*{-1.4cm}
\begin{biography}
[{\includegraphics[width=1in,height=1.25in,clip,keepaspectratio]{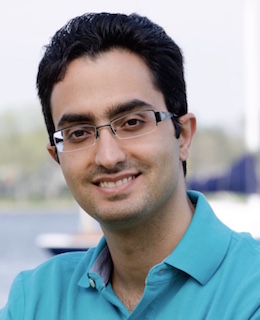}}]
{Armin Zare} (S'08--M'17) received the BSc in Electrical Engineering from Sharif University of Technology, Tehran, Iran, in 2010 and the MSEE and PhD degree in Electrical Engineering from the University of Minnesota, Minneapolis, MN, USA, in 2016. He is an Assistant Professor of Mechanical Engineering at the University of Texas at Dallas, Richardson, TX, USA. He was a Post-doctoral Research Associate in the Ming Hsieh Department of Electrical and Computer Engineering at the University of Southern California, Los Angeles, from 2017 to 2019. His primary research interests are in the modeling and control of complex fluid flows using tools from optimization and systems theory. He was a recipient of the Doctoral Dissertation Fellowship at the University of Minnesota in 2015 and a finalist for the Best Student Paper Award at the American Control Conference in 2014.
\end{biography}

	\vspace*{-1.4cm}

\begin{IEEEbiography}
[{\includegraphics[width=1in,height=1.25in,clip,keepaspectratio]{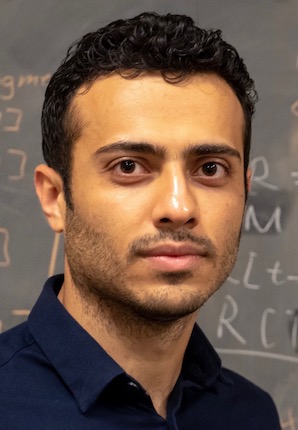}}]
{Hesameddin Mohammadi} (S'17) received the MSc degree from Arizona State University, Tempe, AZ, USA, in 2017 and the BSc degree from Sharif University of Technology, Tehran, Iran, in 2015, both in Mechanical Engineering. He is currently pursuing the PhD degree in the Ming Hsieh Department of Electrical and Computer Engineering at the University of Southern California, Los Angeles, CA, USA. His primary research interests are in large-scale optimization, control, and inference problems.
\end{IEEEbiography}

	\vspace*{-1.4cm}

\begin{IEEEbiography}
[{\includegraphics[width=1in,height=1.25in,clip,keepaspectratio]{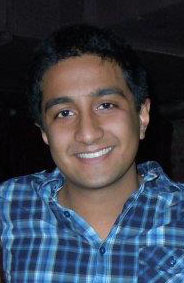}}]
{Neil K.\ Dhingra} (S'10--M'19) received his PhD in Electrical and Computer Engineering from the University of Minnesota, Twin Cities in 2017 where he developed tools for regularized optimization and studied the design of structured controllers for complex large-scale systems. He received his BSE in Electrical Engineering from the University of Michigan, Ann Arbor, in 2010 and has worked with NASA JPL, NASA AFRC, and the WIMS Center. Dr.\ Dhingra is currently a Research Scientist at Numerica Corporation in Fort Collins, CO.
\end{IEEEbiography}

	\vspace*{-1.4cm}
	
	\begin{biography}
[{\includegraphics[width=1in,height=1.25in,clip,keepaspectratio]{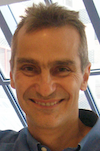}}]
{Tryphon T. Georgiou} (M'79--SM'99--F'00) received the Diploma in Mechanical and Electrical Engineering from the National Technical University of Athens, Greece, in 1979 and the PhD degree from the University of Florida, Gainesville, FL, USA, in 1983.  He is currently a Chancellor's Professor at the Department of Mechanical and Aerospace Engineering, University of California, Irvine, CA, USA. Earlier, he served on the faculty of Florida Atlantic University from 1983 to 1986, Iowa State University from 1986 to 1989, and the University of Minnesota from 1989 to 2016.

He is a recipient of the George S.\ Axelby Outstanding Paper award of the IEEE Control Systems Society for the years 1992, 1999, 2003, and 2017, a Fellow of the Institute of Electrical and Electronic Engineers (IEEE) and the International Federation of Automatic Control (IFAC), and a Foreign Member of the Royal Swedish Academy of Engineering Sciences (IVA).
    \end{biography}

	\vspace*{-1.4cm}
	
    \begin{IEEEbiography}
[{\includegraphics[width=1.in,height=1.25in,clip,keepaspectratio]{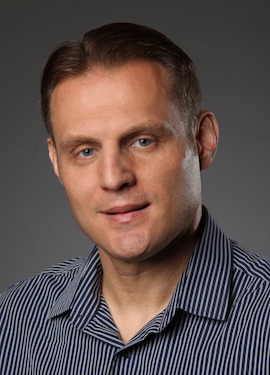}}]
{Mihailo R. Jovanovi\'c} (S'00--M'05--SM'13--F'19) received the PhD degree from the University of California at Santa Barbara in 2004. He is a Professor in the Ming Hsieh Department of Electrical and Computer Engineering and the Founding Director of the Center for Systems and Control at the University of Southern California, Los Angeles, CA. He was a faculty member in the Department of Electrical and Computer Engineering at the University of Minnesota, Twin Cities, MN, from 2004 until 2017, and has held visiting positions with Stanford University and the Institute for Mathematics and its Applications. 

Prof.\ Jovanovi\'c is a Fellow of the American Physical Society (APS) and the Institute of Electrical and Electronic Engineers (IEEE). He received a CAREER Award from the National Science Foundation in 2007, the George S. Axelby Outstanding Paper Award from the IEEE Control Systems Society in 2013, and the Distinguished Alumni Award from the Department of Mechanical Engineering at UC Santa Barbara in 2014. Papers of his students were finalists for the Best Student Paper Award at the American Control Conference in 2007 and 2014.
	\end{IEEEbiography}


\begin{thebibliography}{10}
\providecommand{\url}[1]{#1}
\csname url@rmstyle\endcsname
\providecommand{\newblock}{\relax}
\providecommand{\bibinfo}[2]{#2}
\providecommand\BIBentrySTDinterwordspacing{\spaceskip=0pt\relax}
\providecommand\BIBentryALTinterwordstretchfactor{4}
\providecommand\BIBentryALTinterwordspacing{\spaceskip=\fontdimen2\font plus
\BIBentryALTinterwordstretchfactor\fontdimen3\font minus
  \fontdimen4\font\relax}
\providecommand\BIBforeignlanguage[2]{{%
\expandafter\ifx\csname l@#1\endcsname\relax
\typeout{** WARNING: IEEEtran.bst: No hyphenation pattern has been}%
\typeout{** loaded for the language `#1'. Using the pattern for}%
\typeout{** the default language instead.}%
\else
\language=\csname l@#1\endcsname
\fi
#2}}

\bibitem{boyelgferbal94}
S.~Boyd, L.~E. Ghaoui, E.~Feron, and V.~Balakrishnan, \emph{Linear matrix
  inequalities in system and control theory}.\hskip 1em plus 0.5em minus
  0.4em\relax SIAM, 1994.

\bibitem{dulpag00}
G.~E. Dullerud and F.~Paganini, \emph{{A course in robust control theory: a
  convex approach}}.\hskip 1em plus 0.5em minus 0.4em\relax New York:
  Springer-Verlag, 2000.

\bibitem{fazhinboy01}
M.~Fazel, H.~Hindi, and S.~Boyd, ``A rank minimization heuristic with
  application to minimum order system approximation,'' in \emph{Proceedings of
  the 2001 American Control Conference}, 2001, pp. 4734--4739.

\bibitem{boyvan04}
S.~Boyd and L.~Vandenberghe, \emph{Convex optimization}.\hskip 1em plus 0.5em
  minus 0.4em\relax Cambridge University Press, 2004.

\bibitem{fazhinboy04}
M.~Fazel, H.~Hindi, and S.~Boyd, ``Rank minimization and applications in system
  theory,'' in \emph{Proceedings of the 2004 American Control Conference},
  2004, pp. 3273--3278.

\bibitem{liuvan09}
Z.~Liu and L.~Vandenberghe, ``Interior-point method for nuclear norm
  approximation with application to system identification,'' \emph{SIAM J.
  Matrix Anal. Appl.}, vol.~31, no.~3, pp. 1235--1256, 2009.

\bibitem{jovdhiEJC16}
M.~R. Jovanovi\'c and N.~K. Dhingra, ``Controller architectures: tradeoffs
  between performance and structure,'' \emph{Eur. J. Control}, vol.~30, pp.
  76--91, July 2016.

\bibitem{zarchejovgeoTAC17}
A.~Zare, Y.~Chen, M.~R. Jovanovi\'c, and T.~T. Georgiou, ``Low-complexity
  modeling of partially available second-order statistics: theory and an
  efficient matrix completion algorithm,'' \emph{IEEE Trans. Automat. Control},
  vol.~62, no.~3, pp. 1368--1383, March 2017.

\bibitem{zarjovgeoJFM17}
A.~Zare, M.~R. Jovanovi\'c, and T.~T. Georgiou, ``Colour of turbulence,''
  \emph{J. Fluid Mech.}, vol. 812, pp. 636--680, February 2017.

\bibitem{zargeojovARC20}
A.~Zare, T.~T. Georgiou, and M.~R. Jovanovi\'c, ``Stochastic dynamical modeling
  of turbulent flows,'' \emph{Annu. Rev. Control Robot. Auton. Syst.}, vol.~3,
  May 2020, in press; doi:10.1146/annurev-control-053018-023843; also
  arXiv:1908.09487.

\bibitem{farioa93}
B.~F. Farrell and P.~J. Ioannou, ``Stochastic forcing of the linearized
  {N}avier-{S}tokes equations,'' \emph{Phys.\ Fluids A}, vol.~5, no.~11, pp.
  2600--2609, 1993.

\bibitem{bamdah01}
B.~Bamieh and M.~Dahleh, ``Energy amplification in channel flows with
  stochastic excitation,'' \emph{Phys.\ Fluids}, vol.~13, no.~11, pp.
  3258--3269, 2001.

\bibitem{jovbamJFM05}
M.~R. Jovanovi\'c and B.~Bamieh, ``Componentwise energy amplification in
  channel flows,'' \emph{J. Fluid Mech.}, vol. 534, pp. 145--183, July 2005.

\bibitem{moajovJFM12}
R.~Moarref and M.~R. Jovanovi\'c, ``Model-based design of transverse wall
  oscillations for turbulent drag reduction,'' \emph{J. Fluid Mech.}, vol. 707,
  pp. 205--240, September 2012.

\bibitem{ranzarhacjovPRF19b}
W.~Ran, A.~Zare, M.~J.~P. Hack, and M.~R. Jovanovi\'c, ``Stochastic receptivity
  analysis of boundary layer flow,'' \emph{Phys. Rev. Fluids}, vol.~4, no.~9,
  p. 093901 (28 pages), September 2019.

\bibitem{jovbamCDC01}
M.~R. Jovanovi\'c and B.~Bamieh, ``Modelling flow statistics using the
  linearized {N}avier-{S}tokes equations,'' in \emph{Proceedings of the 40th
  IEEE Conference on Decision and Control}, 2001, pp. 4944--4949.

\bibitem{zarjovgeoCDC16}
A.~Zare, M.~R. Jovanovi\'c, and T.~T. Georgiou, ``Perturbation of system
  dynamics and the covariance completion problem,'' in \emph{Proceedings of the
  55th IEEE Conference on Decision and Control}, 2016, pp. 7036--7041.

\bibitem{HotSke87}
A.~Hotz and R.~E. Skelton, ``Covariance control theory,'' \emph{Int. J.
  Control}, vol.~46, no.~1, pp. 13--32, 1987.

\bibitem{yasskegri93}
K.~Yasuda, R.~E. Skelton, and K.~M. Grigoriadis, ``Covariance controllers: {A}
  new parametrization of the class of all stabilizing controllers,''
  \emph{Automatica}, vol.~29, no.~3, pp. 785--788, 1993.

\bibitem{grikarske94}
K.~M. Grigoriadis and R.~E. Skelton, ``Alternating convex projection methods
  for covariance control design,'' \emph{Int. J. Control}, vol.~60, no.~6, pp.
  1083--1106, 1994.

\bibitem{chegeopav16b}
Y.~Chen, T.~T. Georgiou, and M.~Pavon, ``Optimal steering of a linear
  stochastic system to a final probability distribution, {P}art {II},''
  \emph{IEEE Trans. Automat. Control}, vol.~61, no.~5, pp. 1170--1180, 2016.

\bibitem{linjovTAC09}
F.~Lin and M.~R. Jovanovi\'c, ``Least-squares approximation of structured
  covariances,'' \emph{IEEE Trans. Automat. Control}, vol.~54, no.~7, pp.
  1643--1648, July 2009.

\bibitem{zorfer12}
M.~Zorzi and A.~Ferrante, ``On the estimation of structured covariance
  matrices,'' \emph{Automatica}, vol.~48, no.~9, pp. 2145--2151, 2012.

\bibitem{sumcorlyg16}
T.~H. Summers, F.~L. Cortesi, and J.~Lygeros, ``On submodularity and
  controllability in complex dynamical networks,'' \emph{IEEE Trans. Control
  Netw. Syst.}, vol.~3, no.~1, pp. 91--101, 2016.

\bibitem{tzorahpapjad16}
V.~Tzoumas, M.~A. Rahimian, G.~J. Pappas, and A.~Jadbabaie, ``Minimal actuator
  placement with bounds on control effort,'' \emph{IEEE Trans. Control Netw.
  Syst.}, vol.~3, no.~1, pp. 67--78, 2016.

\bibitem{zhaayosun17}
H.~Zhang, R.~Ayoub, and S.~Sundaram, ``Sensor selection for kalman filtering of
  linear dynamical systems: Complexity, limitations and greedy algorithms,''
  \emph{Automatica}, vol.~78, pp. 202--210, 2017.

\bibitem{ols18}
A.~Olshevsky, ``On (non)supermodularity of average control energy,'' \emph{IEEE
  Trans. Control Netw. Syst.}, vol.~5, no.~3, pp. 1177--1181, 2018.

\bibitem{josboy09}
S.~Joshi and S.~Boyd, ``Sensor selection via convex optimization,'' \emph{IEEE
  Trans. Signal Process.}, vol.~57, no.~2, pp. 451--462, 2009.

\bibitem{liuchefarmasleuvar16}
S.~Liu, S.~P. Chepuri, M.~Fardad, E.~Ma{\c{s}}azade, G.~Leus, and P.~K.
  Varshney, ``Sensor selection for estimation with correlated measurement
  noise,'' \emph{IEEE Trans. Signal Process.}, vol.~64, no.~13, pp. 3509--3522,
  2016.

\bibitem{kekgiawol12}
V.~Kekatos, G.~B. Giannakis, and B.~Wollenberg, ``Optimal placement of phasor
  measurement units via convex relaxation,'' \emph{IEEE Trans. Power Syst.},
  vol.~27, no.~3, pp. 1521--1530, 2012.

\bibitem{rog00}
J.~L. Rogers, ``A parallel approach to optimum actuator selection with a
  genetic algorithm,'' in \emph{AIAA Guidance, Navigation, and Control
  Conference}, 2000, pp. 14--17.

\bibitem{konyatino90}
S.~Kondoh, C.~Yatomi, and K.~Inoue, ``The positioning of sensors and actuators
  in the vibration control of flexible systems,'' \emph{JSME Int. J., Ser.
  III}, vol.~33, no.~2, pp. 145--152, 1990.

\bibitem{hirdokobi00}
K.~Hiramoto, H.~Doki, and G.~Obinata, ``Optimal sensor/actuator placement for
  active vibration control using explicit solution of algebraic {R}iccati
  equation,'' \emph{J.\ Sound Vib.}, vol. 229, no.~5, pp. 1057--1075, 2000.

\bibitem{cherowJFM11}
K.~K. Chen and C.~W. Rowley, ``{$\Ht$} optimal actuator and sensor placement in
  the linearised complex {G}inzburg-{L}andau system,'' \emph{J.\ Fluid Mech.},
  vol. 681, pp. 241--260, 2011.

\bibitem{farlinjovACC11}
M.~Fardad, F.~Lin, and M.~R. Jovanovi\'c, ``Sparsity-promoting optimal control
  for a class of distributed systems,'' in \emph{Proceedings of the 2011
  American Control Conference}, 2011, pp. 2050--2055.

\bibitem{linfarjovACC12}
F.~Lin, M.~Fardad, and M.~R. Jovanovi\'c, ``Sparse feedback synthesis via the
  alternating direction method of multipliers,'' in \emph{Proceedings of the
  2012 American Control Conference}, 2012, pp. 4765--4770.

\bibitem{linfarjovTAC13admm}
F.~Lin, M.~Fardad, and M.~R. Jovanovi\'c, ``Design of optimal sparse feedback
  gains via the alternating direction method of multipliers,'' \emph{IEEE
  Trans. Automat. Control}, vol.~58, no.~9, pp. 2426--2431, September 2013.

\bibitem{masfarvar12}
E.~Masazade, M.~Fardad, and P.~K. Varshney, ``Sparsity-promoting extended
  {K}alman filtering for target tracking in wireless sensor networks,''
  \emph{IEEE Signal Process.\ Lett.}, vol.~19, pp. 845--848, 2012.

\bibitem{liufarmasvar14}
S.~Liu, M.~Fardad, E.~Masazade, and P.~K. Varshney, ``Optimal periodic sensor
  scheduling in networks of dynamical systems,'' \emph{IEEE Trans. Signal
  Process.}, vol.~62, no.~12, pp. 3055--3068, 2014.

\bibitem{polkhlshc13}
B.~Polyak, M.~Khlebnikov, and P.~Shcherbakov, ``An {LMI} approach to structured
  sparse feedback design in linear control systems,'' in \emph{Proceedings of
  the 2013 European Control Conference}, 2013, pp. 833--838.

\bibitem{munpfiwol14}
U.~M{\"u}nz, M.~Pfister, and P.~Wolfrum, ``Sensor and actuator placement for
  linear systems based on {$\cH_2$} and {$\cH_\infty$} optimization,''
  \emph{IEEE Trans. Automat. Control}, vol.~59, no.~11, pp. 2984--2989, 2014.

\bibitem{yualin06}
M.~Yuan and Y.~Lin, ``Model selection and estimation in regression with grouped
  variables,'' \emph{J. R. Stat. Soc. Series B Stat. Methodol.}, vol.~68,
  no.~1, pp. 49--67, 2006.

\bibitem{bampagdah02}
B.~Bamieh, F.~Paganini, and M.~A. Dahleh, ``Distributed control of spatially
  invariant systems,'' \emph{IEEE Transactions on Automatic Control}, vol.~47,
  no.~7, pp. 1091--1107, July 2002.

\bibitem{horjoh12}
R.~A. Horn and C.~R. Johnson, \emph{Matrix Analysis}.\hskip 1em plus 0.5em
  minus 0.4em\relax Cambridge University Press, 2012.

\bibitem{becteb09}
A.~Beck and M.~Teboulle, ``A fast iterative shrinkage-thresholding algorithm
  for linear inverse problems,'' \emph{SIAM {J.} {I}maging {S}ci.}, vol.~2,
  no.~1, pp. 183--202, 2009.

\bibitem{parboy13}
N.~Parikh and S.~Boyd, ``Proximal algorithms,'' \emph{Found. Trends Optim.},
  vol.~1, no.~3, pp. 123--231, 2013.

\bibitem{golstubar14}
T.~Goldstein, C.~Studer, and R.~Baraniuk, ``A field guide to forward-backward
  splitting with a {FASTA} implementation,'' \emph{{\em arXiv:1411.3406}},
  2014.

\bibitem{baucom11}
H.~H. Bauschke and P.~L. Combettes, \emph{Convex analysis and monotone operator
  theory in {H}ilbert spaces}.\hskip 1em plus 0.5em minus 0.4em\relax Springer,
  2011, vol. 408.

\bibitem{bec17}
A.~Beck, \emph{{First-Order Methods in Optimization}}.\hskip 1em plus 0.5em
  minus 0.4em\relax SIAM, 2017, vol.~25.

\bibitem{mohzarsoljovCDC19}
H.~Mohammadi, A.~Zare, M.~Soltanolkotabi, and M.~R. Jovanovi\'c, ``Global
  exponential convergence of gradient methods over the nonconvex landscape of
  the linear quadratic regulator,'' in \emph{Proceedings of the 58th IEEE
  Conference on Decision and Control}, Nice, France, 2019, pp. 7474--7479.

\bibitem{ber82}
D.~P. Bertsekas, \emph{Constrained optimization and {L}agrange multiplier
  methods}.\hskip 1em plus 0.5em minus 0.4em\relax New York: Academic Press,
  1982.

\bibitem{ber99}
D.~P. Bertsekas, \emph{Nonlinear programming}.\hskip 1em plus 0.5em minus
  0.4em\relax Belmont, MA: Athena Scientific, 1999.

\bibitem{nocwri06}
J.~Nocedal and S.~J. Wright, \emph{Numerical Optimization}.\hskip 1em plus
  0.5em minus 0.4em\relax Springer, 2006.

\bibitem{boyparchupeleck11}
S.~Boyd, N.~Parikh, E.~Chu, B.~Peleato, and J.~Eckstein, ``Distributed
  optimization and statistical learning via the alternating direction method of
  multipliers,'' \emph{Found. Trends Mach. Learn.}, vol.~3, no.~1, pp. 1--122,
  2011.

\bibitem{dhijovluoCDC14}
N.~K. Dhingra, M.~R. Jovanovi\'c, and Z.~Q. Luo, ``An {ADMM} algorithm for
  optimal sensor and actuator selection,'' in \emph{Proceedings of the 53rd
  IEEE Conference on Decision and Control}, 2014, pp. 4039--4044.

\bibitem{canwakboy08}
E.~J. Candes, M.~B. Wakin, and S.~P. Boyd, ``Enhancing sparsity by reweighted
  $\ell_1$ minimization,'' \emph{J. Fourier Anal. Appl.}, vol.~14, no. 5-6, pp.
  877--905, 2008.

\bibitem{crohoh93}
M.~C. Cross and P.~C. Hohenberg, ``Pattern formation outside of equilibrium,''
  \emph{Rev. Mod. Phys.}, vol.~65, no.~3, p. 851, 1993.

\bibitem{burkno06}
J.~Burke and E.~Knobloch, ``Localized states in the generalized
  {S}wift-{H}ohenberg equation,'' \emph{Phys. Rev. E}, vol.~73, no.~5, p.
  056211, 2006.

\bibitem{SDPT3}
K.-C. Toh, M.~J. Todd, and R.~H. T{\"u}t{\"u}nc{\"u}, ``{SDPT3-a MATLAB
  software package for semidefinite programming, version 1.3},'' \emph{Optim.
  Methods Softw.}, vol.~11, no. 1-4, pp. 545--581, 1999.

\bibitem{cvx}
M.~Grant and S.~Boyd, ``{CVX}: Matlab software for disciplined convex
  programming, version 2.1,'' \url{http://cvxr.com/cvx}, Mar. 2014.

\bibitem{sum16}
T.~Summers, ``Actuator placement in networks using optimal control performance
  metrics,'' in \emph{Proceedings of the 55th IEEE Conference on Decision and
  Control}, 2016, pp. 2703--2708.

\bibitem{gruzarjovranCDC16}
C.~Grussler, A.~Zare, M.~R. Jovanovi\'c, and A.~Rantzer, ``The use of the $r*$
  heuristic in covariance completion problems,'' in \emph{Proceedings of the
  55th IEEE Conference on Decision and Control}, 2016, pp. 1978--1983.

\bibitem{gru17}
C.~Grussler, ``Rank reduction with convex constraints,'' Ph.D. dissertation,
  Lund University, 2017.

\bibitem{dhikhojovTAC19}
N.~K. Dhingra, S.~Z. Khong, and M.~R. Jovanovi\'c, ``The proximal augmented
  {L}agrangian method for nonsmooth composite optimization,'' \emph{IEEE Trans.
  Automat. Control}, vol.~64, no.~7, pp. 2861--2868, July 2019.

\bibitem{dhikhojovTAC17}
N.~K. Dhingra, S.~Z. Khong, and M.~R. Jovanovi\'c, ``A second order primal-dual
  method for nonsmooth convex composite optimization,'' \emph{IEEE Trans.
  Automat. Control}, 2017, submitted; also arXiv:1709.01610.

\bibitem{leesunsau14}
J.~D. Lee, Y.~Sun, and M.~A. Saunders, ``Proximal {N}ewton-type methods for
  minimizing composite functions,'' \emph{SIAM J. Optim.}, vol.~24, no.~3, pp.
  1420--1443, 2014.

\bibitem{stethepat17}
L.~Stella, A.~Themelis, and P.~Patrinos, ``Forward-backward quasi-{N}ewton
  methods for nonsmooth optimization problems,'' \emph{Comput. Optim. Appl.},
  vol.~67, no.~3, pp. 443--487, 2017.

\bibitem{toi85}
H.~T. Toivonen, ``A globally convergent algorithm for the optimal constant
  output feedback problem,'' \emph{Int. J. Control}, vol.~41, no.~6, pp.
  1589--1599, 1985.

\end{thebibliography}


\vspace*{-1ex}
\end{document}